%% file: main.tex
\documentclass[12pt]{report}

\usepackage[dvips]{hyperref}
\usepackage[numbers]{natbib}
\usepackage{marcmacros}
\usepackage{suthesis-2e}

\begin{document}
\include{header}
\include{abstract}
\include{preface}
\include{intro}
\include{lit}
\include{implement}
\include{examples}
\include{proof}
\include{discuss}
\include{voronoi}
\include{afterword}



\bibliographystyle{ims}
\bibliography{thesis,newbib,bcartbib}

\end{document}

%% file: header.tex
\title{Nonparametric Bayesian Classification}
\author{Marc A. Coram}
\dept{Statistics}
\principaladviser{Persi Diaconis}
\firstreader{Jerome H. Friedman}
\secondreader{Bradley Efron}
\submitdate{August, 2002}
\copyrightyear{2002}
\beforepreface

%% file: abstract.tex
\prefacesection{Abstract}
A Bayesian approach to the classification problem is proposed in which
random partitions play a central role. It is argued that the
partitioning approach has the capacity to take advantage of a variety
of large-scale spatial structures, if they are present in the unknown
regression function $f_0$. An idealized one-dimensional problem is
considered in detail. The proposed nonparametric prior uses random
split points to partition the unit interval into a random number of
pieces. This prior is found to provide a consistent estimate of the
regression function in the $\L^p$ topology, for any $1 \leq p <
\infty$, and for arbitrary measurable $f_0:[0,1] \rightarrow [0,1]$. A
Markov chain Monte Carlo (MCMC) implementation is outlined and
analyzed. Simulation experiments are conducted to show that the
proposed estimate compares favorably with a variety of conventional
estimators. A striking resemblance between the posterior mean estimate
and the bagged CART estimate is noted and discussed. For higher
dimensions, a generalized prior is introduced which employs a random
Voronoi partition of the covariate-space. The resulting estimate
displays promise on a two-dimensional problem, and extends with a
minimum of additional computational effort to arbitrary metric spaces.


%% file: preface.tex
\prefacesection{Acknowledgements}
I thank my parents, who laid the foundations of my character and have
been an unfaltering source of love and encouragement. I also thank my
advisor, Persi Diaconis, whose friendship, support, and insight have been vital
to successfully enduring the Ph.D. process. I also wish to thank the
other members of my examining committee: Jerome Friedman, Bradley
Efron, Susan Holmes, David Siegmund, and Hans Andersen, for their professionalism
and patience. Finally, I thank Hua Tang for putting up with my quirks
all these years.
\afterpreface


%% file: intro.tex
\chapter{Introduction}\label{c:intro}
The binary classification problem is perhaps the simplest regression
problem, but it continues to pose fresh challenges.  In the binary
classification problem, we are given a list of $n$ pairs
$Z_i=(X_i,Y_i)$ each pair drawn independently from an unknown
probability measure $F$. The $X$'s play the role of covariate or
``predictor'' and lie in some abstract space $\mathcal{X}$, while the
$Y$'s are interpreted as a class label and are either 0 or 1. Our goal
is to estimate certain functionals of $F$. Specifically, in the binary
regression problem, we are interested in estimating the regression
function $f:\cX \mapsto [0,1]$. The value of $f$ at a given point $x
\in \cX$ is the conditional probability that $Y=1$ given that
$X=x$. In this way we model the joint distribution $F$ by saying that
to draw an $(X,Y)$ pair from $F$, first draw a covariate $X=x$ from
the marginal distribution of $X$ denoted by $\mu$. Then ``flip'' an
$f(x)$ coin to determine the value of $Y$.

In the classification problem, we are concerned with being able to
predict future $Y$ values. The standard formalization of this task is
that we wish to choose the ``decision rule'' that will minimize the
expected loss incurred; this reduces to the problem of estimating the
set $\{x \in \cX : f(x) > c\}$, for some $c$ that depends upon the
loss (for simplicity, ignore the possibility that $c$ depends on
$x$). There are a great many ways to proceed on each of these
problems, as demonstrated by the vast literature on these
subjects. Some references are given in \secref{s:otherlit}.

In this thesis, I propose a nonparametric Bayesian approach to the
binary classification and regression problems. Specifically, to derive
an estimator, I regard $F$ itself as random. For simplicity, I regard
the marginal distribution of $X$, as known. In this case, putting a
prior distribution on $F$ amounts to putting a prior on $f$. More
generally, one can also put a prior on functions $m$ and suppose that
$\mu(dx)=m(x)\mu_0(dx)$. Some sort of mild restriction, like this one
that all $\mu$ share some dominating measure $\mu_0$, is useful to
avoid technical problems in defining the conditional distribution of
$F$ given the data.

Let $\pi$ denote a prior distribution on $F$, or, more precisely, on
$(f,m)$ pairs. Extend $\pi$ to a joint distribution on $F=(f,m)$ and
the infinite data sequence $(Z_1, Z_2, \dots)$ which, conditionally on
$F$, is drawn independently and identically distributed ($\iid$) from
$F$. Formally, the posterior is the measure $\pi_n(dF) \defined
\pi\left(dF | Z_1=z_1, \dots, Z_n=z_n \right)$. In practice Markov
chain Monte Carlo procedures can be used to generate a sample from the
posterior.

The posterior mean of $f$ is an important summary of the posterior:
its value minimizes the posterior risk under an $\cL^2$ loss. Let
$\widehat{f}$ denote the posterior mean: $\widehat{f}(x)=\int f(x)
\pi_n(df)$. Another important summary of the posterior is the
classification rule which minimizes posterior misclassification
loss. If asked to predict the most likely value of the $Y$'s
corresponding to $X_{n+1}$, \dots, $X_{n+n'}$ all at once, the
decision that would minimize the posterior-expected $0$-$1$-loss is
simply $\delta_i=\1{\widehat{f}(X_i)>\frac{1}{2}}$. Interestingly,
though, if asked sequentially instead of all at once, it is necessary
to update $\widehat{f}$
with each new data point before deciding.

Taking this Bayesian approach assures us that the resulting estimators will have
a clear subjective interpretation. In addition, if the prior $\pi$ is
carefully chosen, the resulting estimators, chosen
indirectly through this Bayesian framework, may have frequentist
advantages over the estimators that might otherwise be proposed. For
example, interesting kinds of shrinkage and averaging occur
automatically within this framework. Subsequent chapters assess the
frequentist performance of these Bayesian estimates by simulation
experiments (\autoref{c:examples}) and theoretically (\autoref{c:proof}).

\section{An Example}\label{s:example}

\begin{figure}
\noindent
\begin{minipage}[t]{\linewidth}
  \centering\epsfig{figure=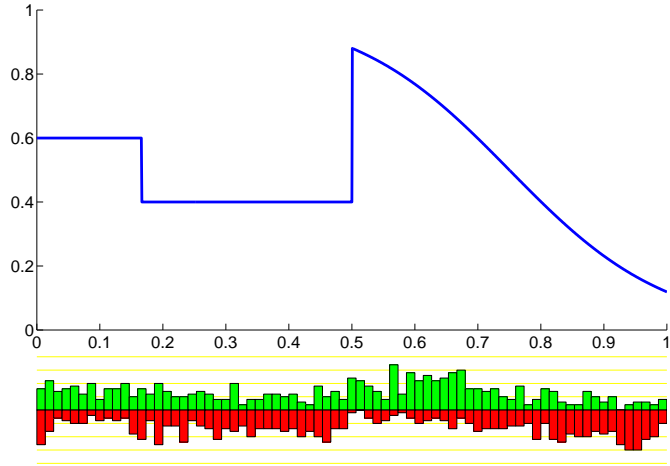,height=3.5in,angle=90}
  \caption{An Example: $f_0(x)$}
  \label{f:f0}
\end{minipage}
\end{figure}

To get started, let us consider the specific case in which $\cX=[0,1]$
and the sampling distribution of the $X_i$, $\mu$, is known to be the
$\U(0,1)$ distribution. Further, let us consider a specific $f=f_0$
which is complicated enough that its estimation should not be too easy
for any of the standard methods; it is piecewise continuous with two
constant regions and a smooth transition region. As shown in
\autoref{f:f0}, $f_0(x)$ is chosen as:
\[
f_0(x)=\begin{cases}
0 \leq x < \frac{1}{6} & 0.6 \\
\frac{1}{6} \leq x \leq \frac{1}{2} &  0.4\\
\frac{1}{2} < x \leq 1 &
\frac{\phi_\sigma(x-\frac{1}{2})}{\phi_\sigma(x-\frac{1}{2})+\phi_\sigma(x-1)}
\end{cases}
\]
where $\phi_\sigma$ is the density of a normal with mean 0 and
standard deviation $\sigma=0.25$.

The two histograms in \autoref{f:f0} summarize a simulated data set of
1024 data points that was drawn from this model. The green histogram
is a histogram of the heads. The red histogram is a histogram of the
tails; it is drawn upside down to facilitate comparison with the
green histogram. To more easily interpret this display, notice that if
we take the sum of the corresponding green and red bins at each point,
we recover a histogram of the marginal distribution, which is
uniform. Furthermore, the ratio of the height of a green bin to the
corresponding red bin represents the empirical odds of a head in that
bin. Near $x=0.5$, for example, notice the sharp transition from
nearly equal green and red bins (on the left) to much longer green
than red bins (on the right). The performance of this posterior mean
estimator is compared with a variety of more
conventional methods in \autoref{c:examples}. This example is used to
illustrate the present approach in the rest of this introduction.

\section{A Nonparametric Prior on Regression Func\-tions}\label{s:priorintro}

To specify a prior $\pi$ over functions $f:[0,1]\mapsto[0,1]$, I explain how
to choose $f$, at random from it. This completely specifies a prior on
the probability distribution $F$, since I consider the marginal
distribution $\mu$ to be known. The prior on $f$ will concentrate on
locally-constant step functions. To choose a step function at random,
first, choose $K$, the number of locally constant intervals, where:
\[
P(K=k)=(1-\alpha)\alpha^{k-1} \quad \textrm{for $k=1 \dots \infty$}
\] 

That is, $K$ is $\Geometric$ with parameter $1-\alpha$. Ultimately, the
choice of $\alpha$ must be specified by the user. For the examples in
this thesis, I have used $\alpha=\frac{1}{2}$ unless otherwise
noted. This choice seems to perform well. For further discussion of
how to choose a prior on $K$ which results in provably
consistent estimators, see \autoref{s:priordisc}.

Now, conditional on $K=k$, choose $V_1, V_2, \dots, V_{k-1}$ \iid\ from
$U(0,1)$. Let $V_{(i)}$ (for $i=1 \dots {k-1}$) be the $i'th$ ordered value
of the $V$'s. This produces $k$ intervals:
\[
I_1=\{0 \leq x < V_{(1)}\}, I_2=\{V_{(1)} \leq x < V_{(2)}\}, \dots,
I_{k}=\{V_{(k-1)} \leq x \leq 1\}
\]
If $k=1$, simply take $I_1=[0,1]$. Finally, conditional on $K=k$,
choose $k$ values $S_i$ (for $i=1 \dots k$) \iid\ from $U(0,1)$. This
generates the random function $f$:

\[
f(x)=\sum_{i=1}^{k} S_i \1{x \in I_{i}}
\]

\section{Sample Results}\label{s:introresults}
\smallfig{unifpost_plotsimple}{-90}{A Sample Result: the Posterior Mean}

Conditioning this prior on the data described in the earlier example
section results in a posterior distribution on regression functions
$f$. Applying the techniques in \autoref{c:implement} to sample from
the posterior results in a long list of sampled functions. Although
the prior specified a $\Geometric(\frac{1}{2})$ prior distribution on
the number of locally-constant pieces $K$ in each function, when drawn
from the posterior, functions typically have at least 5 pieces. Taking
the average of these functions (at every $x$) gives an estimate of the
posterior mean; this is illustrated in
\autoref{f:unifpost_plotsimple}. Further discussion of results like
this can be found in \autoref{c:examples}.

\section{The Partitioning Approach}
Step functions, such as the functions $f$ that the prior $\pi$
concentrates on, have the advantage of great mathematical
simplicity; but, if specified in sufficient detail, they can
approximate general functions. For a multi-dimensional regression
function, it is natural to generalize this idea by considering some
partition of the covariate space $\cX$ into a number of pieces and
constructing a function that takes a different value on each piece. A
procedure that uses a wide variety of geometric shapes to partition
the space might be able to find a partition with the right structure
to approximate the unknown regression function. Ideally, the partition
would be no more complex than necessary to achieve a good
approximation. If the unknown regression function has certain global
features that the chosen partition can be adapted to, this idea becomes very
powerful. Instead of merely ``borrowing strength'' locally, like
ordinary smoothing estimators do, a partition-based estimator borrows
strength across the whole range of a partition element. As a simple
example, if the true regression function does not depend on one of the
covariates, the partition elements do not need to break up space along
this dimension at all; this results in larger partition elements and
more efficient estimation of the success probability on each of the
pieces. Similarly, if the regression function is almost flat in some
large chunk of space, then this whole region can become a single
element. If the level-sets of the regression function have smooth
boundaries, perhaps the partition elements can be chosen to follow
these contours. Finally, since the best partition for the unknown
function is unknown, it makes sense to average together the
approximations found by a variety of partitions that make the data
likely. This is exactly what the posterior mean estimate will do
automatically. Chapter~\ref{c:voronoi} shows one way in which this
idea can be applied using Voronoi partitions. For this partition, a
prior distributes seed-points in the covariate space; each seed
corresponds to a partition element, the one that consists of all points
closer to that seed than any other seed. By placing the seeds
appropriately, the partition can be fine or coarse as needed. The
boundary between partition elements is itself a hyperplane whose
position and orientation can be controlled through the placement of
the seeds.

Other authors have recently considered similar priors with good
success. Their work is discussed in \autoref{s:relatedpriors}. To the
best of my knowledge, this thesis presents the first theoretical
examination of consistency issues for priors of this sort
(c.f. \autoref{c:proof}, with discussion in
\autoref{s:priordisc}). Additionally, I present a detailed assessment
of the empirical performance of my methods on certain novel simulation
experiments. The comparison with bagged CART regression trees in
\autoref{c:examples} is especially interesting.

\section{Outline}\label{s:outline}

Chapter~\ref{c:lit} gives a literature
review. Chapter~\ref{c:implement} shows how to (approximately) compute
samples from the posterior and the posterior
mean. Chapter~\ref{c:examples} trys out the method on examples and
carefully compares its performance with that of a variety of existing
methods. Chapter~\ref{c:proof} gives sufficient conditions on the
prior under which it provides universally consistent estimates of
$f$. Chapter~\ref{c:voronoi} describes a different prior which extends
these ideas to general metric spaces by employing random Voronoi
partitions. Certain modifications are explained that make the proposal
more practical and its performance on an example is shown. Finally,
the afterword gives a philosophical argument that advocates the use of
Kolmogorov-complexity in future statistical thinking. 


%% file: lit.tex
\chapter{Literature}\label{c:lit}

This chapter reviews and discusses the literature on three
subjects. The first section reviews some theoretical results concerning the
frequentist performance of Bayesian procedures. The second section
gives a survey of some of the
work done by authors on related Bayesian efforts. The final section briefly
surveys some salient examples of alternative approaches to the
classification problem. 

\section{Theoretical Results}
The frequentist performance of Bayesian methods is of fundamental
interest in statistics. Given a large sample from a smooth,
finite-dimensional statistical model, the situation is quite well
understood. The Bernstein-von~Mises
theorem~\cite{lecam:yang,freedman:1999} shows that the Bayes estimate and the maximum likelihood
estimate will be close. Furthermore, the posterior distribution of the
parameter vector around the posterior mean is close to the
distribution of the maximum likelihood estimate around the truth: both
are asymptotically normal with mean 0 and the same covariance
matrix. Unfortunately, though, in more general circumstances, such as
those needed for this work, the situation can be much more complex. In
particular, the basic model is based on an infinite hierarchy of finite
dimensional models. Moreover, even for a given finite dimensional
submodel, the dependency of the likelihood function on the parameter
is not smooth; the functions are allowed to take jumps. Consequently, a more
general theory is needed.

This section reviews some of the literature on this subject with
a focus on results that address the question of consistency: i.e. as
the number of data points tends to infinity, will the Bayesian
estimate converge to the true value (in some suitable sense) almost
surely (resp. in probability)?  The literature contains a number of
useful and quite flexible positive results, but also a variety of
interesting negative examples showing that the regularity conditions
under which the theorems hold are not to be taken lightly. A good
introduction to these issues is by Diaconis and
Freedman~\cite{diac:free:1986}. Throughout this section, the reader
may envision a family $\P_\theta(dx)$, a prior $\pi(d\theta)$, and
posterior $\pi(\theta | x_1, \dots, x_n)$, where the $x_i$ are drawn
\iid\ from $\P_{\theta_0}(dx)$. Consistency means that the posterior
concentrates at $\theta_0$ for large samples.

Doob~\cite{doob} established a fundamental result under minimal
regularity assumptions using a martingale convergence
argument. Roughly speaking, the result states that if
consistent estimators exist at all, then a Bayes procedure will
provide an almost surely consistent estimate of the true parameter
$\theta$ under sampling from the $\theta$ distribution for any
$\theta$ in some set $B$ which has prior probability of $1$. Notice,
though, that this does not specify if consistency will obtain at any
{\em particular} point of interest $\theta_0$, unless $\theta_0$
happens to be a point-mass of the prior, or unless it possible to determine
$B$ by some more detailed line of argumentation.

Freedman~\cite{freedman:1963} considered the case in which the
observations are discrete. If the set of possible observations is
finite, the posterior is consistent exactly for parameter values in
the topological support of the prior. The countably infinite case is
more complex. He constructs a class of examples showing that it is
possible to construct a prior which assigns positive mass to every
(weak star) neighborhood of the true parameter value, but for which
the posterior converges to a point mass at some other (chosen)
parameter value. Furthermore, he finds a prior which assigns positive
prior mass to every (weak star) open set of parameters, but for which
the posterior is consistent only at a set of parameters of the first
category. The reader should note that this prior did not assign mass
to all entropy-neighborhoods. This sort of subtle distinction can make
all of the difference and explains the necessity of some such
assumption in the following consistency theorems. He introduces the
``tail-free'' priors for the the countably-infinite case and
demonstrates that these are always consistent.

Lorraine Schwartz~\cite{schw:1965} explored the question of
consistency in a very general setting. She extended Doob's result to a
broad class of loss functions~\cite[lemma 4.2]{schw:1965}.  She also
found sufficient conditions for the posterior to be consistent under
$\iid$ sampling. These conditions, she says, are ``of an essentially
weaker nature'' than the conditions established for the consistency of
maximum likelihood estimators. Nevertheless, she constructs an example
where the maximum likelihood estimate is consistent and the estimates
based on certain priors are not. The example (\cite[example 3]{schw:1965}) involves a simple parametric family of densities
which satisfies Wald's conditions, thereby guaranteeing that the
maximum likelihood estimate will be consistent, but for which the
posterior can be inconsistent. The consistency of the posterior in
this case, is found to depend critically on the amount of mass that
the prior ascribes to small neighborhoods of the true parameter value;
if this mass shrinks too quickly, the prior ``ignores'' the data. One
clever aspect of her construction is the way the densities are
parametrized. Parameter values close to the target value $\theta_0$
correspond to densities that are close to the $\theta_0$-density in an
$\cL^1$ sense, but which are farther and farther away in
Kullback-Leibler discrepancy. In fact, there is only one point in
parameter space (the true parameter) that has Kullback-Leibler
discrepancy from the truth smaller than $\epsilon$, for $\epsilon$
sufficiently small.

Schwartz then shows that the posterior will be consistent under $\iid$
sampling under two basic conditions. First, the prior should have
positive mass on Kullback-Leibler neighborhoods of the true parameter
(defined in \autoref{s:notation} of this thesis), and second, the
model class should not be too rich; specifically, she requires that
uniformly consistent tests of the hypothesis that $\theta=\theta_0$
against the alternative that $\theta$ lies outside a given (open)
neighborhood of $\theta_0$ exist.

It is not always obvious how to verify the later property
directly. Modern authors have employed entropy-type bounds to
guarantee their existence. Ghoshal, Ghosh, and van der
Vaart~\cite{ggv} state a theorem (\cite[theorem 7.3]{ggv}) which
proves that the posterior converges at a certain rate if certain
uniform tests exist (and the prior mass is suitably distributed) and
go on to find a variety of entropy-type conditions that suffice to be
able to construct the necessary tests. Shen and
Wasserman~\cite{shen:wasserman} show related results, requiring
slightly different conditions on how mass needs to be allocated--they
do not a make a connection with testing. Barron, Schervish, and
Wasserman~\cite{barron:schervish:wasserman:1999} find sufficient
conditions for the posterior to be consistent; their results are
reviewed and then used in \autoref{c:proof}.

It should be noted that these various conditions for consistency are
not necessary, but merely sufficient. Nevertheless, it is important to
treat this subject with care because of the variety of examples for
which consistency fails. 

Barron, Schervish, and Wasserman also give an interesting example
where consistency fails. In this example, they show that the prior
puts too much mass on a very rich class of models that will be able to
match any spurious structure that the data might have by chance,
overwhelming the true parameter. Furthermore, lest the reader get the
wrong idea, inconsistency does not only occur in artificial
examples. A series of ``natural'' yet still inconsistent estimators
for the symmetric location problem are discussed by Diaconis and
Freedman~\cite{diac:free:1986}. In addition, the binary regression
example explained in the next section has a natural motivation based on
conditional exchangability.

\section{Related Bayesian Work}\label{s:relatedpriors}

The following subsections contain a review of work by other authors
that is closely related to this thesis. It is followed my a brief synopsis
of the contributions that this thesis makes to the literature.

\subsection{A Dyadic Prior for Binary Regression}\label{s:DFprior}
The most relevant examples for the work of this thesis are the
nonparametric binary regression examples of Diaconis and
Freedman~\cite{diac:free:1993, diac:free:1995}. They use a different
prior; call it $\piDF$, a hierarchical, dyadic prior on $f$. To
describe $\piDF$, let $A_k$ be the set of intervals which result from
partitioning the unit interval into $2^k$ equal pieces. Let $\F_k$ be
the subset of functions which are constant on all intervals $a \in
A_k$. Finally, fix a prior distribution $\kappa$ on the non-negative
integers. Assume, for simplicity, that $\kappa(k) >0$ for all $k$. To
draw $f$ from $\piDF$, draw $K$ from $\kappa$ and then, conditional on
hierarchy-level $K=k$, draw $f$ uniformly at random from $\F_k$. In
effect then, at level $k$ one draws $2^k$ independent $\U(0,1)$ random
variables to describe the success probability on each of the $2^k$
pieces.

They show that for any $\kappa$ and any $f_0$ (except possibly for
$f_0 \equiv \frac{1}{2}$), the posterior estimates are consistent (in
the sense that any $\mathcal{L}^1$ neighborhood of $f_0$ has posterior
probability tending to $1$ a.s.). Remarkably, however, for $f_0 \equiv
\frac{1}{2}$, the posterior can be an inconsistent estimate if the
tail of $\kappa$ is sufficiently heavy.  Specifically, let
$\lambda_k=-\log (\kappa(K \geq k) ) /k $. Then if $\limsup \lambda_k
> \lambdac=2^{-\frac{1}{4}} \approx 0.841$, the posterior is
inconsistent at $f_0 \equiv \frac{1}{2}$. On the other hand, if
$\limsup \lambda_k <\lambdac$, the posterior is consistent for any
$f_0$. To put this in perspective, for $\kappa(k)=(1-\beta)\beta^k$ (a
shifted $\Geometric(1-\beta)$ prior), $\limsup
\lambda_k=-\log(\beta)$. The critical value for $\beta$ is
$\exp(-\lambdac) \approx 0.431$; for larger $\beta$ (longer tails)
inconsistency will occur (but only for $f_0 \equiv \half$).

This result is substantially stronger than the result I have
obtained for my prior $\pi$. In particular, applying the same
(general) method of proof that I employed to prove consistency for
$\pi$ to $\piDF$ yields only the result that $\piDF$ is consistent if
the tails of $\kappa$ drop off at least as fast as those of a
Poisson. (Recall, that at level $k$, $\pi$ only divides $[0,1]$ into
$k$ intervals, but $\piDF$ divides it into $2^k$.) Their method of
proof is direct: using Bernstein's inequality, Poissonization, and
special features of the prior. My method of proof is indirect; it uses
general results that employ entropy-type bounds.

There are striking similarities between $\pi$ and $\piDF$. In fact,
$\pi$ is equivalent to a suitably randomized $\piDF$. To achieve this,
it is not enough to simply randomize the dyadic split points. Instead, recall
that $\piDF$ has an alternative interpretation in terms of binary
sequences. At hierarchy-level $k$, $\piDF$ is uniform over $\F_k$. This
corresponds to independently assigning uniform success probabilities to
each binary sequence of length $k$. Here is an alternative way
to draw $f$ from $\pi$. Draw $g$ from $\piDF$ and interpret $g$ as
function on binary sequences of length $k$ ($k$ depends on $g$). Let
$V_i$ ($i=1, \dots, k$) be $\iid$ $U(0,1)$ random variables. To any
point $u \in [0,1]$ associate the binary random variables
$\eta_i(u)=\II(u \leq V_i)$ ($i=1, \dots, k$). Define $f$ via
$f(u)=g((\eta_1(u), \dots, \eta_k(u) ))$. Note that only a small
fraction of possible binary sequences are realized in this manner (at
level $k$ (which ranges from $0$ to $\infty$ under $\piDF$), $k+1$
sequences out of the full set of $2^k$ possible sequences are achieved).  

\subsection{Bayesian CART}

Two other closely related priors can be described as Bayesian versions
of the CART algorithm. This was pursued by Chipman, George, and
McCulloch, whose prior closely parallels the choices made in the
original CART algorithm~\cite{cgm:1998a, cgm:1998b, cgm:2000a,
cgm:2000b}. Here is a description of their prior when the covariate
space is $\cR^p$. Their prior starts with a root node (which
represents the whole space); this node is then recursively
partitioned in a random way. For each node, randomly choose whether to
split it or not, then choose a coordinate to split on, then choosing a
split point (i.e. the cutoff value) randomly {\em from among the
midpoints between the ordered values of this coordinate}; finally each
leaf node is given an independent regression value. The details of how
these decisions are made differ in their particulars from the ones
that I described in the introduction. In early work, these authors
observed that using MCMC to sample from the posterior of this prior
provides a rudimentary (global) search procedure, which has certain
(apparent) advantages over the {\em greedy} search procedure commonly
implemented in CART-type algorithms. In later work, they examined and
computed the (approximate) posterior mean (working primarily on the
least-squares white-noise regression problem) and found that it had
good performance. They also considered extended priors that modeled
the regression values as additively (not independently)
generated~\cite{cgm:2000b}.

Denison, Mallick, and Smith, independently considered another version
of Bayes\-ian CART~\cite{dms:1998b, dms:1998a, dms:1998c}. For
one-dimen\-sion\-al problems they propose using random splines (the prior
I use is essentially a special case of this prior). They consider some
of the regression examples that are standard in the wavelet literature
and show that their spline methods perform equally well. Additionally,
they propose a Bayesian version of Friedman's MARS which puts a prior on
functions that are constructed by adding together random spline-type
ridge functions. Denison, Adams, Holmes, and Hand discuss the
usefulness of random partitions in this paper~\cite{dahh:2002}.

Very recently, Denison, Holmes, Mallick, and Smith have written a
book~\cite{dhms:2002} which surveys some related Bayesian regression
schemes, including a Bayesian method for (multiple class)
classification using Voronoi partitions that is very closely related
(albeit independent of) the work that I present in
\autoref{c:voronoi}. The book also discusses Bayesian wavelet methods,
and an interesting Bayesian nearest-neighbor prior. As a default
prior, they recommend assuming that every model in a ``single
dimension'' is equally likely, and each dimension is equally probable,
\apriori. This ``flat prior,'' they claim, should serve perfectly well
because of the, ``natural tendency'' for the marginalized likelihood
to penalize complex models:

\begin{quotation}
On the face of it, we might be concerned that the flexible modeling
strategy we advocate might be prone to overfitting the data by adding
too many basis functions. Indeed, many papers found in the literature
advocate explicit priors on the model space that penalize the
dimension of the model. However, throughout this book we argue that
such a measure is unnecessary.  The Bayesian framework contains
a natural penalty against over [sic] complex models, sometimes called
{\it Occam's razor}, which essentially states that a simpler theory is
to be favoured over a more complex one, all other things being
equal.
\end{quotation}

There is no consideration given to the possibility that this might
give rise to inconsistent estimates (e.g. as in the Diaconis and
Freedman non-parametric regression example explained earlier); indeed
there are few theoretical considerations at all in the book. Their
explanation of why the Markov chain techniques that they develop
should actually give meaningful samples from the posterior appeals to
Green's reversible jump~\cite{green:1995}. The explanation given is
vague and ultimately they decide to avoid the issue and appeal to the
fact that their chains are discrete. The chains in
\autoref{c:implement} of this thesis involve a continuous state space
and do not simply avoid this issue by discretizing the continuous
modeling space as these authors seem to do. 

Overall, the book emphasizes main ideas, algorithms, and results. It
seems that for every existing regression technique, they want to
demonstrate that they can make a ``Bayesian'' version of it too. The
book does not emphasize subjectivism, but rather adopts an
``$\cM_{\text{\it open}}$'' perspective to Bayesian modeling: ``we
never believe that the true model lies in the set of possible
models.'' The book does do a good job of supplying default priors for
a wide variety of possible parametric models. Similarly, Denison's
thesis~\cite{denison:1997} emphasizes the wide variety of problems to
which Bayesian partitioning methods of this sort can be applied.

\subsection{Poisson Rate estimates using Random Partitions}
Green~\cite{green:1995}, and
Scargle~\cite{scargle} develop priors on piecewise constant
functions on the real line and $\bbR^d$ using Voronoi
cells. Their priors are quite similar to the ones developed in this thesis,
but are intended to address the problem of estimating the rate function of a
$\Poisson$ process. In principle, one could apply their techniques to
the problem of binary regression by generating an estimate of the rate
function of the ``heads'' process and the ``tails'' process separately
and then combining the results. I do not think that this has been
tried and it seems substantially less ``natural.''

Green applies his method to a coal mining dataset and a synthetic
two-dimen\-sion\-al example. For these example, Green assumes that an
individual cell's rate-para\-met\-er is drawn independently from a
$\Gamma(\alpha,\beta)$ prior. For the one-dimen\-sion\-al case he
advocates a prior which ``probabilistically'' spaces out the
change-point locations; specifically, if there are $j$ change-points,
the ordered locations of the change-points are distributed like the
even order statistics of $2j+1$ independent uniform values. He argues
that this is good because it prevents small change-point intervals from
entering into the posterior. For the two-dimensional example, the
generating points of the Voronoi partition are drawn independently and
uniformly. Green's methods are given, in part, as
examples of his ``reversible jump'' MCMC technique. This technique has
become an accepted part of MCMC practice, but is not accepted by all
experts in MCMC theory because it does not lay down in a
straightforward ``theorem-proof'' manner the necessary conditions and
consequent conclusions. For this reason, detailed verifications for
the chains used in this thesis are given in \autoref{c:implement}.

Scargle's work is applied to astronomical data; he concentrates on the
problem of finding the mode of the posterior, rather than the
posterior mean. Fortunately, he and coworkers have developed a way of
computing this mode in the one-dimensional case exactly and
efficiently using a dynamic programming
approach~\cite{scargle:dynamic}. Instead of giving each cell an
independent value, Scargle gives each cell a (logical) ``color'' and
then associates each unique color with an independent
rate-parameter. This allows him to use a fine partition and then group
``chunks'' back together into more complicated shapes. The way he
forms this partition is also different; in particular his ``prior'' is
data dependent, but not quite in the way of the ``prior'' that I
consider in \autoref{c:voronoi}. Rather, the data is used once and for
all to generate the fine Voronoi partition of space that results from
using all of the data points as generators. These cells are then
``clumped'' (i.e. given a logical color) and the clumps are given an
independent rate parameter.

\subsection{Bayesian ``Image'' Analysis}
M{\o}ller and Skare~\cite{moller:skare} apply their work to reservoir
modeling and connect their work to efforts in Bayesian image analysis
(including Markov random fields). They use a random Voronoi partition
of the data and assign each partition element a random color (in a way
that depends only the colors of neighboring cells). They supply
several further references to work in Bayesian image analysis which
use Voronoi cells. From their perspective, to calculate their
posterior they are simulating from a special ``marked point'' process.
The generators of the Voronoi cells are regarded as point set that has
been drawn from a homogeneous Poisson process of rate $\beta$ on the
unit cube. In the simplest case, the marks or ``colors'' of these
points are just integers from $1$ up to $M$ that have been drawn
independently. More generally, according to their prior, the
conditional distribution of the coloring of cells given is an Ising or
Potts model. The graphical structure of this model is determined by
consideration of which Voronoi cells are neighbors, and the $\theta$
parameter is chosen to reflect their prior belief that neighboring
cells tend to be of the same color. They consider two problems. The
first is a simulation experiment in which a ``true'' binary image is
degraded with Gaussian noise. The second is a three dimensional
reservoir problem based on real data. It is supposed that a certain
three dimensional cube (the reservoir) consists of 4 different types
of rock. The rock types are observed along seven vertical lines,
representing the observations of rock that were made as seven wells
were dug into the reservoir. In both problems, the true object to be
recovered is itself a certain ``coloring'' of space (i.e. rather than
a continuous regression function). For the MCMC computation of their
posterior they apply the birth-death type Metropolis-Hastings
algorithm for point processes, as studied by Geyer and
M{\o}ller~\cite{geyer:moller:1994} and claim that their target
distribution satisfies a local stability condition (see
Geyer~\cite{geyer:1999}, Kendal and
M{\o}ller~\cite{kendal:moller:2000}, and M{\o}ller~\cite{moller:1999})
so that the MCMC is actually geometrically ergodic.

\subsection{Polya Trees}
Finally, Polya trees~\cite{lavine:1992} and especially randomized
Polya trees~\cite{paddock:1999} deserve to be mentioned. The basic
Polya tree puts a prior on distribution functions on the unit
interval. The unit interval is divided recursively in a dyadic binary
way and mass is allocated to each piece of the partition in a
stagewise manner by first determining how much of the mass that is
available will be on the left versus the right half and then
continuing with such determinations layer by layer. Each of these
assignments is ultimately determined by independent $\Beta$ random
variable, whose parameters depend upon its location in the ``tree.''
If a suitable choice of these parameters is made the result prior on
distribution functions concentrates on distributions that are
absolutely continuous with respect to Lebesgue measure. The essential
advantage of Polya trees is that the posterior of Polya tree prior is
easily and analytically computable, being itself another Polya
tree. For randomized Polya trees, the partitioning scheme is
independently ``jittered'' at random in a particular
way~\cite{paddock:1999}. A Hybrid MCMC can be employed to sample from
the randomized Polya tree posterior which uses a Gibbs step to take
advantage of the ease with which the (internal) Polya tree posterior
can be computed. Both methods can be extended (essentially by taking
``direct products'') to put a prior on distributions on the unit cube.

\subsection{The Contributions of this Thesis}

Reviewing the depth and breath of the literature reviewed above may
leave the reader in doubt about the contributions of this
thesis. After all the one-dimensional prior that I consider is
essentially a special case of the univariate spline model and the idea
of using Voronoi partitions is certainly not new, although effective
Bayesian methods using them only started springing up fairly
recently. 

Still there is room for careful analysis. This thesis establishes that
the posterior is consistent under suitable conditions on the prior and
for any measurable regression function (see \autoref{c:proof} for
details): an issue which none of the ``Bayesian CART'' or ``Voronoi
Partition'' authors address at all. This thesis also gives an explicit
Markov chain Monte Carlo algorithm (see \autoref{s:algo}). Broadly
speaking it is a fairly standard birth-death Markov chain as
considered by Geyer and Moller~\cite{geyer:moller:1994}, but the
technicalities of the analysis seem to be somewhat different. This
thesis proceeds to show in detail that it satisfies detailed balance
by direct self-contained argumentation; further, the chain is shown to
have an ergodicity property (see \autoref{s:ergoproof}). These
considerations are often glossed over in modern writing.

On the more practical side, \autoref{c:examples} scrutinizes the
behavior of the posterior mean estimate under a variety of carefully
designed simulation experiments. These experiments both serve to
analyze the posterior mean and to give insight into the relationship
between Bayesian methods and their classical counterparts. See for
example the discussion of CART and bagging in \autoref{s:cartexp}.

\section{Other Approaches}\label{s:otherlit}
The literature on classification and regression methods is huge; the
interested reader is urged to consult good modern books on the subject
like {\em The Elements of Statistical Learning,} by Hastie,
Tibshirani, and Friedman~\cite{HTF:2001}. The following paragraphs
outline some of the methods that have had the most impact upon the author.

In the statistics literature, classical approaches to the
classification and binary regression problem include logistic
regression, Fisher's discriminant analysis, and projection pursuit
methods. Logistic regression specifies that the success probability
regression function is such that its log-odds follows a linear model
with a user specified basis (e.g. by using polynomial or spline
functions of the covariate-data) and estimates the parameters by
maximum likelihood. Model selection is commonly performed using
classical methods to select a subset of the covariate
variables. Fisher's discriminant analysis finds a hyperplane which
``optimally'' separates the two classes using a within versus between
variance criterion. Projection pursuit seeks an interesting linear (or
sometimes nonlinear) projection of the covariate-data onto a lower
dimensional subspace (e.g. $\bbR$). Various criteria have been
proposed to define ``interesting,'' some of which are suitable for the
classification problem. Each of these methods has undergone a variety
of generalizations and tweaks to address a wider range of problems over
the years.

The first {\em general} method to solve the classification problem
automatically was the $k$-nearest neighbor
approach~\cite{cover:hart:1967}. $k$-Nearest neighbor estimates are
known to be universally consistent if $k=k(n) \tendsup \infty$ slowly
enough~\cite{devroye:gyorfi:lugosi:1996}. Their convergence, however,
especially in high dimensional problems, can be slow in
practice~\cite{friedman:1996}.

Local regression methods are a clever extension of
this approach. To predict at a given point, instead of averaging the
values given at the neighbors, they fit a low-order linear model to a
locally-weighted version of the data set~\cite{clev:load:1996}.

Trees~\cite{bfos:1984} and neural nets~\cite{ripley:1996} differ in
that they search through a globally-parametrized class of
functions. In all of these methods, cross-validation is often employed
to estimate frequentist ``out-of-sample'' performance and select a
regularization parameter which governs the trade-off between bias and
variance~\cite{HTF:2001}.

Wavelet methods are in some ways a compromise between the local and
the global approaches mentioned above. They fit an explicit global
linear model to the data, but the basis elements in this model are
carefully constructed to maintain ``localization'' (in space and
frequency domains). They boast powerful asymptotic compression and
approximation properties, computationally efficient transforms, and
can employ special thresholding methods which ``optimally'' choose
which coefficients in the model are kept~\cite{donoho:johnstone:1994}. However,
their practical use seems to remain concentrated on the case of
regularly-spaced regression data. Some recent papers address this
shortcoming~\cite{daubechies}.

Support vector machines (SVMs)~\cite{vapnik:1996} employ a
``kernel-trick'' to reduce consideration of a certain
globally-parametrized model class to consideration of an equivalent
linear model class in an abstract Hilbert space. The
estimated decision rule corresponds to the solution of a convex
optimization problem. This objective function still involves an
unknown regularization parameter. In practice, this parameter is often
chosen by cross-validation, but, in principle, it can be chosen
through consideration of the structural risk minimization (SRM)
paradigm. The advantage of using the SRM paradigm is that one obtains
provably valid confidence statements about the error rate that will
obtain on future data. Moreover, these confidence bounds improve at an
exponential rate in the number of data points. With realistic sample
sizes, however, the bounds are often too crude to be of practical
use. There are hidden connections between SVMs and (1) Bayesian
methods employing Gaussian-process priors on the regression function
(including the generalized spline methods of Wahba~\cite{wahba,wahba:svms}) (2)
projection pursuit regression~\cite{coram:svm:pp}.

Bagging~\cite{breiman96bagging} and
boosting~\cite{freund:schapire:1996a} are meta-algorithms that
``boost'' the performance other classification algorithms (especially
trees) by taking carefully chosen weighted averages of the results of
the boosted (respectively, bagged) algorithm. There are close
connections between boosting and the Lasso penalty, which itself is
closely related to the least angle regression method
(LARS)~\cite{efron:lars}.
\nocite{friedman:mars}


%% file: implement.tex
\chapter{Computing the Posterior}\label{c:implement}

\renewcommand{\v}{\mathbf{v}} \newcommand{\s}{\mathbf{s}} 
\newcommand{\V}{\mathbf{V}} \renewcommand{\S}{\mathbf{S}} 
\renewcommand{\varTheta}{\Xi}
\newcommand{\ppi}{\pi'_\circ}
\newcommand{\shatj}{\widehat{s}_j}
\newcommand{\shatJ}{\widehat{s}_{J(x;\v)}}
\newcommand{\emp}{\between}
\newcommand{\bx}{\underline{x}}

This chapter describes a Markov chain Monte Carlo (MCMC) algorithm
that can be used to (approximately) draw samples from the posterior of
the one-dimensional random step function prior $\pi$. This is
essential, because for a complex prior like $\pi$, analytical
evaluation of properties of the prior is intractable. All computation
about the posterior, therefore, is made through (approximately)
generating a large sample from it. Before describing these algorithms,
it is natural to review the prior and then derive a more refined
mathematical expression for the posterior. This exercise has the side
effect of suggesting a more efficient sampling scheme. An informal
sketch of the MCMC algorithm is then given in
\autoref{s:algo}. Additionally, \autoref{s:calcmean} explains an
efficient way to use these samples to calculate the posterior mean. The
interested reader is invited to download an implementation of these
algorithms and others from the author's web page.

To define the algorithm more mathematically, a brief review of Markov
chains and the Metropolis-Hastings algorithm for general state spaces
is given. Section~\ref{s:simplechain} gives a simple example of such a
chain. Section~\ref{s:localchain} gives a mathematical treatment of
the more complicated MCMC algorithm that was only sketched
previously. It also verifies that the Markov chain satisfies detailed
balance with respect to the posterior. This is done mainly to provide
a thorough and relatively self-contained theoretical analysis of the
Markov chain. 

Many of these calculations are essentially of the same type as those
considered by Green~\cite{green:1995}. Indeed, the birth-death move that I
employ is essentially equivalent to the ones that he calls a
``reversible jump.'' The verification of detailed balance is simple
using his results. One need only observe that the transformations
involved in these jumps essentially only permute coordinates;
consequently the absolute value of the determinant of the Jacobian of
this transformation is identically 1. Section~\ref{s:ergodicity}
reviews two theoretical results that give sufficient conditions for a
Markov chain to produce the intended ergodic
sequence. Section~\ref{s:ergoproof} shows that these results are
applicable so that, for example, the posterior mean that is computed
will indeed approximate the intended curve.

\section{Specification of the Prior in One Dimension}\label{s:priorreview}
A review of notation and the prior is in
order. Section~\ref{s:priorintro} gives an informal description of the
prior $\pi$ and \autoref{s:priorspec} gives a more formal
specification of the prior $\pi$ and its parametrized version
$\pi'$. The parameter space $\Theta=\union_{k=1}^\infty \Theta_k$
where $\Theta_k$ parameterizes the class of functions
$f:[0,1]\goesto[0,1]$ that have $k$ locally-constant regions. Any such
function is (essentially) determined by two vectors $\s$ and
$\v$. Vector $\s$ lists the $k$ values that the function takes on each
region (as enumerated from left to right). Vector $\v$ lists (in no
particular order) the $k-1$ locations at which the function
jumps. This explains the definition:
\begin{align}
\Theta_k \defined \{ (k, \v, \s) ~:~ \v \in [0,1]^{k-1}, \s \in [0,1]^k\}
\end{align}
$\Theta_1$ is a special case because there are no splits in a function
that is everywhere constant. Define $\Theta_1= \{ (1, \emp, \s) ~:~ \s
\in [0,1]^1 \}$, where the {\em symbol} $\emp$ represents an empty list
(which is not considered equivalent to an empty set). This is
consistent with the former definition of $\Theta_k$ if one allows the
notation: $[0,1]^0 \equiv \{~ \emp~ \}$.

To specify the prior $\pi$, first select a probability distribution
$\kappa$ on $\bbN$. $\kappa$ represents the \apriori\ distribution of
the number of regions in the unknown function. In the introduction,
the choice $\kappa \eqd \Geometric(\half)$ was suggested. Assume that
$\kappa(k)>0$ for all $k \in \bbN$, where (technically) $\kappa(k)$ is
shorthand for $\kappa(\{k\})$. To pick a value $\theta \in \Theta$
from the (parametrized) prior $\pi'$, first draw $K \dist
\kappa$. Then, if $K=k$, draw $\S=\s$ uniformly from $[0,1]^k$ and draw
$\V=\v$ uniformly from $[0,1]^{k-1}$. Form $\theta=(k,\s,\v)\in
\Theta_k$. This completes the description of $\pi'$.

For a given point $\theta \in \Theta_k$ it is convenient associate a
number of objects. Let $k(\theta)$, $\v(\theta)$, and $\s(\theta)$
stand for the $k$, $\v$, and $\s$ parts of $\theta$ respectively. Let
$\v_{(i)}$ denote the $i$'th ordered value of $\v$.  Additionally, for
$\theta=(k,\v,\s) \in \Theta_k$ associate the function $f_\theta$
whose splits points and values are determined by $\v$ and
$\s$. Specifically, for $k>1$, let $I_1=[0,\v_{(1)}), I_2=[\v_{(2)},
\v_{(3)}), \dots, I_{k-1}=[\v_{(k-2)},\v_{(k-1)}),
I_{k}=[\v_{(k-1)},1]$. If $k=1$, just let $I_1=[0,1]$. Take
$f_\theta(x)=\sum_{i=1}^k \s_i \1{x \in I_i}$.

I have chosen to work with uniform distributions on the splits and
function values. Both of these choices could be varied, e.g. by using
one-dimensional $\Dirichlet$ distributions for the split locations and
using $\Beta$ distributions for the function values where, perhaps,
the parameters of the $\Beta$ distribution are allowed to depend on
spatial position. Diaconis and Freedman~\cite{diac:free:1995} adopt
this level of generality (for the success probability prior), but I
have not found it useful. What {\em would} be useful (but is avoided
for simplicity of presentation) is to extend from binary
classification to the multi-class case. This can be done by
generalizing the $\Beta(1,1)$ prior into a discrete $\Dirichlet$ prior
on the class probabilities. Allowing dependencies among the parameters
would more substantially complicate the analysis.

\newcommand{\data}{\mathcal D}
\section{Representing the Posterior}\label{s:posterior}
The posterior is just the result of conditioning the prior on the
data. The data is the list $\data \defined (x_1,y_1), \dots, (x_n, y_n)$
where for $i$ from $1$ to $n$, $x_i \in [0,1]$ represents where the
$i$'th point occurred and $y_i \in \{0, 1\}$ represents whether the
$\Bernoulli(f(x_i))$ ``coin'' came up heads or tails.
Denote the posterior distribution on $\theta$ given the data by
$\ppi$:
\begin{align}
\ppi(d\theta) \defined \pi'(d\theta ~|~\data)
\end{align}

Provided that the denominator is non-zero and finite (which it will be) the
posterior $\ppi$ has a density $\phi(\theta)=\phi(\theta; \data)$
with respect to the prior $\pi'$, so that
$\ppi(d\theta)=\phi(\theta)\pi'(d\theta)$, where:
\begin{align}
\phi(\theta)=\frac{L(\theta)}{\int_\Theta L(\theta') \pi'(d\theta')}
\end{align}
The likelihood function $L(\theta)=L(\theta; \data)$ is defined by:
\begin{align}
L(\theta) \defined \prod_{i=1}^n
f_\theta(x_i)^{y_i}\left(1-f_\theta(x_i)\right)^{1-y_i}
\end{align}

Recall that for $\theta=(k,\v,\s)$, $f_\theta(x)=\sum_{j=1}^k \s_j
\1{x \in I_j}$, where the intervals $I_j$ implicitly depend $\theta$,
but only through $\v$. To evaluate $f_\theta(x_i)$, then, is simply a matter
of determining for which value $j$ from $1$ to $k$, $x_i \in I_j$ and
then retrieving that $\s_j$. Call the value of $j$ for which $x \in
I_j$, $J(x; \theta)$. Then
$f_\theta(x_i)=\s(\theta)_{J(x_i;\theta)}$. Consequently, $L$ is simply
a certain product of terms of the form $\s(\theta)_j$ or
$1-\s(\theta)_j$. To collect these together, define:
\begin{align}
N^1_j&=N^1_j(\theta; \data)=(\textrm{\# of data points in $I_j$ labeled 1}) \\
N^0_j&=N^0_j(\theta; \data)=(\textrm{\# of data points in $I_j$ labeled 0})
\end{align}
So that for $\theta \in \Theta_k$, (suppressing the dependence on
$\theta$ and $\data$ from the right hand side) $L(\theta)$ can be
expanded as:
\begin{align}
L(\theta)= \prod_{j=1}^k
\s_j^{N^1_j}(1-\s_j)^{N^0_j}
\end{align}

Notice, then, that for fixed model number $k$ and change-point
locations $\v$, $L(\theta)$ only depends on the data through the
(conditionally) sufficient statistics $N^1_j$ and $N^0_j$ for $j=1,
\dots, k$. Moreover, for $\theta \in \Theta_k$, $L$ is simply the
product of $k$ different binomial likelihood functions. This is
intuitively obvious: If I have already decided exactly where the
change-points are, the only remaining parameters are the success
probabilities $\s$. Furthermore, according to the model, if I get a
data-point $(x,y)$, the $x$ necessarily lands in some interval $I_j$
and, then, the $y$ is just the result of flipping an (independent)
$\s_j$ coin.

Also notice that under the prior the $\S_j$'s are independent
$\U(0,1)$ random variables. From the above discussion, it is apparent
that if we condition on $K=k$ and $\V=\v$, the data simply tell us how
many times each of the $k$ ``coins'' with success probabilities $\s_1$
through $\s_k$ came up ``heads'' and ``tails'' as they were
(collectively) flipped $n$ times.  Consequently, under
the posterior, conditioned on $K=k$, $\V=\v$, and on the values of
$N^1_j$ and $N^0_j$, each $\S_j$ is an independent
$\Beta(N^1_j+1,N^0_j+1)$ random variable. We recover this fact by
direct calculation shortly.

These observations can be used to motivate an MCMC scheme that is
substantially more efficient than the naive one that randomly changes
$k$, $\v$, and $\s$ in the standard Metropolis-Hastings
fashion. Namely, we will only have to use Markov Chain Monte Carlo
steps in order to sample $(k,\v)$ pairs from their marginal under the
posterior. If desired, we can then create a complete sample, including a
realization of $\s$ by sampling from the $k$ independent $\Beta$
random variables whose parameters were explained above. To compute the
posterior mean, $\s$ need not be simulated at all. The mean of the $\Beta$
distribution can be computed analytically. This avoids substantial Monte Carlo
error. For details, see sections~\ref{s:algo} and \ref{s:calcmean}.

Continuing with the calculations, write $\theta=(k,\v,\s)$, and let
$C=\{k\} \times A \times B=\{\theta \in \Theta_k: \v \in A, \s \in
B\}$ for $A, B$ measurable subsets of $[0,1]^{k-1}, [0,1]^{k}$
respectively. Compute:

\begin{align}
\ppi(C)=\pi'(C | \data)
&= \int_C \phi(\theta) \pi'(d\theta) \\
&= \frac{\int_C L(\theta) \pi'(d\theta)}
        {\int_\Theta L(\theta) \pi'(d\theta)}
\end{align}
\begin{align}
\int_C L(\theta) \pi'(d\theta)
&= \kappa(k) \int_{\v \in A} \int_{\s \in B}
       L\left( (k,\v,\s) \right) d\s d\v \\
&= \kappa(k) \int_{\v \in A} \int_{\s \in B}
        \prod_{j=1}^k \s_j^{N^1_j(\v)}(1-\s_j)^{N^0_j(\v)} d\s d\v \\
\intertext{If, in particular, $B$ is the rectangle $[a_1,b_1] \times \dots \times
[a_k, b_k]$ we get:}
\int_C L(\theta) \pi'(d\theta)
&= \kappa(k) \int_{\v \in A} 
        \left[\prod_{j=1}^k \int_{a_j}^{b_j} u^{N^1_j(\v)}(1-u)^{N^0_j(\v)}
        du \right] d\v \\
\intertext{The inner integral is a $\Beta$ integral. Consider the
  special case in which $B=[0,1]^k$ (i.e. $C=\{k\}\times A \times
  [0,1]^k$). Then,}
\int_C L(\theta) \pi'(d\theta) &= \kappa(k) \int_{\v \in A} \rho_k(\v) d\v
\end{align}

Where $\rho_k(\v)$ and also $\rho_k$ are defined for $k\ge 2$ by:
\begin{align}
\rho_k(\v) & \defined \prod_{j=1}^k \frac{N^1_j(\v)! ~
  N^0_j(\v)!}{(N^1_j(\v)+N^0_j(\v)+1)!} \\
\rho_k &\defined \int_{\v \in [0,1]^{k-1}} \rho_k(\v) d\v
\end{align}

For $k=1$, $\v=\emp$ and $I_1=[0,1]$ so that $N^1_1$ and $N^0_1$ are
the total number of heads and tails respectively. In this case define
$\rho_1=\rho_1( \emp ) = N^1_1! N^0_1! / (N^1_1+N^0_1+1)!$.

The posterior probability of the $k$'th model is readily computed:
\begin{align}\label{e:postk}
\ppi( \Theta_k ) = \kappa(k) \rho_k / c
\end{align}
Where the normalizing constant $c$ is the sum: ${\sum_{j=1}^\infty
\kappa(j)\rho_j}$.

Now is a good time to notice that for any $k$ and
$\v \in [0,1]^{k-1}$, and any data set $\data$, $\rho_k(\v)$ and
$\rho_k$ are both positive and $\leq 1$. Consequently, the same holds
for $c$.

Let $\lambda_j$ denote Lebesgue measure
on $[0,1]^{j}$. In the special case that $j=0$, let $\lambda_0$ denote
counting measure on the set $\{ \emp \}$. Then the posterior density of
the change-points $\v$ with respect to $\lambda_{k-1}$, given that
model $k$ holds is $\rho_k(\v) / \rho_k$. Informally:
\begin{align}\label{e:postv}
\ppi(\V \in d\v ~|~ K=k)=\rho_k(\v) / \rho_k
\lambda_j(d\v)
\end{align}

Finally, the posterior probability that $\S$ is in rectangle $B=[a_1,b_1]
\times \dots \times [a_k, b_k]$, given that model $k$ holds and that
the change-points are given by $\v$ is indeed the same as $k$
independent $\Beta$'s:
\begin{align}\label{e:betas}
\ppi( \S \in B ~|&~ K=k, \V=\v)
= \frac{1}{\rho_k(\v)}
  \prod_{j=1}^k \int_{a_j}^{b_j} u^{N^1_j(\v)}(1-u)^{N^0_j(\v)} du \\
&= \prod_{j=1}^k \P\left( \Beta(N^1_j(\v)+1,N^0_j(\v)+1) \in [a_j,
    b_j]\right) du
\end{align}

Consequently, if $\shatj$ denotes the expected value of $\S_j$ under the
posterior, given that model $k$ holds and that the change-points are
given by $\v$, then $\shatj$ is just the mean of the
$\Beta(N^1_j(\v)+1,N^0_j(\v)+1)$ distribution:
\begin{align}\label{e:shatj}
\shatj=\frac{N^1_j(\v)+1}{N^1_j(\v)+N^0_j(\v)+2}
\end{align}

\section{Setup}\label{s:mcmcsetup}
This section sets up some basic definitions and ideas that underlie
the algorithm described in the next section. The first definition,
gives a new meaning to the symbol $\cX$ which will be used throughout
the remainder of this chapter. Elsewhere, $\cX$ still stands for the covariate
space. This should not introduce any confusion.

Write $\cX_k$ for the parameter space formed from $(k,\v)$
pairs with $\v \in [0,1]^{k-1}$ and build up the full parameter space
$\cX$ by taking the countable union:
\begin{align}
\cX &\defined \union_{k=1}^\infty \cX_k \\
\cX_k &\defined \{(k, \v) ~:~ \v \in [0,1]^{k-1}\}
\end{align}
Again, $\cX_1$ is a special case. Define $x_0 \defined (1,\emp)$ and let
$\cX_1= \{ x_0 \}$, so that $\cX_1$ is a singleton.

For convenience, let $k(x)$ stand for the $k$-part of $x$ and let
$\bx$ stand for the $\v$-part of $x$. Bear in mind the nuisance that
for $x \in \cX_k$, $\bx \in [0,1]^{k-1}$. Extend this definition to
sets; i.e. for a subset $C \subset \cX_k$, define $\underline{C}$ as
$\{\bx ~:~ x \in C\}$.

For future reference, endow $\cX$ with the $\sigma$-algebra $\cB$
generated by sets of the form ${k}\times B$ where $k\ge 2$ and $B$ is
some (Borel) measurable subset of $[0,1]^{k-1}$ and by the singleton
set $\cX_1$.  Let $\tau$ denote the extension of Lebesgue measure to
$\cX$. That is, for $k \ge 2$ and $B$ a (Borel) measurable subset of
$[0,1]^{k-1}$ define the $\tau$ measure of a set
$C={k}\times B \subset \cX_k$, as the $k-1$-dimensional Lebesgue
measure of $\underline{C}=B$.  To account for $k=1$, let
$\tau(\cX_1)=1$. 

Finally, combining the results in \autoref{e:postk} and
\autoref{e:postv}, the distribution under $\ppi$ of the random point
$X=(K,V)$ in $\cX$ can be computed. In the present notation, $X$ has
the density $\phi(x)$  with respect to $\tau$:
\begin{align} \label{e:phix}
\phi(x) \defined \kappa(k) \rho_k(\v) /c
\end{align}

\section{An MCMC Algorithm}\label{s:algo}
This section contains a description of a (randomized) computational
algorithm to produce a random sequence $\theta_1, \dots, \theta_M$
drawn from the posterior distribution $\ppi$. A more formal version of
this algorithm will be developed in \autoref{s:localchain}. Finally in
\autoref{s:ergoproof}, it will be shown that the generated sequence
has an ergodicity property. The main consequence of this is that if
$g$ is some measurable function with $\int |g(\theta)| \ppi(d\theta)<
\infty$, one can use the average value of $g$ on the sampled values to
approximate the integral, in the sense that:
\begin{align}
\frac{1}{M} \sum_{j=1}^M g(\theta_j) \tendsto \int g(\theta)
\ppi(d\theta)
\end{align}

The algorithm exploits the observation made in \autoref{s:posterior}
that under the posterior, conditioned on $K=k$, $\V=\v$, and on the
values of $N^1_j$ and $N^0_j$, each $\S_j$ is an independent
$\Beta(N^1_j+1,N^0_j+1)$ random variable. Since $\Beta$ random
variables are easy to simulate and work with analytically, it suffices
to simulate $(k, \v)$ pairs from the (marginalized) posterior, instead
of the full $\theta=(k,\v,\s)$. 

Essentially, the MCMC algorithms that are developed in this chapter
are a minor variation of the usual birth-death MCMC approach that is
often used to simulate point processes. This approach is described by
Geyer and M{\o}ller~\cite{geyer:moller:1994}. They claim geometric ergodicity for
their Markov chain, under suitable restrictions on the sampling
density. One technical distinction from these approaches is that when
simulating a point process, the fundamental object is a subset of
points $\{\v_1, \dots, \v_{k-1}\}$; the theoretical treatment given in
this chapter considers instead mathematical objects of the form
$\bigl(k,(\v_1, \dots, \v_{k-1})\bigr)$.

As far as the computations are concerned, though, there is no
distinction between these formally different objects. Both could be
represented on the computer {\em operationally} as a simple list of
numbers called $x$ with each number in the list specifying the
location of a particular change-point. The algorithm merely assumes
that it can call a function $\phi(x)$ that evaluates to the positive
real number defined by \autoref{e:phix}.

In particular, one can use $k(x)$ to
find out that there are $k-1$ elements in this list. Furthermore the
computer has no problem removing elements from the list all the way
down to the empty list, or (ideally) adding elements one-by-one
indefinitely. To agree with the notation of the previous section, if
$k(x)=k$, then for $1 \leq j \leq k-1$, write $\bx_j$ for the $j$'th
element of the list. Write $x_0$ for the empty list.

Let $M$ be the number of Monte Carlo samples that are desired. The
algorithm simulates the workings of a Markov chain and generates the
realizations: $x_1, x_2, \dots, x_M$. It is assumed that for $j$ from $1$
to $5$, $p_j$ is a positive number, and that these numbers sum to
one. These represent the mixture probabilities with with $5$ component
Markov chains are combined. For my computations, through a mixture of
intuition and trial-and-error, I chose $p$ as in
\autoref{t:pvector}. These values are by no means optimal. (The
irregular numbers quoted here result from standardizing simpler ones
so that they add to 1.)

\begin{table}[h]
\centering
\begin{tabular}{|c|c|c|c|c|}
\hline
 $p_1$  &       $p_2$    &     $p_3$    &     $p_4$  &   $p_5$      \\
\hline
 0.1429  &  0.1429  &  0.2381  &  0.0476  &  0.4762 \\
\hline
\end{tabular}
\caption[MCMC Mixture Probabilities]{{\em MCMC Mixture Probabilities}}
\label{t:pvector}
\end{table}
\begin{algorithm}[H]\label{a:mvcmi}
\caption{}
\begin{algorithmic}
\STATE{Set $x=x_0$.}
\STATE {For $i$ ranging from $1$ to $M$, repeat the following steps:} 
\end{algorithmic}
  \begin{myenumerate}
  \item Pick $J$ at random from $1$ to $5$ with the probabilities
    $p_1$ through $p_5$ respectively.
  \item Pick a ``proposal'' point $y$ by following the subroutine
    specified in  {\bf Action J} (defined below).
  \item Calculate $\alpha=\min(1,\phi(y)/\phi(x))$.
  \item With probability $\alpha$, set $x=y$.
  \item Set $x_i=x$.
  \end{myenumerate}
\end{algorithm}

\begin{action}[H]
  \caption{Add or Delete a Random Coordinate}
  \begin{algorithmic}
\STATE{Set $y=x$.}
\STATE{Flip a fair coin.}
\IF {heads:}
\STATE{Generate $U$ uniformly at random on $[0,1]$.}
\STATE{Add $U$ to the end of $y$.}
\RETURN{$y$.}
\ELSIF{tails:} 
\STATE{Set $y=x$.} 
\IF {$y$ is empty:}
\RETURN{$y$.} 
\ELSE 
\STATE{Permute $y$ randomly.}
\STATE{Delete the last entry from $y$.}
\RETURN{$y$.}
\ENDIF
\ENDIF
  \end{algorithmic}
\end{action}

\begin{action}[H]
  \caption{Randomize a Coordinate}
  \begin{algorithmic}
\STATE{Set $y=x$.}
\IF{$y$ is empty:}
\RETURN{$y$.}
\ELSE
\STATE{Permute $y$ randomly.}
\STATE{Generate $U$ uniformly at random on $[0,1]$.}
\STATE{Replace the last coordinate of $y$ with $U$.}
\RETURN{$y$.}
\ENDIF
\end{algorithmic}
\end{action}

\begin{action}[H]
\caption {Shift a Coordinate}
\begin{algorithmic}
\STATE{Set $y=x$.} 
\IF{$y$ is empty:} 
\RETURN{$y$.}
\ELSE
\STATE{Permute $y$ randomly.}
\STATE{Set $U1=\text{the last element of $y$}$.}
\STATE{Set $U2=U1+\text{a random normal with mean 0 and standard deviation
    0.1}$.}
\IF{$U2$ is in $[0,1]$:}
\STATE{Replace the last coordinate of $y$ with $U2$.}
\ENDIF
\RETURN{$y$.}
\ENDIF
\end{algorithmic}
\end{action}

\begin{action}[H]
\caption {Randomize All Coordinates}
\begin{algorithmic}
\STATE{Set $y=x$.}
\IF{$y$ is empty:}
\RETURN{$y$.}
\ELSE
\FOREACH{coordinate of $y$}
\STATE{Generate $U$ uniformly at random on $[0,1]$.}
\STATE{Replace the current coordinate of $y$ with $U$.}
\ENDFOR
\RETURN{$y$.}
\ENDIF
\end{algorithmic}
\end{action}

\begin{action}[H]
\caption{Shift All of the Coordinates Very Slightly} 
\begin{algorithmic}
  \STATE{Set $y=x$.}
  \IF{$y$ is empty:}
  \RETURN{$y$.}
  \ELSE
  \FOREACH{coordinate of $y$}
  \STATE{Set $U1=\text{the current coordinate of $y$}$.}
  \STATE{Set $U2=U1+\text{a random normal with mean 0 and standard deviation
      0.01}$.}
  \IF{$U2$ is in $[0,1]$:}
  \STATE{Replace the current coordinate of $y$ with U2.}
  \ELSE
  \STATE{Continue to the next coordinate.}
  \ENDIF
  \ENDFOR
  \RETURN{$y$}.
  \ENDIF
\end{algorithmic}
\end{action}

If desired, the resulting $x_1, \dots, x_M$ can be randomly augmented
into a sequence $\theta_1, \dots, \theta_M$. To do so, simply generate
all the necessary $\Beta$ random variables in order to sample $\S$ from its
conditional distribution (c.f. \autoref{e:betas}).

\section{Posterior Mean Calculation}\label{s:calcmean}
There are many potential uses of for the sample of $\theta$ values
that can be approximately drawn from the posterior using the algorithm
in the previous section. For example, for each sampled $\theta$, a
plot of the corresponding $f_\theta$ can be made, and inspecting some
of these can give some idea about how confident to be about the shape
of the unknown regression function.

This section, though, concentrates on estimating the mean of these
functions. Call the resulting function the posterior mean,
$\widehat{f}$. It represents the posterior's best estimate (in an
$\cL^2$ sense) of the unknown regression function. More formally,
define the value of $\fhat$ at a point $u \in [0,1]$ as the expected
value of $f_\theta(u)$ under the posterior on $\theta$:

\begin{align}
\widehat{f}(u)
&=\int_{\Theta} f_\theta(x) \ppi(d\theta) \\
&=\int_{\cX} \widehat{s}_{J(u;x)} \phi(x)\tau(dx) \label{e:fhat}
\end{align}
Where $\shatj$ was defined by \autoref{e:shatj}, $\phi(x)$ was defined
by \autoref{e:phix}, and $\tau$ was defined shortly before $\phi$.

The algorithm from \autoref{s:algo} can be used to approximately
generate a sample $x_1, \dots, x_M$ from $\phi(x)\tau(dx)$ and then
estimate this integral by:
\begin{align}
\frac{1}{M}\sum_{i=1}^M \widehat{s}_{J(u;x)}
\end{align}

\section{Metropolis-Hastings Markov Chains on General Spaces}\label{s:mcmc}
Where does the algorithm described in \autoref{s:algo} come from?
This section addresses the MCMC approach and begins a description of
the larger framework within which algorithms like this one can be
derived and evaluated.

Generally, MCMC techniques suggest how to formulate algorithms (more
specifically Markov chains) that may be useful in order to sample a
stationary ergodic sequence that converges to a given stationary
distribution. This subject is very broad and active. For a review of
the main ideas, see Tierney~\cite{tierney:1994} or
Liu~\cite{liu:2001}. MCMC techniques have opened up to numerical
investigation a wide variety of Bayesian procedures, especially with
the advent of the Gibbs sampler, the Metropolis-Hastings
algorithm~\cite{metrop:1953, hastings:1970}, 
and its extension to the problem of ``model determination'' through
``reversible jump'' MCMC (Green~\cite{green:1995}).

In very general terms (following \cite{athr:doss:seth:1996}), the
Markov chain setup is as follows. Let $\pi$ be a probability distribution on a
measurable space $(\cX,\cB)$. Let $P$ be a transition probability
function on this space, that is, $P$ is a function on $\cX \times \cB$
such that, for each $x\in\cX$, $P(x,\cdot)$ is a probability measure
on $(\cX,\cB)$ and, for each $C \in \cB$, $P(\cdot,C)$ is a
measurable function on $(\cX,\cB)$.  The Markov chain $X_0,X_1,X_2,
\dots$ is generated as follows. We fix a starting point $X_0=x_0$,
generate an observation $X_1$ from $P(X_0,\cdot)$, generate an
observation $X_2$ from $P(X_1,\cdot)$, and so on.

$P$ will be constructed to obey detailed balance with respect to
$\pi$. Namely:
\begin{align}
\int_A P(x,B) \pi(dx) = \int_B P(x,A) \pi(dx) \text{~for every $A,B \in \cB$}
\end{align}
This condition is very convenient because although it will be easy to
construct chains that satisfy it, it is also powerful. In particular,
(by choosing $B=\cX$) it implies that $\pi$ is {\em an} invariant
measure for the Markov chain, that is,
\begin{align}\label{e:invariant}
\pi(A)=\int_\cX P(x,A) \pi(dx) \text{~for every $A \in \cB$}
\end{align}
The goal is to choose $P$ so that $\pi$ is the {\em unique} invariant
measure; and, moreover that the Markov chain will produce an ergodic
sequence of observations from $\pi$. For this goal, detailed balance
is a useful (although not necessary) ``first step.''

\newcommand{\tP}{\widetilde{P}}
\newcommand{\tp}{\widetilde{p}}
\newcommand{\tmu}{\widetilde{\mu}}

Suppose $\tP$ is some transition probability function which
(presumably) does not satisfy reversibility with respect to $\pi$. Let
$\tmu(dx,dy) \defined \pi(dx)\tP(x,dy)$. Suppose that $\tmu$ is
absolutely continuous with respect to some symmetric $\sigma$-finite measure
$\mu$. Specifically, suppose that $\mu$ is a measure on the measurable
space $(\cX \times \cX, \sigma(\cB\times\cB))$ that satisfies $\mu(A
\times B)=\mu(B \times A)$ for all $A,B \in \cB$. In simple cases, one
can choose $\mu$ to be a product measure; for example $\pi(dx)\pi(dy)$
or $\lambda(dx)\lambda(dy)$ where perhaps $\pi$ has a density with
respect to $\lambda$. If no such measure is readily available, one
may take $\mu(dx,dy) = \half \tmu(dx,dy)+ \half \tmu(dy,dx)$;
i.e. $\mu(A \times B)=\half \tmu(A \times B)+\half \tmu(B \times A)$
for $A,B \in \cB$. Let $\tp$ be a version of the Radon-Nikod\'ym
derivative of $\tmu$ with respect to $\mu$ so that
$\tmu(dx,dy)=\tp(x,y)\mu(dx,dy)$. If $\tp(x,y)$ happens to be
symmetric; i.e. $\tp(x,y)=\tp(y,x)$ for all $x,y \in \cX$, then $\tP$
already satisfies detailed balance with respect to $\pi$. To verify
this, compute that for every every $A,B \in \cB$:
\begin{align}
\int_A \tP(x,B) \pi(dx)
&= \int_{x \in A} \int_{y \in B} \pi(dx)P(x,dy) \\
&= \int_{(x,y) \in A \times B} \tp(x,y) \mu(dx,dy) \\
&= \int_{(x,y) \in A \times B} \tp(y,x) \mu(dx,dy) \\
&= \int_{(x,y) \in B \times A} \tp(x,y) \mu(dx,dy) \\
&= \int_{x \in B} \int_{y \in A} \pi(dx)P(x,dy) \\
&= \int_B \tP(x,A) \pi(dx)
\end{align}

When $\tp$ is not symmetric, there is no trouble in
constructing the closely related symmetric function
$p(x,y)=\min(\tp(x,y),\tp(y,x))$.  Does this suggest how to construct a
probability transition function $P$ based on $\tP$ that satisfies
detailed balance with respect to $\pi$? Yes, fortunately it
does. Define: 
\Align{
\alpha(x,y)=\begin{cases}
\min(1,\tp(y,x)/\tp(x,y)) & \text{if $\tp(x,y)>0$} \\ 1 & \text{if
$\tp(x,y)=0$}
\end{cases}
} Then (check) $p(x,y)=\alpha(x,y)\tp(x,y)$ for all $x,y \in \cX$.
This suggests defining $Q(x,dy)=\alpha(x,y)\tP(x,dy)$, so that
$\pi(dx)Q(x,dy)=\alpha(x,y)\tp(x,y)\mu(dx,dy)=p(x,y)\mu(dx,dy)$.  The
only problem is that $Q(x,\cdot)$ is not a probability (in general),
but a sub-probability. To account for the ``forgotten'' mass, set
$h(x)=1-Q(x,\cX)$ for all $x\in\cX$. And define
$P(x,dy)=Q(x,dy)+h(x)\delta_x(dy)$. Here, $\delta_x(\cdot)$ stands for
the Dirac measure at $x$. In expanded form, this defines: 
\Align{
P(x,dy)=\alpha(x,y)\tP(x,dy)+\left[ \int_\cX (1-\alpha(x,z))\tP(x,dz)
\right]\delta_x(dy) 
}

In conclusion, one can now easily verify that this defines a
probability transition function $P$ which satisfies detailed balance
with respect to $\pi$ and which is a
simple modification of $\tP$. Indeed, $P$ is simply a modification of $\tP$
that sometimes ``holds'' instead of taking the transition that $\tP$
proposes. Suppose that $Y=y$ is a particular value drawn from $\tP(x,
\cdot)$. That is, suppose that $\tP$ has ``proposed'' the transition
from $x$ to $y$. Then $P$ ``accepts'' this transition with probability
$\alpha(x,y)$, but holds with probability
$1-\alpha(x,y)$. Furthermore, because $\alpha(x,y)$ only depends on
the ratio $\tp(y,x)/\tp(x,y)$ it is sufficient to be able to compute
$\tp(x,y)$ up to an unknown constant factor.

Clearly, then, $P$ is a computationally simple modification of $\tP$;
the only caveat is that $\alpha(x,y)$ may often be very small or even
$0$ so that in the extreme degenerate case in which $\alpha
\identically 0$, $P$ is the Markov chain that always holds. Indeed,
this Markov chain is reversible with respect to any distribution, but
it certainly does not serve the larger goal of producing an ergodic
sequence of realizations from $\pi$. Similar problems can occur if
$\tP$ is not transitive or is otherwise unsuitable. For these reasons
the ergodicity conditions from \autoref{s:ergodicity} are needed.

\section{A Simple Markov Chain}\label{s:simplechain}
A simple example is in order to make these ideas more concrete. This
chain will not be as efficient (in practice) as the local-move Markov
chain developed in the next section.

In words, this will be the chain that stays fixed (holds) at its
current value $x \in \cX$ until a new value $y \cX$ drawn from the prior is
accepted; $y$ will always be accepted if $y$ makes the data more
likely (i.e. higher predictive probability under our model); otherwise
it is accepted with a probability reflecting the ratio of the
predictive probabilities. In this way, the chain readily walks
``uphill,'' but, with just the right probability (because of the
detailed balance condition that will be shown) it also walks ``downhill.''

Recall the notation from \autoref{s:mcmcsetup} that defined the
measurable space $(\cX, \cB)$ and denote the posterior distribution
on this space by $\ppi$. For convenience, recall that for any
$x\in\cX$, $x=(k,\v)$ and write $\kappa(x)=\kappa(k(x))$ and
$\rho(x)=\rho(\v(x))$ so that the prior $\pi$ on points $y \in \cX$
has density $\kappa(y)$ with respect to $\tau$. For any $x \in \cX$
and any $B \in \cB$, let $\tP$ be defined by $\tP(x,B)=\int_{y \in B}
\kappa(y)\tau(dy)$. That is, $\tP$ is the probability transition
function that (without reference to $x$) samples a new point $y\in\cX$
from the prior. Now expand $\ppi(dx)P(x,dy)=\phi(x)\kappa(y)
\tau(dx)\tau(dy)=(\rho(x)/c)\kappa(x)\kappa(y)\tau(dx)\tau(dy)$.
Conveniently then, this distribution has a density with respect to the
product measure $\kappa(x)\tau(dx)\kappa(y)\tau(dy)$. To agree with
the notation in the previous section, let $\mu$ denote this product
measure and let
$\tp(x,y)=\rho(x)/c$. Recall that
each of these terms is always positive. Because of this and a
convenient cancellation, the expression for $\alpha$ becomes
\[
\alpha(x,y)=\min(1,\rho(y)/\rho(x))
\]
Finally, as before, define $Q(x,dy)=\alpha(x,y)\tP(x,dy)$,
$h(x)=1-Q(x,\cX)$, and set $P(x,dy)=Q(x,dy)+h(x)\delta_x(dy)$.

\newcommand{\gdown}{g_\downarrow}
\newcommand{\gup}{g_\uparrow}
\newcommand{\lupx}{\lambda^x_\uparrow}
\newcommand{\ldownx}{\lambda^x_\downarrow}

\comment{
\newcommand{\gdown}{g_-}
\newcommand{\gup}{g_+}
\newcommand{\lupx}{\lambda^x_+}
\newcommand{\ldownx}{\lambda^x_-}
}
\section{A Local-Move Markov Chain}\label{s:localchain}
This section gives a formal definition of the Markov chain type
algorithm that was explained in \autoref{s:algo}. It then shows that
this chain satisfies detailed balance with respect to $\ppi$. To introduce this
more useful, but more complicated Markov chain on $(\cX, \cB)$, some
notation is needed. When $j \ge 2$, and $\v \in [0,1]^j$, write $\v_-$
for the vector $(\v_1, \dots, \v_{j-1}) \in [0,1]^{j-1}$ which leaves
off the last coordinate of $\v$. Consider some point $x=(k,\v) \in
\cX$. Let $\gdown(x)$ stand for the point in $\cX$ that is ``one level
down'' from $x$ with the last coordinate of $\v$ having been
removed. That is, when $k \ge 3$ and $\v=(\v_1, \dots,
\v_{k-1})\in[0,1]^{k-1}$, let $\gdown(x)= (k-1, \v_-)$. For $x=(2,
(\v_1)) \in \cX_2$, let $\gdown(x)=(1,\emp)$. For $x=(1,\emp) \in
\cX_1$, there is no further down to go, and so let $\gdown(x)=x$. Let
$\ldownx(\cdot)=\delta_{\gdown(x)}(\cdot)$ denote the Dirac measure at
$\gdown(x)$.

Similarly, let $\gup(x,u)$ stand for the point in $\cX$ that is ``one
level up'' from $x$, where the last coordinate is filled in with
$u$. That is, for $k \ge 2$ and $\v=(\v_1, \dots, \v_{k-1})$,
$\gup(x,u)=(k+1, (\v_1, \dots, \v_{k-1},u))$. For $x=(1,\emp)$, let
$\gup(x,u)=(2,(u))$. Let $\lupx(\cdot)$ denote the
distribution of $\gup(x,U)$ where $U$ is uniformly distributed on
$[0,1]$. 

Let $\lambda_{j}(\cdot)$ denote Lebesgue measure on $\bbR^j$ and for
$B \in \cB$ where $B \subset \cX_k$, let $\underline{B}=\{ \v \in [0,1]^{k-1}
~:~ (k,\v) \in B \}$.

\newcommand{\w}{\mathbf{w}}

To define the transition probability function $P$ on $(\cX, \cB)$
satisfying detailed balance with respect to $\ppi$, I first define
various transition probability functions $\tP_j(x,dy)$; then, for a
generic function $\alpha_j(x,y)$, define:
\begin{align}
Q_j(x,dy)=\alpha_j(x,dy)\tP_j(x,dy)
\end{align}
Next the $\alpha_j$'s are chosen so that for every $j$, $Q_j$ satisfies
the appropriate detailed balance formula:
\begin{align}
\int_{x \in A} & \pi(dx)  \int_{y \in B} Q_j(x,dy)
= \int_{x \in B} & \pi(dx)  \int_{y \in A} Q_j(x,dy)
\end{align}
It is easily verified then that $P_j=Q_j(x,dy)+[1-Q_j(x,\cX)]
\delta_x(dy)$, is a transition probability which satisfies detailed balance.
Finally, set:\(P(x,dy)=\sum_j p_j P_j(x,dy)\).

Generally the proposals I consider are ``symmetric'' and so $\alpha_j(x,y)$
will work out to be $\min(1,\phi(y)/\phi(x))$ in each case. An
exception is in \autoref{e:balance4}, but when the proposal density is
simple, (as it is in the case of interest) it also reduces to the
previous case.

To begin, set $\tP_1(x,dy)=\half \lupx(dy) + \half \ldownx(dy)$. This
transition probability function represents the chain that adds or
deletes coordinates from the vector $\v$ randomly.

To choose $\alpha_1$, first compute the left and right hand sides of
the detailed balance equation for $Q_1$. Let $A, B \in \cB$ with $A
\subset \cX_j$ and $B \subset \cX_k$. There are three cases of
interest for $j$ and $k$: (1) $k=j+1, j\ge 2$, (2) $j=1,k=2$, (3)
$j=1,k=1$. These account for all the possibilities because if $j<k$,
simply replace the roles of $j$ and $k$; if $k>j+1$ or $j=k\neq 1$
both sides will evaluate to $0$. If $\alpha_1$ can be chosen to set the
left and right sides equal for every such case, detailed balance is
proven because general $A,B \in \cB$ can be decomposed into these
component subsets.

Suppose that $k=j+1$, and calculate:
\begin{align}
\int_{x \in A} & \pi(dx)  \int_{y \in B} Q_1(x,dy) \\
&= \int_{x \in A} \int_{y \in B} \phi(x)\alpha_1(x,y) \tau(dx) [\half
  \lupx(dy) + \half \ldownx(dy)] \\
&= \int_{x \in A} \int_{y \in B}
   \half \phi(x)\alpha_1(x,y) \tau(dx) \lupx(dy) \\
&= \int_{\v \in \underline{A}} 
   \int_{u \in C(\v;B)}
   \half \kappa(j)\rho_j(\v) \alpha_1((j,\v),\gup((j,\v),u))
   \lambda_{j-1}(d\v) du \\
&= \int_{\w \in \underline{D}} 
   \half \kappa(j)\rho_j(\w_-) \alpha_1((j,\w_-),\gup((j,\w_-),u))
   \lambda_{j}(d\w) \\
&= \int_{x \in D} 
   \half \phi(\gdown(x)) \alpha_1(\gdown(x),x) \tau(dx)
\end{align}
Where:
\begin{align}
C(\v;B) = {\{u \in [0,1]~:~ \gup( (j,\v), u) \in B \}}
\end{align}
\begin{align}
\begin{split}
\underline{D} &= \{ \w \in [0,1]^j ~:~ \w_- \in \underline{A}, w_j \in C(\w_-;B)
 \}
\end{split} \\
&= \{ \w \in \underline{B} ~:~ \w_- \in \underline{A} \}
\end{align}
\begin{align}
D &= \{ x=(k,\w) \in C_k ~:~ \w \in \underline{D} \} \\
  &= \{ x \in B ~:~ \gdown(x)\in A\}
\end{align}

Similarly, compute:
\begin{align}
\int_{x \in B} & \pi(dx)  \int_{y \in A} Q_1(x,dy) \\
&= \int_{x \in B} \int_{y \in A} \phi(x)\alpha_1(x,y) \tau(dx) [\half
  \lupx(dy) + \half \ldownx(dy)] \\
&= \int_{x \in B} \int_{y \in A}
   \half \phi(x)\alpha_1(x,y) \tau(dx) \ldownx(dy) \\
&= \int_{x \in B} 
   \half \phi(x)\alpha_1(x,\gdown(x)) \tau(dx) \1{\gdown(x) \in A} \\
&= \int_{x \in D} 
   \half \phi(x)\alpha_1(x,\gdown(x)) \tau(dx)
\end{align}

Suppose that $j=1,k=2$, so that $A=\{x_0\}$ or $A=\empty$
where $x_0=(1,\emp)$ and any $x \in B$ is of the form $(2,(\v_1))$. Then: 
\begin{align}
\int_{x \in A} & \pi(dx)  \int_{y \in B} Q_1(x,dy) \\
&= \int_{x \in A} \int_{y \in B} \phi(x)\alpha_1(x,y) \tau(dx) [\half
  \lupx(dy) + \half \ldownx(dy)] \\
&= \int_{x \in A} \int_{y \in B}
   \half \phi(x)\alpha_1(x,y) \tau(dx) \lupx(dy) \\
&= \int_{y \in B}
   \half \1{x_0\in A}\phi(x_0)\alpha_1(x_0,y) {\lambda^{x_0}_\uparrow}(dy) \\
&= \int_{\v_1 \in \underline{B}}
   \half \1{x_0\in A}\phi(x_0)\alpha_1(x_0,(2,\v_1)) d\v_1 \\
&= \int_{x \in D}
   \half \phi(x_0)\alpha_1(x_0,x) \tau(dx)
\end{align}

Similarly, compute:
\begin{align}
\int_{x \in B} & \pi(dx)  \int_{y \in A} Q_1(x,dy) \\
&= \int_{x \in B} \int_{y \in A} \phi(x)\alpha_1(x,y) \tau(dx) [\half
  \lupx(dy) + \half \ldownx(dy)] \\
&= \int_{x \in B} \int_{y \in A}
   \half \phi(x)\alpha_1(x,y) \tau(dx) \ldownx(dy) \\
&= \int_{x \in B} 
   \half \phi(x)\alpha_1(x,\gdown(x)) \tau(dx) \1{\gdown(x) \in A} \\
&= \int_{x \in D} 
   \half \phi(x)\alpha_1(x,x_0) \tau(dx)
\end{align}

The only remaining case to compute is where $j=k=1$ and the only case
of interest here is the one in which $A=B=\{x_0\}$. Even this is
trivial, because the left and right sides of this symmetric case must
surely match.

Recalling that $\phi(x)>0$ for all $x \in \cX$, for $x,y
\in \cX$ detailed balance is achieved upon defining:
\begin{align}
\alpha_1(x,y)=\min(\phi(y)/\phi(x),1)
\end{align}

\newcommand{\lchj}{\lambda^x_{\cdot j}}
\newcommand{\lch}{\lambda^x_{\cdot}}

To define $\tP_2$ additional notation is needed. For $k \geq 2, 1 \le j \le
k-1$ let $\lchj(\cdot)$ denote the distribution on $(\cX,\cB)$ under
which the $j$'th coordinate of the $\v$ in $x=(k,\v)$ is replaced with
a uniform random variable. For $k \ge 2$ and $x \in \cX_k$, let
$\lch(dy)=\frac{1}{k-1}\sum_{j=1}^{k-1} \lchj(dy)$, i.e. the
probability distribution that picks a coordinate of the $\v$ in $x$ randomly
and then changes it to a uniformly random value. For $x=x_0$, let
$\lch(dy)=\delta_{x_0}(dy)$.

Now define $\tP_2(x,dy)=\lch(dy)$. More generally, define:
\begin{align}
\tP_3(x,dy)=\Psi(x,y)\lch(dy)+h_\Psi(x)\delta_x(dy)
\end{align}
Where $\Psi$ is any (measurable) non-negative function on $\cX \times
\cX$ satisfying for all $x,y \in \cX$: (1) $\Psi(x,y)=\Psi(y,x)$, (2)
$\int_{\cX_k} \Psi(x,y)\lch(dy)+h_\Psi(x)=1$ for some non-negative
$h_\Psi(x)$.  For example, if $x=(k,\v_x)$, $y=(k,\v_y)$, and
$\phi_\sigma$ represents the univariate normal density with standard
deviation $\sigma>0$, take $\Psi(x,y)=\phi_\sigma(||\v_x-\v_y||_2)$ and
$h_\Psi(x)=1-\int_{\cX_k} \Psi(x,y)\lch(dy)$. Fortunately, there will be
no need to compute $h_\Psi$; it is enough to note that it is
non-negative.

To check the detailed balance for $Q_3$ it is sufficient to consider
$A,B \in \cB$ where $A,B \subset \cX_k$. 

\begin{align}
\int_{x \in A} & \pi(dx)  \int_{y \in B} Q_3(x,dy) \\
&= \int_{x \in A} \int_{y \in B} \phi(x)\alpha_3(x,y) \tau(dx)
  [\Psi(x,y)\lch(dy)+h_\Psi(x)\delta_x(dy)] \\
\begin{split}
&= \int_{x \in A} \int_{y \in B} \phi(x)\alpha_3(x,y)\Psi(x,y)
   \tau(dx) \lch(dy) \\
&\qquad\qquad +\int_{x \in A} \phi(x)\alpha_3(x,x) h_\Psi(x)\1{x \in B} \tau(dx)
\end{split}\\
\begin{split}
&= \int_{x \in A} \int_{y \in B} \phi(x)\alpha_3(x,y)\Psi(x,y)
   \tau(dx) \lch(dy) \\
&\qquad\qquad +\int_{x \in A \intersect B}
   \phi(x)\alpha_3(x,x) h_\Psi(x) \tau(dx)
\end{split}
\end{align}

\newcommand{\wx}{\widetilde{X}}
\newcommand{\wy}{\widetilde{Y}}

This second term in the latter expression does not change if we
replace $A,B$ with $B,A$ and consequently needs no further
consideration. To expand the first term, recall that in this situation
$\lch$ just picks one of the $k-1$ coordinates and randomizes it.  Let
$S_j^{k}(\w)$ stand for the modified version of $\w$ in which the $j$'th
and $k$'th coordinates of $\w$ are swapped. 
In suggestive notation, let:
\begin{align}
&\wx=\wx(\w)=(k,w_-) \\
&\wy=\wy(\w)=\wx(S_j^k(\w)) \\
&D_j= \{\w \in [0,1]^k ~:~ \wx(\w) \in A,~ \wy(\w) \in B \}\\
&D'_j= \{\w \in [0,1]^k ~:~ \wx(\w) \in B,~ \wy(\w) \in A \}
\end{align}
Notice that $\w \in D_j \iff S_j^k(\w) \in D'_j$. 

Then, suppressing the $\w$ dependence from $\wx$ and $\wy$, the first
term expands into:
\begin{align}
\frac{1}{k-1} \sum_{j=1}^{k-1} \label{e:BA}
\int_{\w \in D_j}& 
   \phi(\wx)
   \alpha_3(\wx,\wy)
   \Psi(\wx,\wy)
   \lambda_k(d\w) \\
&=
\frac{1}{k-1} \sum_{j=1}^{k-1}
   \int_{\w \in D_j'}
   \phi(\wy)
   \alpha_3(\wy,\wx)
   \Psi(\wy,\wx)
   \lambda_k(d\w)\\
&=
\frac{1}{k-1} \sum_{j=1}^{k-1} \label{e:AB}
\int_{\w \in D_j'} 
   \phi(\wy)
   \alpha_3(\wy,\wx)
   \Psi(\wx,\wy)
   \lambda_k(d\w).
\end{align}

In summary then, (using \autoref{e:AB})
\begin{align}
\int_{x \in A} & \pi(dx)  \int_{y \in B} Q_3(x,dy) \\
&=
\frac{1}{k-1} \sum_{j=1}^{k-1}
\int_{\w \in D_j'}
   \phi(\wy)
   \alpha_3(\wy,\wx)
   \Psi(\wx,\wy)
   \lambda_k(d\w) + \mathit{const}
\end{align}
While, using \autoref{e:BA} with the roles of $A$ and $B$ interchanged
so that ``$D_j$''=$D'_j$:
\begin{align}
\int_{x \in B} & \pi(dx)  \int_{y \in A} Q_3(x,dy) \\
&=
\frac{1}{k-1} \sum_{j=1}^{k-1}
\int_{\w \in D'_j}
   \phi(\wx)
   \alpha_3(\wx,\wy)
   \Psi(\wx,\wy)
   \lambda_k(d\w) + \mathit{const}
\end{align}

Again we achieve balance using:
\begin{align}
\alpha_3(x,y)=\min(\phi(y)/\phi(x),1)
\end{align}

Define $\tP_4(x,dy)=\xi_x(y)\tau(dy)$. Where for any $x \in \cX$,
$\xi_x(\cdot)$ is a density with respect to $\tau$.
Accordingly:
\Align{
\ppi(dx)\tP_4(x,dy) &= \phi(x)\xi_x(y)\tau(dx)\tau(dy) 
=\tp(x,y)\tau(dx)\tau(dy)
}
By the arguments in \autoref{s:mcmc} the corresponding transition
probability function $P_4$ satisfies detailed balance with respect to
$\ppi$ provided that:
\begin{align}\label{e:balance4}
\alpha_4(x,y)=
\begin{cases}
\min(1,[\phi(y)\xi_y(x)]/[\phi(x)\xi_x(y)] & \text{if $\phi(x)\xi_x(y)>0$} \\
1 & \text{if $\phi(x)\xi_x(y)=0$}
\end{cases}
\end{align}

For example, for $x \in \cX_k$, choose $\xi_x(y)=\1{y \in \cX_k}$.
And for this $\xi_x$, the familiar choice
$\alpha_4(x,y)=\min(1,\phi(y)/\phi(x))$ works equally well.

Finally, define $\tP_5(x,dy)$ for $x \in \cX_k$ as the composition of
the $k-1$ separate transition probability functions $\tP_5^j$ applied
sequentially from $j=1$ to $j=k-1$. Each $\tP_5^j$ is just
intended to move the $j$'th coordinate of $\v_x$ by a small amount (or
hold). In effect, then, $\tP_5$ moves all of the coordinates of $x$
randomly, but technically speaking some subset of them may hold on any
given step. For $x \in \cX_k$ and $1 \leq j \leq k-1$ set
$\tP_5^j(x,dy)=\Psi(x,y)\lchj(dy)+h_\Psi(x)\delta_x(dy)$ for some $\Psi$
and $h_\Psi$ satisfying the same constraints as made when defining
$\tP_3$. I omit the verification that by choosing
$\alpha_5=\min(1,\phi(y)/\phi(x))$, as usual, $P_5$ satisfies detailed
balance.

This concludes a consideration of each chain $P_1$ through $P_5$. Each
was shown to satisfy detailed balance with respect to $\ppi$.
Consequently, their mixture $P$ also satisfies detailed balance with
respect to $\ppi$.

As an aside, notice that one may compose any of the $P_j$
with a random permutation since $\phi$ itself is invariant with
respect to permuted $\v$. Doing so will allow $P_1$ to add and (by
pre-composing) delete coordinates in arbitrary locations.

\section{Markov Chain Convergence Theory}\label{s:ergodicity}

Early results on the ergodicity of Markov chains on general state
spaces used a condition known as the Doeblin condition. It implies
that there exists an invariant probability measure to which the Markov
chain converges at a geometric rate, from any starting point.\footnote{This section
  reviews some material given in a paper by Athreya, Doss, and
  Sethuraman~\cite{athr:doss:seth:1996}}

\begin{thm}[Doob (1953)~\cite{doob:1953}]
Suppose that the Markov chain on the measure space $(\cX, \cB)$ generated
by probability transition function $P$ satisfies the Doeblin
condition that there is a probability measure $\lambda$ on
$(\cX,\cB)$, an integer $k$, and an $\epsilon>0$ such that:
\[
P^k(x,C) \geq \epsilon \lambda(C) \quad\text{for all $x \in \cX$ and
  all $C \in \cB$}
\]
Then there exists a unique invariant probability measure $\pi$ such
that for all $n$,
\[
\sup_{C \in \cB} \left| P^n(x,C)-\pi(C)\right| \leq (1-\epsilon)^{(n/k)-1}
\qquad\text{for all $x \in \cX$}
\]
\end{thm}

An easy corollary is that if $P$ satisfies the conditions of the
theorem, and was already known to be reversible with respect to some
specific distribution, then that same distribution must be the unique
stationary distribution $\pi$. The Doeblin condition is quite strong
and rarely holds in applications.

Athreya, Doss, and Sethuraman~\cite{athr:doss:seth:1996} prove an ergodicity
result for general state spaces whose conditions hold much more
broadly and remain reasonably easy to check. An abbreviated version is
given below

\begin{thm}[Athreya, Doss, and Sethuraman (1996)]\label{t:ath}
Suppose that the Markov chain $\{X_n\}$ with transition function
$P(x,C)$ has an invariant probability measure $\pi$, that is
\autoref{e:invariant} holds. Suppose that there is a set $A\in\cB$, a
probability measure $\lambda$ with $\lambda(A)=1$, a constant
$\epsilon>0$, and an integer $n_0 \ge 1$ such that:
\Align{\label{c:canreach}
\pi(\{ x: P_x(X_m \in A \text{~for some~} m\ge 1)>0 \})=1
}
and
\Align{
P^{n_0}(x,\cdot) \ge \epsilon \lambda(\cdot) \quad\text{for each $x\in
  A$}
}
Further suppose that either $n_0=1$ or that
\Align{
\mathrm{g.c.d.}\{m ~:~\text{there is an $\epsilon_m>0$ such that
  $P^m(x,\cdot) \ge \epsilon_m\lambda(\cdot)$ for each $x\in A$}\}=1
}
Then there is a set $D \in \cB$ such that
\Align{
\pi(D)=1\quad\text{and}\quad \sup_{C\in\cB} |P^n(x,C)-\pi(C)|
  \tendsto 0~~\text{for each $x \in D$.}
}
Let $f(x)$ be a measurable function on $(\cX,\cB)$ such that
$\int\pi(dy)|f(y)|<\infty$. Then
\Align{
P_x\left\{
\frac{1}{n} \sum_{j=1}^n f(X_j) \tendsto \int \pi(dy)f(y)
\right\}=1 \quad \text{for $[\pi]$-almost all $x$}
}
\end{thm}

\comment{
A simple corollary is:
\begin{cor}\label{t:cor:ath}
Suppose that $P$ is a (finite) mixture of probability transition
function on $(\cX,\cB)$. That is $P(x,dy)=\sum_j p_j\P_j(x,dy)$, for
some $p_j \ge 0$ with $\sum_j p_j=1$. Further suppose that for every
$j$, $\pi$ is an invariant measure for $P_j$. And that for some $j^*$,
with $p_{j^*}>0$, $P_{j^*}$ satisfies the conditions of Theorem~\ref{t:ath},
then $P$ itself satisfies these conditions.
\end{cor}
\begin{proof}
Let $A$, $\epsilon$, and $n_0$ that satisfy the conditions of
Theorem~\ref{t:ath} for $P_{j^*}$. Consider some $x$ for which
$P_{j^*}(E_x)>0$ where $E_x$ is the event $\{X_m \in A \text{~for some
  $m \ge 1$}\}$. Let $T$ be the (random) time until $X_m \in A$
\end{proof}
}

\section{Convergence Results}\label{s:ergoproof}
A short proof suffices to show that theorem~\ref{t:ath} applies to
the local-move Markov chain. It is assumed that mixing probability
$p_1>0$. The proof takes advantage of the atom $x_0=\bigl(1,\emp\bigr)$.

Let $P$ denote the local-move Markov chain from
\autoref{s:localchain}. It was already verified there that $P$
satisfies detailed balance with respect to $\ppi$, and it follows that
$P$ has invariant probability measure $\ppi$. Let $A=\cX_1$, the
singleton set $\{x_0\}$. Let $\lambda(\cdot)$ be counting measure on
$\cX_1$. Let $n_0=1$.  Set $\epsilon=P(x_0, \cX_1)$, i.e. the chance
of holding at this atom. This quantity is positive (because
$P(x_0, \cX_1)\ge p_1 P_1(x_0,\{x_0\}) \ge p_1\tP_1(x_0,\{x_0\})>0$).

This verifies all of the conditions of theorem~\ref{t:ath} except
condition~\ref{c:canreach}. It suffices to show that for any $k \in
\bbN$ and any $x \in \cX_k$, $P_x(X_m \in A~\text{for some $m \ge
1$})>0$. It has already been shown that for any $x \in \cX_1$,
$P_x(X_m \in A~\text{for some $m \ge 1$})>0$, since $\cX_1$ is the
singleton set considered earlier and because this quantity is greater
than $P(x_0,A)$, which was positive. Suppose, for induction, that for
any $k \le k_0$ and any $x \in \cX_k$, $P_x(X_m \in A~\text{for some
$m \ge 1$})>0$. Consider any $x \in
\cX_{k_0+1}$.
\begin{align}
P(x,\cX_{k_0})&=p_1P_1(x,\cX_{k_0})\\
&=p_1\int_{y\in \cX_{k_0}} \left[\half \alpha_1(x,y)\lupx(dy) + \half \alpha_1(x,y)\ldownx(dy)+h_1(x)\delta_x(dy)\right]\\
&=\half p_1\int_{y\in \cX_{k_0}} \alpha_1(x,y) \ldownx(dy)\\
&=\half p_1\int_{y\in \cX_{k_0}} \alpha_1(x,y) \delta_{\gdown(x)}(dy)\\
&=\half p_1 \alpha_1(x,\gdown(x)) \\
&=\half p_1 \min(\phi(\gdown(x))/\phi(x),1) 
\end{align}
Now, both $\phi(x)$ and $\phi(\gdown(x))$ are positive and finite, so
their ratio is as well, and so $P(x, \cX_{k_0})>0$. This proves that
for any $x \in \cX_{k_0+1}$,
\begin{align}
P_x(X_m \in A~~\text{for some
$m \ge 1$})>0
\end{align}
All the conditions are now verified.

The same argument goes through without modification for the simple
Markov chain from \autoref{s:simplechain}.

In summary the preceding sections have proven the following theorem:
\begin{thm}
Let $(\cX,\cB)$ be the measurable space defined in
\autoref{s:mcmcsetup}. Let $P$ be the probability transistion
function on this space that is defined as a mixture of the component
probability transistion functions $P_1$, through $P_5$, explained in
\autoref{s:localchain}, with mixture weights $p_1$, through $p_5$
positive and summing to 1. (This Markov chain is a formal version of the
Algorithm given in \autoref{s:algo}.)  Then, $P$ satisifies detailed
balance with respect to the distribution $\ppi$ defined by
\autoref{e:phix}. Furthermore, $P$ satisfies the conditions of
Theorem~\ref{t:ath} and therefore, $\ppi$ is the unique invariant
distribution of $P$. Indeed, there is a set $D \in \cB$ such that
\Align{
\ppi(D)=1\quad\text{and}\quad \sup_{C\in\cB} |P^n(x,C)-\ppi(C)|
  \tendsto 0~~\text{for each $x \in D$.}
}
Let $f(x)$ be a measurable function on $(\cX,\cB)$ such that
$\int\ppi(dy)|f(y)|<\infty$. Then
\Align{
P_x\left\{
\frac{1}{n} \sum_{j=1}^n f(X_j) \tendsto \int \ppi(dy)f(y)
\right\}=1 \quad \text{for $[\ppi]$-almost all $x$}
}
\end{thm}


%% file: examples.tex
\chapter{Examples}\label{c:examples}

This chapter describes the results of a variety of simulation
experiments that I have conducted to better understand and evaluate
the performance of the posterior mean estimates based on the prior
$\pi$. In all of these experiments, except where specifically noted,
the prior $\kappa$ on the number of steps $K$ is taken to be
$\Geometric(\half)$. Other $\Geometric$ priors are considered in
\autoref{s:geomparam} and \autoref{s:largedata}. A few examples in
\autoref{c:discuss} consider $\Poisson$ priors. There, they are
discussed in relation to the convergence theory proven in
\autoref{c:proof}. Most of this chapter, however, concerns evaluating
how efficient the $\Geometric(\half)$ posterior mean estimate is by
comparing it with a wide variety of competing estimation
procedures. The first section establishes the standard format for
these experiments. It compares the posterior mean estimate with CART
and bagged CART estimates and interprets the results. Other methods
are considered and compared in \autoref{s:otherpop}, namely: a Lasso example
that is connected with bagging, three smoothers, and
some wavelet-based estimates. Finally, the dyadic Diaconis and Freedman binary
regression prior is compared in
\autoref{s:dyadcomp}. Section~\ref{s:predprob} takes a step back to
analyze the interaction between the data and the model by inspecting
how the predictive probability changes as a function of where splits
are placed. The final sections experiment with smaller and larger data
set size, and also evaluate the performance of the posterior mean on a more challenging
regression problem. 

\section{Comparison with CART and Bag\-ged CART}
{ \newcommand{\cp}{{\em cp}} \newcommand{\X}{\cX} \newcommand{\x}{x}
Since $\pi$ puts a prior on piecewise constant regression functions,
it is natural to ask how its performance compares with conventional
estimators that employ piecewise constant approximations.  Of
particular interest is the {\em Classification and Regression Tree}
(CART) algorithm~\cite{bfos:1984}, and the closely related bagged CART
algorithm. These methods are briefly reviewed in the next two
sections. The impatient reader should skip ahead to
\autoref{s:cartexp}, where the methods are compared with the
posterior mean estimate on simulated data sets.  I have not
``filtered'' the experimental data at all: these are the originally
simulated data sets in their original order. Technically, there has
been a certain amount of filtering in the results because I did try
using a variety of settings for the CART and bagged CART methods that
I do not discuss. The choices that are presented are among the more
standard and better performing possibilities. As far as the posterior
mean results, these are not filtered at all, except that arguably I
would not have a thesis if the results were not interesting.

\subsection{CART Review}
The CART algorithm prescribes how to select a ``tree'' that represents
a good estimate of the unknown regression function (or classification
rule). The ``tree'' terminology connotes the fact that the covariate
space $\X$ is recursively partitioned with each piece assigned its own
regression value: this recursive partition can be naturally associated
(by inclusion)
with a graph-theoretic tree. For completeness, a brief description of
the CART algorithm is in order. The reader should bear in mind that
CART and related algorithms have been in use for many years now and so
there are a variety of possible tweaks and alternatives that I do not
discuss. The CART algorithm also has important advantages that the
following discussion will not address. Its regression estimate, for
example, is unaffected if individual coordinates of the covariate
vector are rescaled, so that the estimate does not depend upon the units of
measurement. It is capable of dealing with large data sets because of
its efficient implementation. It is also very easy to interpret
(although this is a risky business since the estimates can sometimes
change dramatically with new data).

Suppose for simplicity that the covariate space $\X$ is $\bbR^d$
(specifically, the case $\X=\bbR^1$ is of interest at present). In its
basic form, the CART algorithm uses coordinate-aligned splits. That
is, if a certain subset $A$ of $\X$ is being partitioned, it will be
partitioned into the two subsets $\{\x \in A ~:~ \x_j \le c\}$ and
$\{\x \in A ~:~ \x_j > c\}$ for some choice of $j \in \{1, \dots, d\}$
and $c \in \bbR$. CART proceeds to construct a partition of $\X$ in a
greedy manner. That is, it begins by finding the binary partition of
$\X$ that maximizes a certain splitting criterion and (proceeding
recursively) all subsequent splits are subordinate to this
one. Ultimately, the splitting criterion is chosen by the user, and I
have chosen to use the ANOVA criterion. To explain this criterion
consider that at any given stage of partitioning there is a certain
class of possible real-valued functions that are constant on each
partition element. Among this class, the function that minimizes the
mean of squared residuals to the response data $y_1, \dots, y_n$ (MSE)
is clearly the one whose value on any given element of the partition
is the mean of the response values whose covariates ``hit'' this
element. The split that is considered best is the split that results
in the greatest possible reduction in this measure of residual
error. Having chosen a split, the CART algorithm continues,
recursively, to split the resulting subsets. Implicitly, it is
building up a ``tree'' of subsets at each stage with $\X$ forming the
root. The recursion terminates whenever there is only a single data
point in the current partition element. Call the resulting binary tree
of subsets the ``full'' tree. CART then proceeds to ``prune'' this
tree. This operation depends critically upon a complexity parameter
\cp, which must itself be chosen. Typically \cp\ is chosen by
cross-validation with the restriction that it not be smaller than some
user specified value. For these experiments I choose \cp\ using
10-fold cross-validation and the one-standard-deviation selection
rule. This means that if the best achievable cross-validated measure
of error (searching over all possible values of \cp) is {\em xerr} and
the sample standard deviation of {\em xerr} is {\em xstd}, the
selected value of \cp\ will be the largest value whose {\em xerr} does
not exceed {\em xerr}$+${\em xstd}.  To prune the full tree using
parameter \cp, each pair of leaves is considered in turn; if the pair
does not improve the splitting criterion by at least the value \cp, it
is removed. Ultimately, every pair of leaves might be removed,
resulting in the tree consisting only of the root node. Finally, the
pruned tree corresponds to the estimated regression function. It is
constant on each element of the (pruned) partition and its value on a
given element is the mean of the response values there.

}

\comment{
C4.5 Tutorial
P. Winston, 1992.
T. Mitchell, "Decision Tree Learning", in T. Mitchell, Machine Learning, The McGraw-Hill Companies, Inc., 1997, pp. 52-78.

P. Winston, "Learning by Building Identification Trees", in
P. Winston, Artificial Intelligence, Addison-Wesley Publishing
Company, 1992, pp. 423-442.

Breiman, L., Friedman, J., Olshen, R., and Stone, C. (1984)
Classification and Regression Trees", Wadsworth.  }

\subsection{Bagging Review and Discussion}\label{s:bagging}
A comparison with the bagging procedure \cite{breiman96bagging} is also
relevant. Bagging~\cite{breiman96bagging} is a meta-algorithm that can
(hypothetically) be applied to any existing classification or
regression technique in an effort to improve them. This thesis focuses
on bagged CART. Essentially, the
bagging idea is just to take many bootstrap resamples of the data
set, apply some existing technique to each resampled data set, and
then take an average of all of the resulting regression estimates. For
completeness, to form a bootstrap resample of a dataset with $n$ items
$z_i=(x_i,y_i)$: (1) Independently, choose $n$ integers $N_1, \dots,
N_n$ uniformly at random from $1$ to $n$. (2) Form the new data set:
$z_{N_1},z_{N_2}, \dots, z_{N_n}$.

It is sometimes stated~\cite{HTF:2001} that bagging is approximately a
non-parametric Bayesian procedure, but I think that this is a
misleading claim.  Bagging, or more specifically bootstrapping,
approximates the behavior of a Bayesian who has a (limiting) Dirichlet
prior (as in Rubin's Bayesian bootstrap~\cite{rubin}). This is not really a
prior for several reasons.  Besides the fact that this limiting prior
is improper, the ``prior'' depends on the data. This is not just in a
partial sense (such as the prior I use in \autoref{c:voronoi} in which
the ``prior'' depends on the covariates but not the response) or even
in the manner of empirical Bayes procedures which choose some
parameters of the ``prior'' by looking at the data. Indeed, this prior
specifies the law of all possible data {\em relative} to the empirical
distribution. This might make sense if the data were multinomial, but
for data on a continuous space it is quite problematic. If taken
literally it specifies that all future data will consist of elements
drawn from the current data set. Notably, for classification, this
means that unless the data set happens to contain a head and a tail
for each case, then the predictive distribution that
the bootstrap prior corresponds to excludes the possibility that the
missing flip will ever occur. Indeed, this bootstrap prior can never
directly say anything about future datapoints whose covariates are
new. Furthermore, there is no meaningful model involved. For example,
the prior $\pi$ implicitly builds in information that says that
if a certain region has more heads than tails, then future points in
and around this region probably will as well. The bootstrap prior says
nothing like that explicitly; if it says that {\em in effect} it is
only because of the fitting method that is forced upon it.

Finally, there is nothing Bayesian about using CART, so how can bagged
CART be a non-parametric Bayesian procedure? Why would a Bayesian who
believed in a bootstrap prior use CART or neural nets or whatever when
he could easily compute his own (very bizarre) posterior and get
conclusions directly? Arguably, he might do so in order to get new
information to guide his decision because he has observed that CART
(say) has worked well on other problems so it probably will work on
this one as well. In this sense, bagging could be said to model the
behavior of a Bayesian who had a limiting Dirichlet prior (that
magically was supported on the data set itself), who then computed the
posterior of his ``prior,'' and who then goes to seek the opinion of
an ``expert.'' Cleverly, he does so, not only for the dataset actually
received, but for a multitude of datasets of size $n$ that are about
equally likely under his posterior. In this way, he finds out what
CART would think in a variety of situations that he subjectively
considers as possibilities. So far so good. Now he ought to weight
these opinions according to how credible they seem in light of the
data and his prior opinion about when CART works and when it
doesn't. Instead, he now effectively forgets his (Dirichlet) posterior
and puts a flat prior on the CART regression results themselves. With
his ``{\em newly found prior},'' he calculates the decision that
minimizes squared error loss, the mean. So, in summary, bagging
approximates the behavior of a forgetful Bayesian who looks at the
data first, {\em then} formulates his ``prior'' and posterior, then
{\em ignores} them, except to ask CART what {\em it's opinion} would
be in the cases that he thinks are likely, then promptly forgets the
data altogether as well as any priors or posteriors of his own that he
might have held (recently), so he makes up a new uniform {\em prior}
on all the results that he got back from CART; finally he computes the
mean according to his latest prior and reports his ``findings.''

In any case, however dubious as a ``Bayesian'' procedure, bagging
works. One could say that this is because it reduces modeling bias or
because it eliminates certain instabilities in CART: both of these
arguments make perfect sense to Bayesians and frequentists alike. It
is clear, however, that it's not always a good idea, especially when
the procedure already has low ``instability.'' In this case, bagging
mostly adds noise and reduces the effective size of the dataset
somewhat.

I also found one example where bagging was
disastrous. Following the bagging procedure strictly, I fed bootstrap
resampled data sets into CART, but the result was a mess of
indecipherable noise with splits everywhere. This was true even though
the CART procedure gave a reasonable estimate on the original data
set. Why? Because in the CART step I used cross-validation and this is
problematic because CART is working on a bootstrapped sample. Certain
repeated data points wind up in both the test and training sets. As a
consequence the CART algorithm has less data that it thinks and also
thinks that is not over-fitting when it uses a model with too many
splits; the pruning procedures became completely ineffective. 

Some authors argue that this problem is easy to avoid by simply not
bothering to cross-validate within CART, and instead using the full
trees. This may be true in some cases, or if one implicitly uses an
effective default pruning rule (and not actually full trees), but it
failed on my example. To correct the cross validation, some authors
``tell'' CART which points are repeated so that the whole case is
either left in or left out. Instead, for my experiments, I simply used
a hand-tuned complexity parameter lower bound. For the disastrous
experiment, the lower bound was set to $0$; when increased to $0.005$
the estimates were still quite rough, but usable. For the experiments
that I present in the next section, I set to the lower bound to
$0.01$, which (admittedly) is the default for the RPART implementation
of CART. Still, unless this default is universally good, this leaves a
tuning parameter to be set, and I shudder to think about the bizarre
computations involved in choosing it by cross-validating bagged CART.

For the experiments involving bagged CART, I used $100$ bootstrap
resamples. This is a larger number of bootstrap resamples than is
generally considered necessary for bagging. By informal
Rao-Blackwell-type reasoning one would think that this only serves to
reduce Monte Carlo error.

\subsection{Comparative Simulation Experiment}\label{s:cartexp}
\bigfig{unif_vs_cart1}{90}{Comparative Simulation Experiment: Run
  1}{{\em Simulation Experiment: Run 1}\qquad The CART, Bagged CART, and
  posterior mean estimates are compared with the true regression curve on a
  simulated data set ($\alpha=\frac{1}{2}, n=1024$). CART misses the
  left hand change-point entirely on this example; bagged CART seems to
  edge out the posterior mean near $0.75$ but they are quite similar
  over all.\\
{\centering
Key:
 \blue{True $f$},
 \magenta{Posterior Mean},
 \black{CART},
 \dcyan{Bagged CART}}
}
\bigfig{unif_vs_cart2}{90}{Comparative Simulation Experiment: Run
  2}{{\em Simulation Experiment: Run 2}\qquad The CART, Bagged CART, and
  posterior mean estimates are compared with the true regression curve on a
  second (equivalently) simulated data set ($\alpha=\frac{1}{2},
  n=1024$). The posterior mean and bagged CART estimates track each
  other quite closely again.\\
{\centering
Key:
 \blue{True $f$},
 \magenta{Posterior Mean},
 \black{CART},
 \dcyan{Bagged CART}}
}

This experiment involves $10$ simulated data sets each containing
$1024$ data points that were created in the manner explained in the
introductory chapter. In this and in future sections these will be
referred to as experimental runs $1$-$10$. Briefly, one simply chooses
$1024$ random uniforms for the $x$-values, and then flips $1024$
independent coins with the success probability of the coin being given
by the function $f_0$ that is indicated by the blue curve in
figures~\ref{f:unif_vs_cart1}, and~\ref{f:unif_vs_cart2}. These
figures also include a red and a green histogram (drawn upside
down). These are histograms of the of the $x$-values for which the
coin came out heads or tails respectively. There are $75$ bins in each
histogram, and to keep track of the counts, yellow lines are drawn for every 5 events. In these figures, the
posterior mean estimate ($\alpha=\half$) is in magenta, the CART
estimate is in black, and the bagged CART estimate is in cyan. The
first two experiments are shown in a large format, and the latter $8$
are grouped together in \autoref{f:unif_vs_cart_3-10}.

To be specific, in these experiments to compute the CART and bagged
CART estimates, I used the RPART library in S-Plus with default
control settings, 10-fold cross validation, the 1-standard-deviation
rule to control pruning, and the ANOVA method. I apply the RPART
regression algorithms to binary-classification data by assigning
the values $1.0$ and $0.0$ to the two categories. The complexity
parameter was restricted to be at least $0.01$ (this is the default
behavior).

\subsection{Observations}
Here is a summary of notable features in the results of the first two
experiments. Some of these features are explained by the subsequent
discussion. In experimental run 1 (\autoref{f:unif_vs_cart1}), the CART
estimate only has $4$ steps, ``leaving out'' an important split on the
left side.  In experimental run 2 (\autoref{f:unif_vs_cart2}), CART has $5$
splits in about the right places. In both figures, of course, CART
retains its jagged appearance, but it does do an admirable job of
finding the locations where the true curve has change-points. All three
estimates share some features with CART, they all make a fairly abrupt
change at essentially the same place, somewhere near the location of
the true change at $\half$. They all have substantially more trouble
with the change at $\frac{1}{6}$ and this makes sense because it is
much easier to detect the difference between two coins with success
probabilities $0.8$ and $0.4$ than between two coins with
probabilities $0.6$ and $0.4$.  Surprisingly, at least in experimental run
1, CART seems to be a bit more similar to the posterior mean estimate
than to its own bagged version.

Comparing the posterior mean with bagged CART in the first figure,
notice that bagged CART smoothes out some steps that the posterior mean
leaves in. Notably this happens near $0.75$ where bagged CART comes
closer to the truth. For some reason this does not happen near $0.8$
where bagged CART is more blocky than the posterior mean. Bagged
CART's smoothing was a disadvantage near
$\half$, though. Here the posterior mean makes a much sharper
transition and also has an extra ``blip'' up in the correct
direction. Looking at the plateaus, bagged CART has been pulled
closer to $\half$ (away from the truth) than the posterior mean. This
is probably due to averaging in a large number of CART trees that omit
the left hand split.

Comparing the posterior mean with bagged CART in the second figure,
some of the features have remained, but not all. In experimental run 2, the
posterior mean takes a much smoother descent on the right than it did
before. In this case, the posterior mean is closer to
the truth over all. In the first, bagged CART seemed to have an
edge. The posterior mean has a blip near $0.1$ that bagged CART does not.

\subsection{Some Explanations}
Overall, though, in both experimental runs, the bagged CART estimate and the
posterior mean estimate seem quite similar. Consequently, despite the
very different way in which the estimates are arrived at, on this
example at least the computations achieve a similar result. This is
remarkable, especially considering that CART and bagged CART are the
result of years of careful problem specific work and tweaking. The
posterior mean estimate represent an enormous amount of work too, but
most of the work is of a general nature (e.g. MCMC techniques) and not
problem specific. Moreover, this prior is very naive (by design) there
are many ways to modify it that would improve performance on this
example by problem specific tweaking. For example, one could allow
both linear regions and constant regions, or impose a suitable (but
not too restrictive) dependency among the success probabilities that
would help smooth the right side. It may also make sense to space out
the locations of the splits explicitly. Rather, this prior is very
flat. The choice of $\Geometric(\half)$ as the hierarchy prior is not
the result of years of experience, but is, in fact, the first thing I
tried. The prior $\pi$ cannot yet be recommended in general, but that
it even performs modestly well ``out of the box'' on this example is
an excellent defense of the Bayesian approach. 

\mediumfig{bootstrap_unifpost}{-90}{The Bagged Posterior Mean Estimate}
{
{\em The Bagged Posterior Mean Estimate: }
Bagging the posterior mean estimate on the data from experimental run
1 produces the red curve which does not seem to change the posterior
mean estimate much except to add in some irregular wiggles.\\
{\centering Key: \blue{True $f$}, \red{Bagged Posterior Mean}, \magenta{Posterior Mean}, \dcyan{Bagged CART}
\gray{Posterior Mean on a Fifteen Bootstrap Resamples}
}
}

Why, is it, though, that the results are so similar? If the results
were identical one could hope to prove a theorem about why this was
so, but since they are only similar, and since cross-validated CART is
not easily amenable to mathematical analysis (much less its bagged
variant), I can only speculate that they are similar because they both
average together roughly the same functions. They arrive at similar
functions because, after all, they use the same data
set. Additionally, under ordinary circumstances there are bound to be
similarities between the estimates that CART gives and the estimates
that the posterior mean gives (even if I might argue that the
posterior mean estimates are preferable). Consequently, if (ignoring
decision theoretic discipline), I decided to bag the posterior mean
estimate, it ``follows'' that the {\em bagged} posterior mean estimate
should be close to the bagged CART estimate. Since I believe that the
posterior mean estimate is fairly stable under subsampled data
(stability under bootstrap resampling is perhaps more questionable,
but both questions suggest interesting future research), I conclude
that the posterior mean estimate should be close to the bagged CART
estimate. If the reader is skeptical that the posterior mean has
stability under subsampling, they may be interested in the examples
given in \autoref{s:samplesize}, these substantiate this claim (but do
not address it specifically).

The above argument is merely heuristic, of course, so it is reasonable
to ask: what does happen if the posterior mean calculation is bagged?
The answer is shown in \autoref{f:bootstrap_unifpost}. The gray curves
show the results of computing the posterior mean on fifteen
bootstrap-resampled data sets. The red curve is their average: the
``bagged posterior mean.'' Looking at the gray curves, they have many
wobbles and spikes, so it seems that the posterior mean is more
sensitive to the repeated observations that occur in a bootstrap
sample than CART is. As a result, the red bagged posterior mean curve
is itself rather wobbly.

\subsection{Situations in which the Estimates Differ}
It is possible to construct examples where the posterior mean estimate
would differ more dramatically from the CART and bagged CART
estimates. One need only consider circumstances in which CART will reliably
perform in a rather special way.

For a first example, recall the CART works by pruning a full tree and
that this tree is selected in a greedy manner. The first split that
CART chooses, for example, is usually a very important split, but
there are certainly cases where choosing it in a greedy way is
suboptimal if one considers the global search for the ``best tree.''
To some extent bagging improves CART's ability to find useful trees
because sometimes a resampled data set suggests a different splitting
order. Considerations like this are not even an issue for the
posterior proper (because it is a theoretical construct), but they are
an issue for the actual estimates that get produced by MCMC. Still the
issues are different and generally, MCMC methods will perform a much
more ``global'' search over tree space than CART does. This search
ability was emphasized as an advantage of Bayesian CART by~\cite{cgm:1998a}. In
summary, then, in a circumstance where the greedy search has problems,
the posterior mean estimate may avoid those problems, and,
consequently, give a rather different answer than bagged CART. In my
experiments this sort of thing showed up (to a small extent) when I
increased the sample size to $8192$. This
experiment is discussed in \autoref{s:largedata}. The CART estimates
preferred a split in the middle of the right-hand slope that the
posterior mean avoided.

A second example occurs if the $x$-axis is
transformed. One-dimensional CART estimates are invariant with respect
to this, while the posterior mean is not. In some senses this property
is desirable; it seems ``scientific.'' It is not always desirable
though: suppose there is a large amount of data from $0$ to $0.1$ and
from $0.9$ to $1$ but no data in between. Roughly, CART will treat
this in the same fashion as if there were no gap; essentially it only
looks at the ordered values. If it splits the gap at all, it
will split in exactly the middle of the gap: between the
rightmost of the left-hand data and the leftmost of the right-hand
data (ignore the fact that this is not {\em strictly} invariant under
transformations). This effect does not go away under bagging (although
it might get smoothed a bit as the endpoints of the gap
change).

However, the posterior mean estimates will be quite
different. Notably, it will make a smooth transition from the
regression value on the left to the value on the right. Additionally,
because this gap is especially large and the prior specifies that
split points are put down uniformly, it is quite likely that the gap
will be split at least once or twice. If the gap is split twice or
more, then the middle intervals will include no data at all. When
there is no data (or even little data), the prior kicks in to specify
that it thinks that the success probability is uniformly distributed
from $0$ to $1$ and consequently that the mean value of the success
probability on this middle interval is $\half$. Because cases like
this get averaged in, the posterior mean should show some shrinkage to
$\half$ on the gap. Shrinking to $\half$ is, of course, not always
ideal, but in such a case one ought simple to modify the prior and/or
loss function. To construct an example with a gap, I took the data
from experimental run 1 and transformed it so that the left half of
the data now lies in $[0, 0.1]$ and the right half lies in $[0.9,
1]$. The result of computing the posterior mean is shown in
\autoref{f:unifpost_playgap} \mediumfig{unifpost_playgap}{-90}{Data
with a Gap}{{\em Transformed Data with a Gap:} For this example, the
data from experimental run 1 was transformed nonlinearly to create a
gap. Notice the smooth transition that the posterior mean makes along
the gap.}

These features make obvious intuitive sense. Furthermore, if a $95\%$
subjective confidence interval were formed by asking at each point $x
\in [0,1]$ for the smallest interval containing $95\%$ of the
posterior mass, it would grow wider in the middle.  In contrast, if
naive confidence bands were formed around the CART or bagged CART
estimates (by using bootstrapping perhaps) they would make rather
little sense.  

\subsection{Posterior Mean Behavior}
In the first experimental run, CART left out an
``obvious'' split, while in the second it put it in. Sometimes CART
will also ``add in'' splits that it should not have; this is very
sensitive to the particular data set and the pruning rules that are
used. Bagged CART averages all of these together (in some special way
that is hard to formulate, except algorithmically) to arrive at its
smoother curve.

In contrast, I imagine that the posterior mean is considering each of
these possibilities and giving them an appropriate weight before
averaging, in order to give its best estimate. Examining the posterior
mean curve closely, one can see that wherever CART takes a step, the
posterior mean also moves more abruptly than normal. Both estimates
are, after all, both looking for steps, and the locations that CART
chooses are bound to be special parts of the data set, often
containing a run of heads on one side and a run of tails on the
other. Surely this feature would stand out to both methods.

The reverse is not true, however. Consider the ``bump'' in the
posterior mean, just to the right of $\half$.  This results because
the posterior considered functions with additional splits in this area
and gave them weight. It is not, after all, the result of averaging a
large quantity of individually pruned trees, but the result of
averaging over all trees (in principle) with appropriate weights. CART
trees, if allowed to have extra splits would have included this one as
well. Along this line, a modest improvement to the bagged CART
procedure might result if in addition to using different bootstrapped
datasets, one sometimes used different pruning criteria as well, and
then averaged the less penalized trees in together with the more
penalized ones (with appropriate weights).

As an aside, the posterior sometimes becomes more ``sure'' about the
location of a split than it ``really ought to.'' This happens because
of the mismatch between the truth and the prior (or perhaps more
accurately, because of the mismatch is between $\pi$ and my own
subjective prior). The posterior is doing the absolutely optimal thing
if the prior $\pi$ is true (indeed, it is an admissible estimator and
there is no easy way to tell if the others are or not), but according
the prior, functions like the one pictured in blue with a smooth
transition are impossibly rare. When faced with data that could have
resulted from a smooth transition, but might also credibly be
created by a step function with two steps, the posterior only
considers the latter possibility. This is especially apparent when
there is a run of heads and then a run of tails occurs by chance (as
it is bound to do from time to time). The posterior will concentrate
more tightly around this cut-point than makes sense if one considers a
smooth transition to be a credible alternative explanation. The
posterior does, of course, allow for the possibility that there are
multiple splits, but if the success probabilities on each side are
reasonably similar, there may not be enough data to make this
possibility stand out and the single split will remain the most
prominent feature of the posterior mean. This results in a stair-step
appearance that does not go away with larger sample sizes (see
\autoref{s:samplesize}), although the steps tend to get
smaller. Indeed, it appears that the stair-step shape grows more
prominent for larger $n$. Perhaps this is because a larger data set
also is more likely to have at least a few very long runs.

\subsection{Experimental Runs 3-10 and a Summary}
The results from experimental runs $3$ through $10$ are shown
in~\autoref{f:unif_vs_cart_3-10}. By and large, the same observations
made before above apply to these examples. Bagged CART and the posterior mean
track each other quite closely although each occasionally takes a
``wobble'' that the other does not. Broadly speaking both estimators
still retain visible traces of the CART-type functions that they are
averaging together. In particular, both usually have some
``stepiness'' in their appearance; the averaging ``softens'' this, but
does not eliminate it. This is especially visible in experimental run
10 in which both estimates also follow the CART estimate quite
closely. CART leaves out the left hand split on $5$ of the $10$
experiments.
\begin{figure}
\centering
\begin{minipage}[t]{5.5in}
\begin{minipage}[t]{2.5in}
\centering
\minifigx{unif_vs_cart3}{90}{Run 3}\\
\minifigx{unif_vs_cart4}{90}{Run 4}\\
\minifigx{unif_vs_cart5}{90}{Run 5}\\
\minifigx{unif_vs_cart6}{90}{Run 6}
\end{minipage}
\begin{minipage}[t]{2.5in}
\centering
\minifigx{unif_vs_cart7}{90}{Run 7}\\
\minifigx{unif_vs_cart8}{90}{Run 8}\\
\minifigx{unif_vs_cart9}{90}{Run 9}\\
\minifigx{unif_vs_cart10}{90}{Run 10}
\end{minipage}
\end{minipage}
\caption[Comparative Simulation Experiment: Runs 3-10]{
{\em Simulation Experiment: Runs 3-10}\\
{\centering
Key:
 \blue{True $f$},
 \magenta{Posterior Mean},
 \black{CART},
 \dcyan{Bagged CART}}
}
\label{f:unif_vs_cart_3-10}
\end{figure}

\mediumfig{unif_vs_cart_overview}{90}{Comparative
  Scatterplot}{{\em Comparative Scatterplot:} The results of the ten
  experimental runs are summarized by $\cL^2$-norm error (also known
  as $\sqrt{\mathit{MSE}}$). On the $x$-axis is the error committed by
  the posterior mean estimate. On the $y$-axis is the error committed
  by the \black{CART} or \dcyan{bagged CART} estimates respectively.
}
\begin{table}[h]~\label{t:cartsummary}
\centering
\begin{tabular}{|r|p{1.5cm}|p{1.5cm}|p{1.5cm}|}
\hline
  & CART & Posterior Mean & Bagged CART \\
\hline
1 &   0.0891 &   0.0522 &   0.0527 \\
2 &   0.0776 &   0.0478  &  0.0568 \\
3 &   0.0893 &   0.0587  &  0.0710 \\
4 &   0.0996 &   0.0660  &  0.0737 \\
5 &   0.1101 &   0.0659  &  0.0683 \\
6 &   0.1130 &   0.0609  &  0.0799 \\
7 &   0.0779 &   0.0604  &  0.0593 \\
8 &   0.0907 &   0.0669  &  0.0617 \\
9 &   0.1009 &   0.0720  &  0.0718 \\
10&    0.0671&    0.0531 &   0.0587 \\
\hline
$\widehat{\mu}$    &    0.0915  &  0.0604 &   0.0654\\

$\widehat{\sigma}$ &    0.0147  &  0.0076 &   0.0088\\
\hline
\end{tabular}
\comment{
!!Make a table of this data!!
qq =
     cart     postmean baggedcart
1    0.0891    0.0522    0.0527
2    0.0776    0.0478    0.0568
3    0.0893    0.0587    0.0710
4    0.0996    0.0660    0.0737
5    0.1101    0.0659    0.0683
6    0.1130    0.0609    0.0799
7    0.0779    0.0604    0.0593
8    0.0907    0.0669    0.0617
9    0.1009    0.0720    0.0718
10    0.0671    0.0531    0.0587
>> mean(qq)
ans =
    0.0915    0.0604    0.0654
>> std(qq)
ans =
    0.0147    0.0076    0.0088}
\caption[Numerical Summary]{Numerical Summary of $\cL^2$-norm errors
  on the ten experimental runs}
\label{t:baggedcart}
\end{table}

A numerical summary of these ten experiments can be made by computing
the $\mathcal{L}^2$-norm of the error between the estimated curve and
the truth. This summary is given in table~\ref{t:cartsummary} and
illustrated by the scatter-plot in
\autoref{f:unif_vs_cart_overview}. The black points compare CART to
the posterior mean. The cyan points compare bagged CART to the
posterior mean. In each case, the x-axis is the $\L^2$-norm error of
the posterior mean, and the y-axis is that of the
competitor. Obviously, small numbers are preferred, and because the
black points lie exclusively above the identity line on these ten
experiments, the posterior mean is preferable here. The performance of
the bagged CART and the posterior mean estimates is much closer, although
the posterior mean's performance is slightly better on average.

\section{Comparison with Other Popular Methods}\label{s:otherpop}
This section compares the posterior mean estimate with a variety
popular techniques on the data from experimental run
1. Throughout this section, as usual, the true curve is plotted in
blue and the posterior mean with $\Geometric(\half)$ prior is plotted
in magenta.

\subsection{Smoothers}
\mediumfig{Ctable_smoothplot}{-90}{Three Smoothers}{Three Smoothers:\\
{\centering Key: \blue{True $f$}, \magenta{Posterior Mean},
\black{Loess}, \green{Smoothing Spline}, \red{Kernel Smoother} } }

Figure~\ref{f:Ctable_smoothplot} shows the results from running three standard
smoothers. The results are little surprise. The three smoother's estimates are
quite similar over all on this example. They over-smooth the jumps, but
partially make up for this by giving a smoother approximation on the
smooth half. They also all take a turn at the ends; this is consistent
with the data which happens to behave somewhat unusually there.

Loess (plotted in black) fits a locally-weighted linear
regression at each point to make its estimate. Smoothing splines (in
green) use efficient computational tricks to compute the regression
curve that optimizes a tradeoff between small MSE and small
integrated second derivative. Gaussian kernel smoothers (plotted in
red) take a weighted average of response values near a point to
predict. For each method, I chose a smoothing parameter that seemed to
give results that were about as good as possible.
\subsection{LARS/Lasso/Boosting}
\mediumfig{trylasso}{-90}{A Lasso Estimate}{A Lasso Estimate:\\
{\centering Key: \blue{True $f$}, \magenta{Posterior Mean},
\green{Lasso (Cp)} }
}

The green curve in figure~\ref{f:trylasso}, shows the result of using
a Lasso penalized regression. This was particularly easy to do using 
using the Least Angle Regression (LARS)
software~\cite{efron:lars}. Like any linear regression, the results depend on
what basis is chosen. For this example, the basis is constructed by
considering the function $\1{x \ge c}$ for $1200$ values of $c$ evenly
spaced from $0$ to $1$. To encode a datapoint with covariate $x \in
[0,1]$ under this basis, evaluate the $1200$ functions at $x$ and
pack the results into a vector: this vector becomes the covariate that the
regression uses. Finally, Lasso regression resembles ordinary regression
except for one critical difference: the regression parameter $\beta$ is
penalized by $\lambda ||\beta||_1$, where $\lambda$ is a tuning
parameter. For this example, the parameter $\lambda$ was chosen using
a Cp criterion.

Constructed in this way, the Lasso regression estimate should be quite
similar to the estimates that would be arrived at using other important
techniques such as ``Boosting stumps'' and the least angle regression
method. The similarity between these different methods is discussed
in~\cite{efron:lars}.

As can be seen from the figure, the Lasso estimate is has some
appealing features. It is piecewise constant, but it also takes a
fairly large number of steps and spaces them out usefully along the
smooth transition on the right side of the figure. It does a
reasonable jump of ``detecting'' the two change-points. On the other
hand, it does not go as low as it should from $\frac{1}{6}$ to
$\half$, nor as high as it should to the right of $\half$.  
\subsection{Wavelets}
\mediumfig{trywavelet}{-90}{Wavelet Estimates}{Wavelet Estimates:\\
{\centering Key: \blue{True $f$}, \magenta{Posterior Mean},
\green{Wavelet 1}, \red{Wavelet 2}, \black{Wavelet 3} }
}

Figure~\ref{f:trywavelet}
compares some estimates that were based on wavelet techniques. The
red, black, and dotted curves show various wavelet reconstructions of
the regression curve. Overall, I think the results are quite
disappointing; artifacts from the particular basis used show through
clearly into the estimate. Since the data are not regularly spaced, some
accommodation is necessary to use conventional software
(e.g. Wavelab). Algorithms exist that apply directly to irregularly
spaced data, but I did not successfully locate any working
implementations. Instead, the dotted curve shows the wavelet
reconstruction that results from simply using the ordered covariate
values as if they were regularly spaced and then extrapolating back to
the irregularly spaced reality.

There are a number of problems with this approach and some recent work
has developed more sophisticated schemes. For the red and black
curves, I tried one of the simplest~\cite{daubechies}\comment{silverman?} which recommends
using direct linear interpolation to produce values on a fine grid and
then applying wavelet shrinkage methods to this gridded data. The
estimates were computed by the \texttt{wden} function in Matlab. It
has a large number of options, most of which (quite frankly) I do not
understand. These include the choice of a threshold selection rule
from among four options, the choice to use hard or soft thresholding,
an option titled ``multiplicative threshold scaling'' with three
options, the choice of the level at which to compute the coefficients,
and finally the name of the wavelet family to use (many options). For
someone with as little experience with wavelet methods as me, these
options are not a feature but a drawback. Furthermore, experimenting
with the different options, they all seemed to make a difference. A
reasonable, but not heroic effort was made to choose working parameter
values. In any case, for the illustrated curves, the Matlab commands
that were used are:

{\centering
\begin{verbatim}
Xi=seq(min(X),max(X),2^12); % a fine grid of X-values
Yi=interp1(X,Y,Xi);         % on which to interpolate the response
Yhat1 = wden(Y, 'heursure','s','mln', 8,'sym8'); % 1: green
Yhat2 = wden(Yi,'heursure','s','one',10,'db4' ); % 2: red
Yhat3 = wden(Yi,'heursure','s','one',10,'haar'); % 3: black
\end{verbatim}
}

\section{Comparison with Dyadic Prior}\label{s:dyadcomp}
\mediumfig{postdyad_plot}{-90}{Dyadic Posterior}{{\em Dyadic Posterior:}
The posterior mean of the Diaconis and Freedman dyadic prior on the
data from experimental run 1 is drawn in green}
For comparison,
\autoref{f:postdyad_plot} shows the result of computing the posterior
mean resulting from the Diaconis and Freedman dyadic binary regression
prior~\cite{diac:free:1995}. Like the prior $\pi$ studied in this
thesis, this prior chooses a random partition and assigns independent
success probabilities to each partition element. However, this prior
uses a dyadic partition, splitting the data into $2^k$ equal pieces
for some $k \ge 0$. The consistency of the resulting estimates is
guaranteed for any choice of hierarchy prior, except perhaps when the
true regression function is identically $\half$ and the data is pure
noise. For this case, certain priors will be consistent and others
will be inconsistent. It was of interest, then to try one on the
boundary. For this reason the prior on hierarchy level $K$ that was
used assigns: $P(K=k)=(1-\beta)\beta^k$ for
$\beta=\exp(-2^{-\frac{1}{4}}) \approx 0.431$. In fact, for this data,
the results were stable over a wide range of choices of $\beta$.

Since the prior is dyadic, it has no trouble at all nailing the split
at $\half$; of course, it does not have such good luck for the
split at $\frac{1}{6}$. The posterior strongly favors a model with
around $8$ steps.

\section{Dependence on the Parameter of the Geometric Prior}\label{s:geomparam}
\landscapefig{unifpost_plottalk}{-90}{The Posterior Mean for a variety
  of Geometric Priors}{{\em The Posterior Mean for a variety of
  Geometric Priors}: Consider the data of
  experiment 1 and compute the posterior under a
  $\Geometric(1-\alpha)$ prior.}

In this section, consider a departure from the $\Geometric(\half)$
  hierarchy prior. Figure~\ref{f:unifpost_plottalk} shows what the
  posterior mean estimate would be for the data from experimental run
  1, if a $\Geometric(1-\alpha)$ prior is used on the number of steps
  $K$. As usual, the true response curve is indicated in blue, and the
  posterior mean estimate for $\alpha=\frac{1}{2}$ is drawn in
  magenta; the posterior means for other values of $\alpha$ are also
  drawn. As can be seen in the figure, as $\alpha$ ranges from small
  values (short tail prior, dotted black curve) to large values (long
  tail prior, solid black curve) there is not so much difference in
  the posterior mean estimate except that certain small bumps and
  wiggles that are suppressed for the small $\alpha$ values become
  visible for the larger values. Indirectly, this experiment also
  provides a check on the stability of the Monte Carlo estimates of
  the posterior mean; it is unlikely that there would be such close
  agreement among these independently computed estimates if the MCMC was not working reasonably well. To make a
  more detailed comparison, it would be sensible to pool the sampled
  regression curves and use an importance sampling technique to
  compute combined results. I do not pursue this here.

\mediumfig{unifpost_rangehist}{-90}{The Posterior on Model
Size}{{\em The Posterior on Model Size} Consider the data of
  experiment 1 and employ a $\Geometric(1-\alpha)$ prior. A histogram of the number of
  sampled regression functions that had $k$ steps is shown for three
  values of $\alpha$: 0.10 (top), 0.50 (middle), 0.99 (bottom)}

It is also of interest to know how many
splits the posterior is using. Figure~\ref{f:unifpost_rangehist} shows
the posterior on the number of steps $K$ for three
$\Geometric(\alpha)$ hierarchy priors: $\alpha=0.1,$ $0.5$, and $0.9$
respectively. The difference between the priors has a more pronounced
effect here. Under a $\Geometric(1-\alpha)$ prior with $\alpha=0.1$,
complex models are quite rare and consequently, the posterior on step
size shown in the top panel has a fairly short tail. Notice, though,
that even for this conservative model, the likelihood has been able
to overwhelm the conservative prior enough to shift most of the
posterior mass onto models with $5$ splits.

For the bottom panel, $\alpha=0.99$, which corresponds to a rather
slowly decaying tail. Notice, though, that the tail of the
posterior is not nearly this long: models with a large number of steps
$k$ {\em can} fit the data well, but also
have a large number of parameters: $2k-1$. This tends to down-weight
them as a group. When the posterior is marginalized to yield the
posterior distribution on the number of steps, an account is taken of
``how many'' of these more complicated models give a good fit as well
as how good the fit itself is. Because of this trade of, models with
$9$ steps are the most common.

The middle panel, corresponds to the $\Geometric(\half)$ prior that
has been the subject of so many experiments. Notice that models with
fewer than $5$ steps have almost no effect on the posterior mean; most
of the mass is on models with $5$ to $9$ steps: $6$ is the most
common choice.

\section{The Predictive Probability Surface}\label{s:predprob}
\mediumfig{lik_surface}{-90}{A Marginal Likelihood Surface}{{\em A
  Marginal Likelihood Surface}: A slice of the posterior for the $3$
  change-point case ($k=4$) is shown. One change-point is fixed at
  $0.491$ (not shown) and the location of the other two vary along
  the $x$ and $y$ axes of plot respectively.}

To better understand the interaction of the data and the posterior on
experimental run 1, reference to \autoref{f:lik_surface} is useful. It
demonstrates just how spiky, multi-modal, and (in particular)
non-normal the posterior's density can be. It shows a slice of the
posterior for the $3$ change-point case (i.e. four steps). One
change-point is fixed at $0.491$ (not shown) because this was the most
likely location for a single split; this location accounts for the
jump at $\half$ in $f_0$. The other two change-points are allowed to
range over the $x$ and $y$ axes of the plot. More technically, what is
plotted is a self-normalized version of the function $\phi(x)$,
defined by \autoref{e:phix} where $x=\bigl(4,(0.491,\text{\it
x-coord},\text{\it y-coord})\bigr)$. Essentially $\phi(x)$ computes
the likelihood that the splits occur at a given location; it
marginalizes out the different possible choices for the success
probabilities. The height of the surface follows $\phi(x)$, and the
color follows $\log(\phi(x))$, so that the small-probability structure
is also visible. The $x$ and $y$ axes are symmetric, of course,
because splitting at $a$ and $b$ is the same as splitting at $b$ and
$a$. Similarly the function is largest along horizontal and vertical
``bands.'' This is because when one split is in a particularly
fortuitous place, it tends to improve the fit over all, even if the
placement of the second split is suboptimal. The highest two peaks
(near the opposite corners), represent splitting on the left (in the
vicinity of $\frac{1}{6}$) to take care of the jump that is there, and
on the right (in the vicinity of $0.75$) to split the smooth
transition region into its higher and lower halves. The secondary peaks near
the far corner, represent splitting the smooth transition in two
places and ignoring the left half (recall that this was the choice
made by CART on this data).

\section{Behavior on a Small Data Set}\label{s:tinydata}
\mediumfig{postunif_plot}{-90}{Small Data Set Experiment}{{\em Small
  Data Set Experiment:} On this five element data set, the posterior
  mean curve in magenta is the weighted (c.f. inset histogram)
  composition of the colored curves. These curves represent the mean
  contribution from models with a fixed number of steps which ranges
  from $1$ to $15$} It is interesting to see how the posterior
  responds to individual data points; this is most easily seen in a
  very small dataset. In \autoref{f:postunif_plot} I consider a data
  set with 5 data points, that, as usual are shown by the two
  histograms. In this case there are two heads on the left and three
  tails on the right. The posterior mean (magenta) is seen to respond
  in a smooth way, except at the data points where it remains
  continuous but takes a (small) sharp turn. The other colored curves
  represent the posterior mean when the number of steps that the
  function is allowed to take is fixed a-priori. The flat black line,
  for example, is the posterior mean when no splits are allowed (one
  step). The red curve, on the other hand is the posterior mean if 14
  splits (15 steps) are required. Notice how the red curve,
  especially, drifts back towards $\half$ for $x$-values that are not
  close to the data points. This happens because with 14 splits and 5
  data points, there are bound to be many empty intervals and those,
  necessarily, fall back on the prior. The posterior mean curve for a
  $\Geometric(\half)$ prior on the number of steps is the weighted
  average of these curves, where the weights (as percentages) have
  been tabulated by the inset histogram. Since the $1$-split model
  fits this data so especially well, it winds up contributing about as
  much as the contant model (which cannot be ruled out with so little
  data) to the final result. The
  contribution of the higher models is not forgotten, though; notice
  how the posterior mean (magenta) drifts back to $\half$ slightly on
  the left of $0.2$ and the right of $0.8$. Overall the posterior mean
  is conservative; it does not, for example, split the data in half at
  $0.5$ and declare the left hand mean to be $1$ and the right hand
  mean to be $0$.

\section{Behavior on a Large Data Set}\label{s:largedata}
\mediumfig{unifpost_8192compare}{-90}{Large Data Set Experiment}{{\em
    Large Data Set Experiment:} The data from experimental run 1 is
    extended to a large data set of size $8192$\\
{\centering
Key:
 \blue{True $f$},
 \magenta{Posterior Mean ($\alpha=\half$)},
 \green{Posterior Mean ($\alpha=0.995$)},
 \black{CART},
 \dashedblack{min-{\it xerr} CART}
} 
}

The data from experimental run 1, consisted of $1024$ $(x,y)$ pairs
with the $x$-values drawn uniformly from $[0,1]$ and the $y$-values
drawn $\Bernoulli(f_0(x))$. This section answers the question: how
does the posterior change if this data set is enlarged to have $8192$
datapoints by generating additional data from this model?  As usual
the true curve is drawn in blue and the posterior mean for the
$\Geometric(\half)$ prior is drawn in magenta. The data set, as is
visible from the histograms, is getting rather large. Pleasantly, the
posterior mean estimate is also giving an accurate estimate of the
true curve. Somewhat disappointingly, though, the step functions on
which the prior is based have not gone away. Although they are
smaller, they are clearly visible in the magenta curve.

The CART estimate for this data is shown in black. A couple of
observations need to be made. First of all, because the
$1$-standard-deviation pruning rule was used, the CART curve only has
$6$ splits. Looking through the full tree, the model that minimizes
the cross-validated error has two additional splits and is drawn by
the dotted black line. This model agrees quite closely with the
posterior mean estimate, but two of its splits were pruned away when
the $1$-standard-deviation pruning rule was used because the standard
deviation of the {\it xerr} is not small enough. It might be selected
automatically if a more intensive cross-validation were used. Finally,
note that both CART curves minor artifacts (when compared to the
posterior mean) that result from the greedy nature of the full tree.

One might suppose that with this much data, it ought to be possible to
fit a model with many more steps. This does not seem to be true (at
least for a prior that models the success probabilities {\em
independently}). For example, if the full CART tree is manually pruned
to have only one or two additional splits beyond the $8$ that were
used by the minimum {\it xerr} model, the additional splits visibly
degrade the fit. To understand this, consider that the right half of
the data should contain around $4096$ points. By chance, they will not
be (quite) evenly distributed over this half but for simplicity suppose
that these ``$4096$'' points are divided evenly into the $6$ intervals
selected by the larger CART model. This leaves around $680$ points in
each partition cell.  Recall that the variance of a $\Binomial(n,p)$
random variable is $n p(1-p)$, so that the standard deviation of the
estimated success probability on each of these partition cells is
going to be around $0.02$. Considering that the regression estimator
has to not only detect a difference, but also locate a good choice of
split, and optimize the accuracy of the estimated success probability,
it does not seem too unreasonable that the average jump in success
probability between neighboring cells is around $0.1$.

\comment{
variance of the difference
between two success probability estimates based on two \iid\
$\mathit{Binomials}$ would have variance $2p(1-p)/n$. For $p=0.6$,
then, the standard deviation of the difference would be about
$0.027$. To have {\em approximately} $50\%$ chance of detecting such a
difference (at a $95\%$ significance level), the means need to be
about 2 standard deviations apart. Roughly then, one could expect a
$50\%$ chance of detecting a change in height, with this amount of
data, if the success probability on two neighboring intervals differed
by about $0.053$. Looking at the $6$ steps, they seem to descend in
success probability by about $0.1$ each time.
standard deviation such random variables is for $p=0.6$, the standard
deviation of an estimate of the success probability based on this
random variable is about $0.016$. Using this heuristically, it seems
that the change in the success probabilities from step to step over
the transition region
}

Also shown (by the green curve) is the result of computing the
posterior mean when a $\Geometric(1-\alpha)$ prior is used for
$\alpha=0.995$. Interestingly, this long tail makes little
difference and the green curve barely peaks out from under the magenta
one.

\section{The Effect of Sample Size}\label{s:samplesize}
\begin{figure}
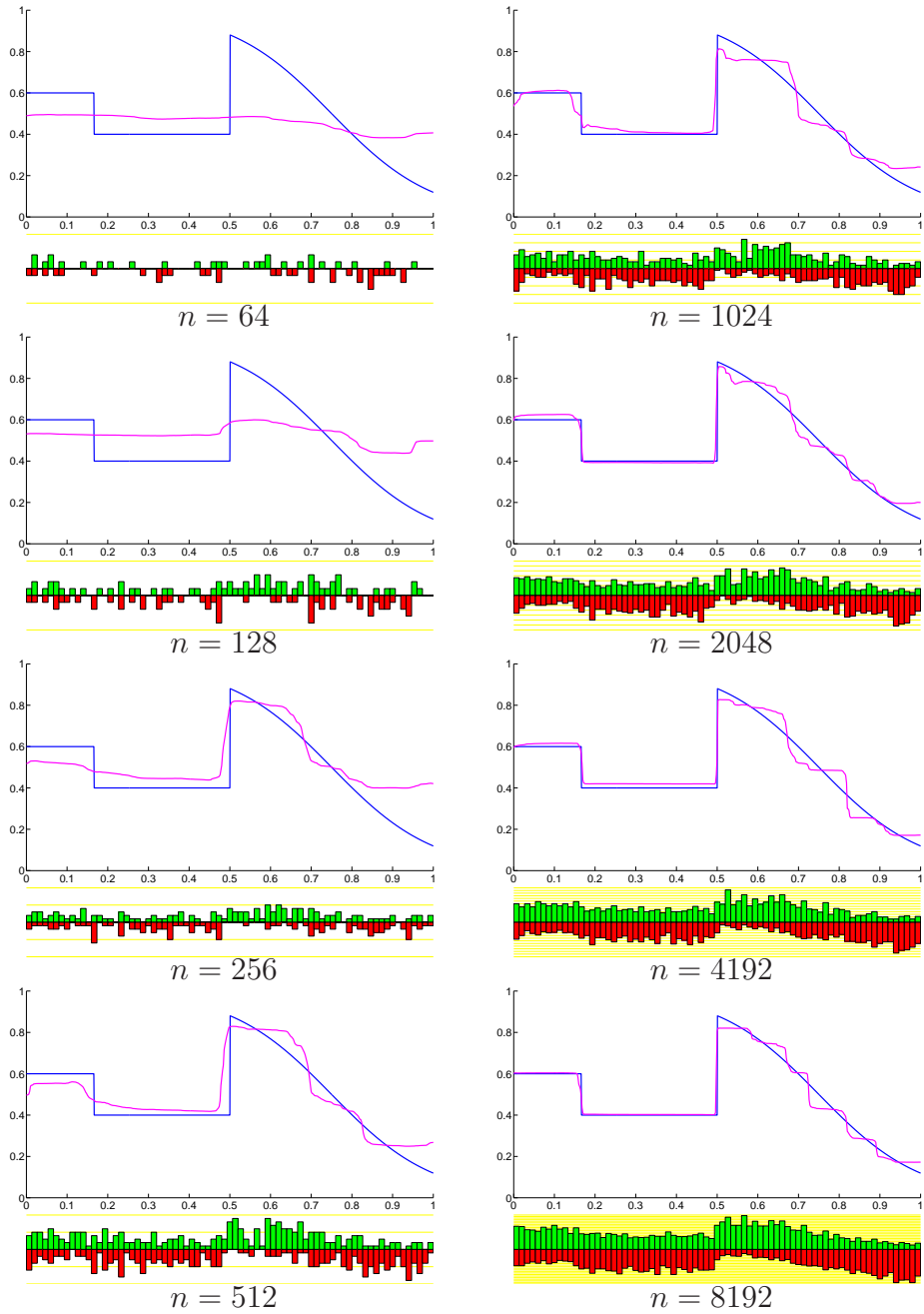

\centering
\begin{minipage}[t]{5.5in}
\begin{minipage}[t]{2.5in}
\centering
\minifigx{unifpost_playn_64}{-90}{$n=64$}\\
\minifigx{unifpost_playn_128}{-90}{$n=128$}\\
\minifigx{unifpost_playn_256}{-90}{$n=256$}\\
\minifigx{unifpost_playn_512}{-90}{$n=512$}
\end{minipage}
\begin{minipage}[t]{2.5in}
\centering
\minifigx{unifpost_playn_1024}{-90}{$n=1024$}\\
\minifigx{unifpost_playn_2048}{-90}{$n=2048$}\\
\minifigx{unifpost_playn_4096}{-90}{$n=4192$}\\
\minifigx{unifpost_playn_8192}{-90}{$n=8192$}
\end{minipage}
\end{minipage}
\caption[The Effect of Sample Size]{{\em The Effect of Sample Size:}
The posterior mean is computed as sample size $n$ runs through a wide
range
}
\label{f:playn}
\end{figure}

The previous section developed a extended version of experimental run
1 that contained $8192$ data points. In figure~\ref{f:playn}, this
data is analyzed in more detail. Smaller data sets are formed by
taking the first $n$ data points for $n$ ranging from $64$ to $8192$
by powers of two. It is very pleasant to see how the posterior mean
incrementally grows closer to the truth. At first, the steps on the
left are almost ignored (with so little data, any pattern they contain
could have resulted from noise), but gradually they fill in. The
transitions near the change-points in $f_0$ become very sharp and the
estimates of the smooth transition steadily improve.

\section{The Effect of Sample Size: a Harder Example}
\begin{figure}
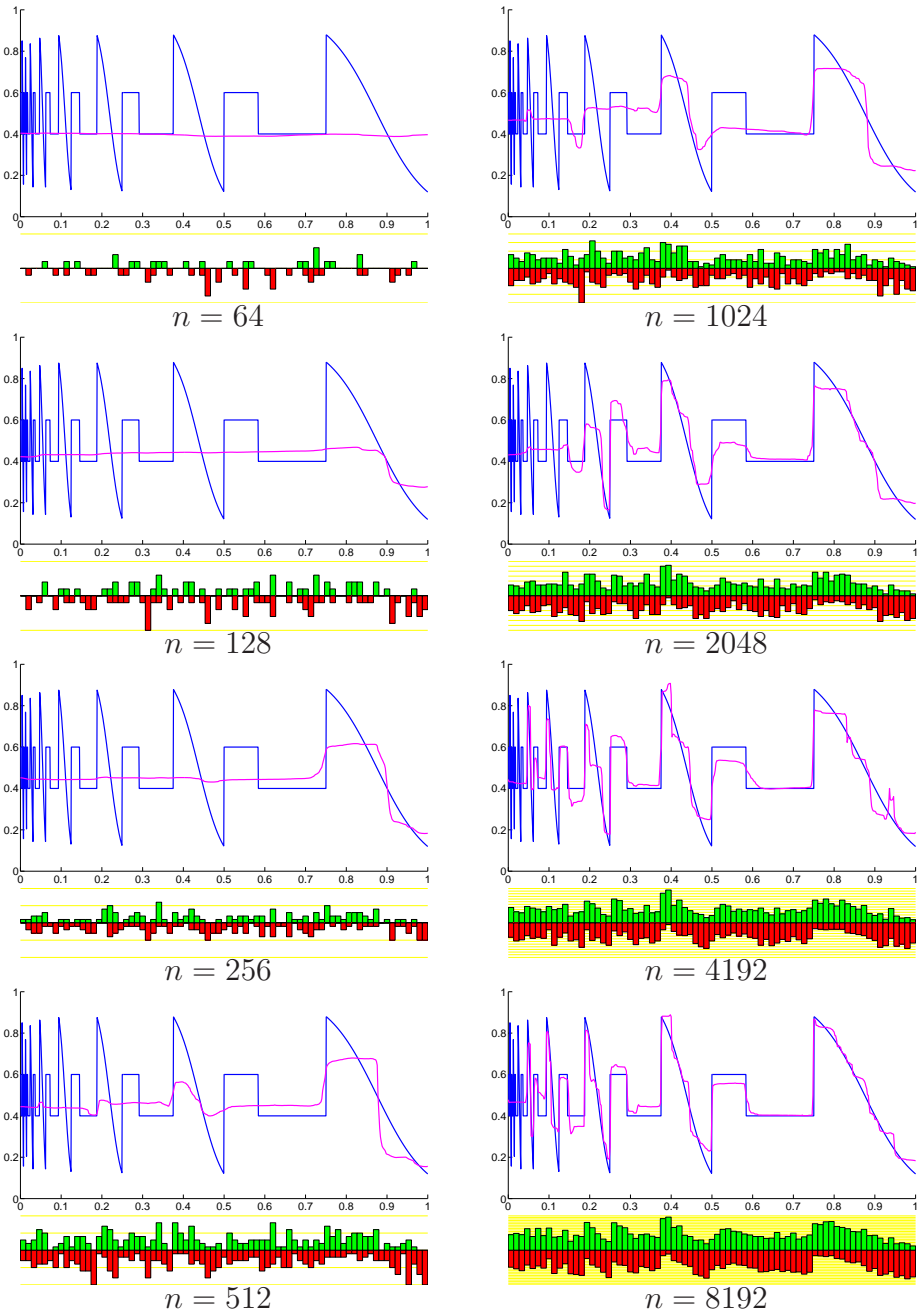

\centering
\begin{minipage}[t]{5.5in}
\begin{minipage}[t]{2.5in}
\centering
\minifigx{unifpost_playxx_64}{-90}{$n=64$}\\
\minifigx{unifpost_playxx_128}{-90}{$n=128$}\\
\minifigx{unifpost_playxx_256}{-90}{$n=256$}\\
\minifigx{unifpost_playxx_512}{-90}{$n=512$}
\end{minipage}
\begin{minipage}[t]{2.5in}
\centering
\minifigx{unifpost_playxx_1024}{-90}{$n=1024$}\\
\minifigx{unifpost_playxx_2048}{-90}{$n=2048$}\\
\minifigx{unifpost_playxx_4096}{-90}{$n=4192$}\\
\minifigx{unifpost_playxx_8192}{-90}{$n=8192$}
\end{minipage}
\end{minipage}
\caption[The Effect of Sample Size: a Harder Example]
{{\em The Effect of Sample Size: a Harder Example} The posterior mean is computed as sample size $n$ runs through a wide
range on this more challenging problem}
\label{f:playnxx}
\end{figure}

For the final example, I consider a much harder regression
function. It is depicted in blue in \autoref{f:playnxx}. It was formed
by taking multiple copies of the original $f_0$ and shrinking them to
half their size repeatedly. A data set with $8192$ point is simulated
and the posterior mean is calculated for subsets of increasing
size. As had been hoped, the features get filled in as the data size
increases incrementally. The larger features rise above the noise
first and then the smaller, so that for this regression function, the
approximation seems to grow better as $n$ increases on the right
first, but then steadily spreads to the left. Importantly, for the
larger $n$, the posterior mean concentrates on models with many more
splits than it did for the easier data since. Presumably this is
because more complex models are necessary to increase the
likelihood. To some extent this makes sense because even though this
example is more complex than the previous one, some aspects of the
regression function are relatively easy to detect (e.g. the large
jump) and there are more of these features available for analysis.  It
is interesting to compare the results with the former dataset for
$n=4096$ with the results {\em on the right half} of the current
dataset for $n=8192$. They should be quite comparable because the
regression function for this more complicated model on $[\half, 1]$ is
just a rescaled version of the original model, and in both cases there
should be about $4096$ datapoints available. To appearances, the two
results are quite similar except that the latter result does not do as
well on the small interval from $0.5$ to around $0.6$. Perhaps it is
being confused by the low success probability region immediately to
the left of $0.5$. Conducting a similar comparison between the
original $n=1024$ example the right half of the $n=2048$ example, the
latter result seems substantially inferior. On the other hand, it is
not so bad when compared with the original result on experimental run
5. Furthermore, the performance of the posterior mean is quite good
considering the overall increase in the difficulty of this
problem. Even more amazing, considering the popular state-of-the-art
methods, it does all of this without any tuning or
cross-validation.



%% file: proof.tex
\chapter{Consistency}\label{c:proof}
\nocite{wong:shen:1995}

This chapter establishes conditions under which the prior for
one-dimensional classification that was introduced in
\autoref{s:priorintro} is a consistent estimator of the true
regression function. The consistency
of the posterior is proven using a result by Barron, Schervish, and
Wasserman~\cite{barron:schervish:wasserman:1999}. Before reviewing
their theorem, I pause to introduce some notation. I also present a
lemma that shows that the original conditions given in their theorem
are equivalent to some others that may be easier to check. Finally, I
specify the prior $\pi$ more formally and complete a proof of
consistency.

\section{Notation and the Basic Theorem}\label{s:notation}
Let $\cF$ be the class of all (Borel) measurable functions
$f:[0,1]\mapsinto\{0,1\}$. Write $\mu$ for the uniform distribution on
$[0,1]$ and $\eta$ for counting measure on the set $\{0,1\}$. Write
$\cZ$ for the product space $[0,1] \times \{0,1\}$, and call the
product measure $\nu$. To any $f \in \cF$ there is a corresponding
density on $\cZ$ with respect to $\nu$ that we denote by $\wf$:
\begin{equation}\label{e:wf}
\wf(x,y)=f(x)\1{y=1}+(1-f(x))\1{y=0}
\end{equation}
For notational convenience, we may write either $\wf(x,y)$ for $x \in
[0,1]$ and $y \in \{0,1\}$ or $\wf(z)$ for $z=(x,y) \in \cZ$. Write
$F=F_f$ for the distribution on $\cZ$ whose density with respect to $\nu$ is
$\wf$. In words, the consequence of this construction is that sampling
a point $Z$ from $F$ is the same as choosing an $X$ uniformly on
$[0,1]$ and then (conditionally on $X=x$) determining $Y$ by flipping an
``$f(x)$'' coin.

Write $\wcF$ for the class of densities formed by considering $\wf$
for every $f\in\cF$.

Let $d$ denote the Hellinger distance on $\wcF$:
\begin{equation}
d(\wf,\wf')=\bigl\{ \int
\bigl(\wf(z)^{\frac{1}{2}}-{\wf'(z)}^{\frac{1}{2}}\bigr) \nu(dz)
\bigr\}^{\frac{1}{2}}
\end{equation}

And let $D$ denote the Kullback-Leibler discrepancy on $\wcF$ (employing
the usual convention that the integrand is interpreted as 0 whenever
$\wf(z)=0$):
\begin{equation}
D(\wf,\wf')=\int \log{\frac{\wf(z)}{\wf'(z)}\wf(z)} \nu(dz)
\end{equation}

We are concerned with posterior consistency and so the mass
that the prior or posterior ascribes to certain small sets containing
the true parameter $\wf_0$ is of interest. For any $\epsilon>0$, we define two such
``neighborhoods:''
\begin{align}
K_\epsilon &= \{ \wf \in \wcF ~|~ D(\wf_0,\wf) \leq \epsilon \} \\
H_\epsilon &= \{ \wf \in \wcF ~|~ d(\wf_0,\wf) \leq \epsilon \}
\end{align}

The ``richness'' of the parameter space is also an important quantity;
to make this precise we supply the following definitions.

\begin{definition}
Consider a class of functions $\cA$ that are densities with respect to
dominating measure $\nu$. We say that the collection of functions
$\{\wf_1^U, \dots, \wf_r^U\}$ is a $\delta$-upper bracketing of
$\cA$ if:
\begin{enumerate}
\item for every $a \in \cA$ there exists $i$ such that $a \le \wf_i^U$
a.e. $[\nu]$
\item every $\wf_i^U$ satisfies $\int \wf_i^U(z) \nu(dz) \le 1 + \delta$
\end{enumerate} 

Furthermore, write $\cH(\cA, \delta)$ for the $\delta$-upper metric
entropy of $\cA$ which is the logarithm of the size of the smallest
possible $\delta$-upper bracketing of $\cA$ (infinity, if no finite
bracketings exist).
\end{definition}

The following result can now be stated. It gives conditions on the
prior that are sufficient to ensure that the posterior will
concentrate on $H_\epsilon$ for any $\epsilon>0$.

\begin{thm}[Barron, Schervish, and Wasserman 1999]
Let $\nu$ be a $\sigma$-finite measure on a measurable space $(\cZ,
\cB)$, where the $\sigma$-field $\cB$ is separable. Let $\pi$ be a
probability distribution on $\wcF$, a class of probability densities
with respect to $\nu$. Endow $\wcF$ with the Borel $\sigma$-field induced by the Hellinger metric $d$. Let $\wf_0$ be
a certain chosen density with respect to $\nu$ and write $F_0$ for the
corresponding distribution on $\cZ$. Let $Z_1, Z_2, \dots$ be drawn
$\iid$ from $F_0$. In the notations explained above, further assume
that for every $\epsilon>0$:
\begin{enumerate}
\item $\pi(K_\epsilon)>0$
\end{enumerate}
\begin{enumerate}
\setcounter{enumi}{1}
\item There exists a sequence $\{\cA_n\}_{n=1}^\infty$ of measurable subsets of $\wcF$
  and positive, real numbers $d_1,d_2,c,$ and $\delta$ such that:
  \begin{enumerate}
  \item $\pi(\cA_n^c) \leq d_1 \exp( -d_2 n)$ for all $n$ sufficiently large
  \item $\cH(\cA_n, \delta) \leq c n $ for all $n$ sufficiently large
  \item $c < \left( (\epsilon - \sqrt{\delta})^2 - \delta \right)/2$
  \item $\delta < \epsilon^2/4$
  \end{enumerate}
\end{enumerate}
Then, with probability $1$ [under $F_0^\infty$ measure], Bayes
theorem applies for all $n$; i.e. for any measurable $B \subset \wcF$
and any $n$:
\[
\pi(B | z_1, \dots, z_n)=
\frac{\int_B \prod_{i=1}^n \wf(z_i) \pi(d\wf)}
{\int_{\wcF} \prod_{i=1}^n \wf(z_i) \pi(d\wf)}
\]
And for any $\epsilon>0$:
\[
\lim_{n \tendsto \infty} \pi(H_\epsilon | z_1, \dots, z_n) = 1
\]
\end{thm}

\vskip 1cm
\renewcommand{\labelenumi}{\Roman{enumi}}
\renewcommand{\labelenumii}{\alph{enumii}} 

The second condition of this theorem seems rather technical. To
restate this theorem in simpler terms we prove an elementary lemma
which shows that these conditions (item {\it I} in the lemma) can be
expressed in three other equivalent forms.

\begin{lemma}
The following conditions, given in the notation defined above, are
equivalent:

\begin{enumerate}
\item For all $\epsilon>0$, there exists
  $\{\cA_n\}_{n=1}^\infty$, a sequence of measurable subsets of $\wcF$
  and positive, real numbers $d_1,d_2,c,$ and $\delta$ such that:
  \begin{enumerate}
  \item $\pi(\cA_n^c) \leq d_1 \exp( -d_2 n)$ for all $n$ sufficiently large
  \item $\cH(\cA_n, \delta) \leq c n $ for all $n$ sufficiently large
  \item $c < \left( (\epsilon - \sqrt{\delta})^2 - \delta \right)/2$
  \item $\delta < \epsilon^2/4$
  \end{enumerate}
\item There exists a sequence of positive real numbers
  $\{\delta^i\}_{i=1}^\infty$, with $\delta^i \tendsdown 0$, and
  $\left\{\{\cA_n^i\}_{n=1}^\infty\right\}_{i=1}^\infty$, a sequence
  of sequences of measurable subsets of $\wcF$ such that:
  \begin{enumerate}
  \item for all $i$, $\limsup_n~ \log(\pi\left((\cA_n^i)^c)\right)/n < 0$
  \item $\limsup_i \limsup_n \cH(\cA_n^i, \delta^i)/n = 0 $
  \end{enumerate}
\item For all $\epsilon>0$ there exists
  $\{\cA_n\}_{n=1}^\infty$, a sequence of measurable subsets of $\wcF$
  and a positive, real number $\delta \leq \epsilon$ such that:
  \begin{enumerate}
  \item $\limsup_n~ \log(\pi\left((\cA_n)^c)\right)/n < 0$
  \item $\limsup_n \cH(\cA_n, \delta)/n \leq \epsilon $
  \end{enumerate}
\item\label{cond3} For all $\epsilon>0$ there exists
  $\{\cA_n\}_{n=1}^\infty$, a sequence of measurable subsets of $\wcF$
  such that:
  \begin{enumerate}
  \item $\limsup_n~ \log(\pi\left((\cA_n)^c)\right)/n < 0$
  \item $\limsup_n \cH(\cA_n, \epsilon)/n \leq \epsilon $
  \end{enumerate}
\end{enumerate}
\end{lemma}
\begin{proof}
{\it I} $\implies$ {\it II}: 

First, notice that the condition from
{\it I.a} that ``there exist $d_1>0$ and $d_2>0$ so that $\pi(A_n^c)
\leq d_1 \exp(-d_2 n)$ for $n$ sufficiently large'' implies that
$\limsup_n \log( \pi(\cA_n^c) )/n \leq -d_2 < 0$. Conversely, if
$\limsup_n \log( \pi(\cA_n^c) )/n = \alpha <0$, then for all
sufficiently large $n$, $\log( \pi(\cA_n^c))/n \leq \alpha/2$ and
$\pi(\cA_n^c) \leq \exp(\frac{\alpha}{2} n)$ for all sufficiently large $n$.

Choose a sequence of $\epsilon^i>0$, $\epsilon^i \tendsdown 0$ and use
{\it I} to establish the existence of $\{\cA_n^i\}, d_1^i, d_2^i, c^i,
\delta^i$ satisfying {\it I} for $\epsilon^i$. Since $\delta^i <
(\epsilon^i)^2/4 \tendsdown 0$, we can find a subsequence on which
$\delta^i \tendsdown 0$. Without loss of generality, assume we already
have such a subsequence. The above reasoning, then, establishes {\it
II.a} for this sequence. Since {\it I.a} implies that $\limsup_n
\cH(\cA_n^i, \delta^i)/n \leq c^i$ for all $i$, to establish {\it
  II.b}, we need only show that $\limsup_i c^i=0$. Condition {\it I.c}
constrains $c^i <  \left( (\epsilon^i - \sqrt{\delta^i})^2 - \delta^i
\right)/2$. Viewing this as a function of $\delta^i$, notice that it
is monotonically decreasing on $[0, (\epsilon^i)^2/4]$ so that,
necessarily, $c^i < (\epsilon^i)^2/2$, the value of this function
at $\delta^i=0$.

{\it II} $\implies$ {\it I}: 

Consider some $\epsilon>0$. Note that $\left( (\epsilon - \sqrt{\delta^i})^2 - \delta^i
\right)/2 \tendsup \epsilon^2/2$ as $i \tendsto \infty$.  Find an $i$
sufficiently large so that this expression is at least
$\epsilon^2/3$. Then find a subsequent $i$ sufficiently large so that
$\limsup_n \cH(\cA_n^i, \delta^i)/n \le \epsilon^2/5.$  Choose
$c=\epsilon^2/4,$ and $\delta=\delta^i$, and {\it I} is proven.

{\it II} $\iff$ {\it III}:
This is a straightforward exercise in nitpicking.

{\it III} $\iff$ {\it IV}:
Observe from the definition that $\cH(\cA_n,\delta)$ is a non-increasing function of
$\delta$.

\end{proof}

\renewcommand{\labelenumi}{\arabic{enumi}}
\renewcommand{\labelenumii}{\alph{enumii}} 

\begin{thm}[Corollary to Barron, Schervish, and Wasserman 1999]\label{t:main}
Let $\nu$ be a $\sigma$-finite measure on a measurable space $(\cZ,
\cB)$, where the $\sigma$-field $\cB$ is separable. Let $\pi$ be a
probability distribution on $\wcF$, a class of probability densities
with respect to $\nu$. Endow $\wcF$ with the Borel $\sigma$-field induced by the Hellinger metric $d$. Let $\wf_0$ be
a certain chosen density with respect to $\nu$ and write $F_0$ for the
corresponding distribution on $\cZ$. Let $Z_1, Z_2, \dots$ be drawn
$\iid$ from $F_0$. In the notations explained above, further assume
that for every $\epsilon>0$:
\begin{enumerate}
\item $\pi(K_\epsilon)>0$
\end{enumerate}
There exists a sequence $\{\cA_n\}_{n=1}^\infty$ of measurable subsets of
$\wcF$,  so that:
\begin{enumerate}
\setcounter{enumi}{1}
\item $\limsup_n~ \log(\pi\left((\cA_n)^c)\right)/n < 0$
\item $\limsup_n~ \cH(\cA_n, \epsilon)/n \leq \epsilon $
\end{enumerate}
Then, with probability $1$ [under $F_0^\infty$ measure], Bayes
theorem applies for all $n$; i.e. for any measurable $B \subset \wcF$
and any $n$:
\[
\pi(B | z_1, \dots, z_n)=
\frac{\int_B \prod_{i=1}^n \wf(z_i) \pi(d\wf)}
{\int_{\wcF} \prod_{i=1}^n \wf(z_i) \pi(d\wf)}
\]
And for any $\epsilon>0$:
\[
\lim_{n \tendsto \infty} \pi(H_\epsilon | z_1, \dots, z_n) = 1
\]
\end{thm}

\vskip 1cm

In words, this is what we have required: 1) the prior must put
positive mass on all Kullback-Leibler neighborhoods of $\wf_0$. 2) We
must be able to choose an increasing sequence of subsets of the
parameter space so that the $n$'th of these sets captures all but
exponentially much of the prior mass 3) This sequence must not grow in
``complexity'' too quickly. We conclude that the posterior will concentrate
on the subset of the parameter space which is Hellinger close to the
true parameter.

\section{Specification of the Prior}\label{s:priorspec}

We describe $\pi$ as a distribution on $\wcF$ by means of first
describing a parametric prior $\pi'$. Let $\kappa$ be a distribution
on the positive integers. Let $\Theta_k=\{(k,\v,\s)~:~ \v \in
[0,1]^{k-1},\s \in [0,1]^k \}$ and let $\Theta=\union_{k=1}^\infty
\Theta_k$. Let $\pi'$ be the distribution on $\Theta$, that can be described
by first picking $K$ according to $\kappa$; and then, conditional on $K=k$,
picking a point $\theta \in \Theta$ uniformly from $\Theta_k$. That is, $\v$ and
$\s$ are chosen independently and uniformly from the appropriate unit
cubes. Let $\v_{(i)}$ denote the $i$'th ordered value of $\v$.  Now to
any $\theta \in \Theta$ associate the function $f_\theta \in \cF$ given
by the following construction. For $k>1$, let $I_1=[0,\v_{(1)}),
I_2=[\v_{(2)}, \v_{(3)}), \dots, I_{k-1}=[\v_{(k-2)},\v_{(k-1)}),
I_{k}=[\v_{(k-1)},1]$. If $k=1$, just let $I_1=[0,1]$. Take
$f_\theta(x)=\sum_{i=1}^k \s_i \1{x \in I_i}$. Finally, then, to draw
$\wf$ from $\pi$, draw $\theta$ from $\pi'$, construct $f_\theta$, and
form $\widetilde{f_\theta}$ as in \autoref{e:wf}.

\comment{

\begin{equation}
f_\theta(x)=
\begin{cases}
k=1 \mathtext{~or~} x < \v_{(1)} & \s_1 \\
k>1 \mathtext{~and~} x \geq \v_{(k-1)} & \s_{k} \\
j \in \{2, \dots, k-1\}, \v_{(j-1)} \leq x < \v_{(j)} & \s_j
\end{cases}
\end{equation}
}

\section{A Consistency Proof}\label{s:proof}

This section proves that the prior $\pi$ just described is consistent
in the sense that, under repeated sampling, the posterior mass will
concentrate on a Hellinger neighborhood $H_\epsilon$ for any $\epsilon>0$.
The proof uses three lemmas to establish the conditions of
Theorem~\autoref{t:main}.

The first lemma shows that $\pi$ puts mass on all
Kullback-Leibler neighborhoods $K_\epsilon$.
\begin{lemma}\label{l:pack}
Let $\pi$ be the prior distribution on $\wcF$, the class of densities
w.r.t $\nu$, that was described above. Assume that $\kappa$, the
hierarchy prior, assigns positive mass to every natural number. Let
$\widetilde{f_0}$ be an arbitrary density in $\wcF$. Then, for any
$\epsilon>0$, $\pi(K_\epsilon)>0$.
\end{lemma}

\begin{proof}
Let $f_0$ be the corresponding function in $\cF$. For some
$0<\delta_1<\frac{1}{16}$, to be determined later, let:
\[
\mathcal{G}=\{g \in \F : ~~ ||g-f_0||_1 \le 2\delta_1 \text{~and~}
\delta_1 \le g(x) \le 1-\delta_1 ~~(\forall x) \}
\]

Let $A_g=\{x \in [0,1]~:~ |g(x)-f_0(x)| \leq \delta_2\}$ where
$\delta_2 \defined \sqrt{2\delta_1}$. By Chebyshev's inequality, for
any $g \in \mathcal{G}$, the Lebesgue measure of ${A_g}^c$ is 
smaller than $\delta_2$.

\begin{figure}
\centering\epsfig{figure=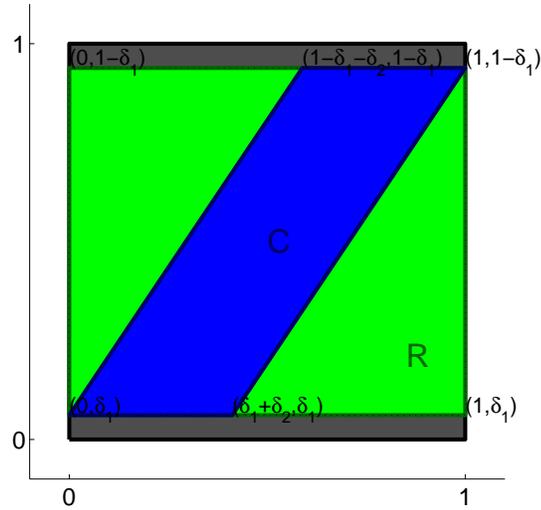,height=7cm,angle=0}
\caption{The sets $C$ and $R$}
\label{f:D_is_convex}
\end{figure}

In other words, as shown in \figref{f:D_is_convex}, if $x\in A_g$, this
restricts the pair $\left(f_0(x),g(x)\right)$ to lie in the convex set
$C$ (blue) whose extreme points are $(0,\delta_1)$,
$(\delta_1+\delta_2,\delta_1)$, $(1,1-\delta_1)$,
$(1-\delta_1-\delta_2,1-\delta_1)$. For $x \notin A_g$, we still have a
restriction on $g(x)$, so that the pair $\left(f_0(x),g(x)\right)$
lies in the rectangle $R$ (green) with vertices $(0,\delta_1)$,
$(1,\delta_1)$, $(1,1-\delta_1)$, $(0,1-\delta_1)$. Let $D_1(p,p')$
denote the Kullback-Leibler discrepancy between the $\Bernoulli(p)$
and $\Bernoulli(p')$ distributions.
\[
D_1(p,p') \defined p \log \frac{p}{p'} + (1-p) \log \frac{1-p}{1-p'}
\]
Then for any $g \in \mathcal{G}$ we can bound
$D(\widetilde{f_0},\widetilde{g})$, as follows:
\begin{align}
D(\widetilde{f_0},\widetilde{g}) &=   \int_{A_g} D_1(f_0(x),g(x)) + \int_{{A_g}^c} D_1(f_0(x),g(x)) \\
&\le
   \int_{A_g} \sup_{(a,b) \in C} D_1(a,b)+
   \int_{{A_g}^c} \sup_{(a,b) \in R} D_1(a,b) \\
\intertext{and because $D_1(a,b)$ is convex in the pair $(a,b)$
  \cite{cover}, the supremum is achieved at the vertices so that the
  above is bounded by:}
&\le
   \max_{(a,b) \in \{\text{vertices of~} C\} }D_1(a,b) +
   \delta_2 \max_{(a,b) \in \{\text{vertices of~} R\}} D_1(a,b)
   \\
\intertext{Using the symmetry $D_1(a,b)=D_1(1-a,1-b)$}
&\le
   \max\left(
   D_1(0,\delta_1), D_1(\delta_1+\delta_2,\delta_1)
   \right)
   +\delta_2 \max\left(
   D_1(0,\delta_1),D_1(0,1-\delta_1)
   \right)
   \\
\begin{split}
&= \max \biggl( -(1-\delta_1)\log(1-\delta_1), \\
& \qquad\qquad (\delta_1+\delta_2)
\log(\frac{\delta_1+\delta_2}{\delta_2}) +
(1-\delta_1-\delta_2)\log(\frac{1-\delta_1-\delta_2}{1-\delta_1})
\biggr) \\
& \qquad + \delta_2( -\log(\delta_1))
\end{split}\\
&\le
   \max\left(
   -\log(1-\delta_1),
   2\log(2)\delta_2 
  \right)
   +\delta_2( -\log(\delta_1))   \\
\intertext{For small $\delta_1$, the last term is the most
   important. To simplify the first term, verify that for $0 \le \delta \le \frac{1}{2}$, $-\log(1-\delta) \le
  \delta+\delta^2 \le \frac{3}{2}\delta$}
&\le 
   \max\left(
   \frac{3}{2}\delta_1,2\log(2)\delta_2
   \right)
   +\delta_2( -\log(\delta_1))
   \\
&\le 
   \frac{3}{2}\delta_2 + \delta_2( -\log(\delta_1) )\\
&= 
   \sqrt{2 \delta_1}( -\log(\delta_1) +\frac{3}{2})
\end{align}

Which tends to zero as $\delta_1 \tendsdown 0$. For a sufficiently
small choice of $\delta_1$, then:
\[
\widetilde{\mathcal{G}}=\{\widetilde{g}~:~ g \in \mathcal{G}\}
\subset K_\epsilon=\{\wf \in \wcF~:~ D(\wf_0,\wf)<\epsilon\}.
\]

It remains to show that $\widetilde{\mathcal{G}}$ has positive prior
mass. To do this, we first find a step function $g_0 \in \mathcal{G}$
which approximates $f_0$ and then show that $\mathcal{G}^*$, a set of
perturbations of $g_0$, remains in $\mathcal{G}$ and that
$\widetilde{\mathcal{G}^*}$ has positive prior mass.

Observe that since $f_0$ is Lebesgue-measurable there exist two
increasing sequences of step functions $h_i^{+}$ and $h_i^{-}$ for
which $||f_0 - (h_i^{+}-h_i^{-})||_1 \tendsto 0$. This is, in fact,
the basis of a common construction of the Lebesgue
integral~\cite{haaser:sullivan:1971}. Consequently, we can find a step
function $h \in \F$ for which $||f_0-h||_1 \le \delta_1/2$. Let $g_0$
be the function obtained by modifying $h$ so that it always remains in
$[\delta_1,1-\delta_1]$,
i.e. $g_0(x)=\max(\min(h(x),1-\delta_1),\delta_1)$, so that $g_0 \in
\mathcal{G}$ and $||f_0-g_0||_1 \le \frac{3}{2}\delta_1$.

Now, parameterize $g_0$ in the manner of \secref{s:priorspec}
(changing it at a set of measure 0 if necessary) so that $k<\infty$ is
the number of locally constant regions in $g_0$ and the vectors
$\mathbf{v}=(v_i)_{i=1}^{k-1}$ and $\mathbf{s}=(s_i)_{i=1}^k$
signify the change-points and success probabilities of $g_0$.

If we then perturb each $s_i$ towards the value $\frac{1}{2}$ by no
more than $\delta_1/4$ we obtain a new function $g'$ which remains in
$\mathcal{G}$ and satisfies $||g_0-g'||_1 \le \delta_1/4$. If, in
addition, we perturb the $v_i$'s by no more than $\delta_1/(4(k-1))$
we obtain $g''$ which also remains in $\mathcal{G}$, satisfying
$||g_0-g''||_1 \le \delta_1/2$ or $||f_0-g''||_1 \le 2 \delta_1$.

Denote the class of functions thus obtained by $\mathcal{G}^*$. Then
$\pi(\widetilde{\mathcal{G}}) \ge \pi(\widetilde{\mathcal{G}^*}) \ge
\kappa{\{k\}}(\delta_1/4)^k(\delta_1/(2(k-1)))^{k-1}>0$.
\end{proof}

\begin{remark} A shorter proof can be based on the martingale
sequence formed from conditional expectations of $f_0$.
\end{remark}

For the second lemma, consider $\F_m$, the class of functions which
are constant except possibly for as many as $m-1$ change-points.
\begin{definition}
Using the notation from
\autoref{s:priorspec}, let $\F_m=\{f_\theta ~:~ \theta \in
\Theta_m\}$. Let $\wcF_m$ be the associated class of densities with
respect to $\nu$:
\[
\widetilde{\F}_m=\{\wf(x,y)=f(x)\1{y=1}+(1-f(x))\1{y=0} : f \in \F_m
\}
\]
\end{definition}

\begin{lemma}\label{l:bracket}
The $\delta$-upper bracketing entropy of the class $\widetilde{\F}_m$
is no more than $(2m-1)\log(\lceil m/\delta \rceil)$.
\end{lemma}

\begin{proof}
Fix $a,b$ positive integers. Partition $[0,1]$ into $a$ equal intervals
$I_1, \dots, I_a$. Consider the class $H_{m;a,b}$ of all functions
$\wf^U$ that can be formed in the following way: Choose some $m-1$ of
these intervals. Let $C$ be the union of the chosen intervals and
let $C_1, \dots, C_k$ be the nonempty subintervals of $[0,1]$ formed
by subtracting $C$ ($k \le m$). Finally choose $B_1, \dots B_k \in
\{1, \dots, b\}$. Construct the function $\wf^U:[0,1] \times \{0,1\}
\mapsto [0,1]$ by:
\[
\wf^U(x,y)=
\begin{cases}
1 & \text{if $x \in C$} \\
B_i/b & \text{if $x \in C_i$ and $y=1$} \\
1-(B_i-1)/b & \text{if $x \in C_i$ and $y=0$}
\end{cases}
\]

It is easy to see that, by appropriate choices, for any $\wf \in
\widetilde{\F}_m$ we can find an $\wf^U \in H_{m; a,b}$ that is greater
than or equal to it globally. Furthermore, the integral of $\wf^U$ is
less than or equal to $1 + 1/b+(m-1)/a$. The size of $H_{m;a,b}$ is
$\le {\binom{a}{m-1}} b^m \le a^{m-1}b^m$. By choosing $a=b=\lceil m
\delta^{-1} \rceil$, we have shown that
$\exp(\cH(\widetilde{\F}_m,\delta)) \le (\lceil m \delta^{-1}\rceil)^{2m-1}$.
\end{proof}

The final lemma establishes a result about the tail of $\Poisson$
random variables.

\begin{lemma}\label{l:poisstail}
If $K \dist \Poisson(\lambda)$ for some $\lambda>0$, then for any $k>\lambda$: 
\[
\P(K \geq k) \leq
e^{-\lambda}\frac{\lambda^k}{k!}\left(\frac{k}{k-\lambda}\right)
\]
And, consequently, for any $0<\beta<1$ there is a $k_0$ sufficiently
large so that for every $k \ge k_0$,
$\P(K \ge k) \le \exp(-\beta k \log(k))$.
\end{lemma}
\begin{proof}
\begin{align*}
\P(K \ge k) &= e^{-\lambda} \sum_{i=k}^\infty \frac{\lambda^i}{i!} \\
& = e^{-\lambda} \frac{\lambda^k}{k!} \sum_{i=k}^\infty
  {\left(\frac{\lambda}{k}\right)}^{i-k} \left[ \frac{k!k^{i-k}}{i!}
    \right]
\end{align*}
Now, the bracketed term is no more than 1 and we are left
  with a geometric series whose factor, $\lambda/k$ is less
  than $1$ by our assumption so that:
\begin{align*}
\P(K \ge k) & \le e^{-\lambda} \frac{\lambda^k}{k!} \left( 1-\frac{\lambda}{k}
\right)^{-1} \\
& = e^{-\lambda} \frac{\lambda^k}{k!} \left( \frac{k}{k-\lambda}
\right)
\end{align*}

So that for $k$ sufficiently large, from Stirling's formula:
\begin{align*}
\log(\P(K \ge k)) &\le 
 -\log(k!)-\lambda+k\log(\lambda)+\log\left(\frac{k}{k-\lambda}\right) \\
&\le -\beta k \log(k)
\end{align*}
\end{proof}

The main result can now be established.
\begin{thm}\label{t:mine}
Suppose that the prior $\pi$, described in \autoref{s:priorspec} is
based on the hierarchy prior $\kappa$. Suppose that $\kappa$ gives positive
probability to every natural number and that its tail satisfies:
$\kappa(\{k \in \bbN ~:~ k \ge j \}) \le \exp(-\beta j \log(j))$ for
all $j$ sufficiently large and some $\beta>0$. Let $f_0$ be an
arbitrary measurable function from $[0,1]$ into $[0,1]$. Suppose that
$Z_1, Z_2, \dots$ are drawn $\iid$ from the distribution $F_0$, which
has density $\wf_0$ with respect to $\nu$. Equivalently, suppose that the
$Z_i$'s ($Z_i=(X_i, Y_i)$) are drawn as follows: the $X_i$'s are drawn
independently and uniformly from $[0,1]$; and, conditional on $X_i=x_i$,
$Y_i$ is an independent $\Bernoulli(f_0(x_i))$ random variable. Then,
in the notation of \autoref{s:notation}, for any $\epsilon>0$,
$\pi(H_\epsilon|z_1, \dots, z_n) \tendsto 1$ as $n \tendsto \infty$
[${F_0}^\infty$-a.s.]. Specifically, for any $1 \leq p < \infty$, and
any $\epsilon>0$, the posterior mass on the set $\{\wf \in \wcF ~:~
||f-f_0||_p<\epsilon \}$ tends to $1$ as $n \tendsto \infty$
[${F_0}^\infty$-a.s.].
\end{thm}
\begin{proof}
Let $m(n) \defined \lfloor \alpha n / \log(n) \rfloor$. Choose the sequence
$\{\cA_n\}$ as $\cA_n=\wcF_{m(n)}$. Then, by Lemma~\ref{l:bracket},  $\cH(\cA_n, \epsilon) \le
(2m(n)-1)\log(\lceil m(n)/\epsilon \rceil)$. Observe that $m(n)
\tendsup \infty$ as $n\tendsto\infty$. For $n$ large enough, then:
\begin{align*}
\cH(\cA_n, \epsilon) 
&\le 2 m(n)\log(\alpha n/\log(n) \epsilon^{-1}+1) \\
&\le 2 \alpha n/\log(n) [\log(n/\log(n)) +\log(\alpha)-\log(\epsilon)+1] \\
&\le 3 \alpha n \frac{\log(n/\log(n))}{\log(n)}
\end{align*}
Choosing $\alpha=\epsilon/3$, we have proven that 
$\limsup \cH(\cA_n, \epsilon)/n \le \epsilon$.

Now calculate that for all sufficiently large $n$, 
\begin{align*}
-\log(\pi({\cA_n}^c)) 
&= -\log(\kappa(\{k \in \bbN ~:~ k \ge m(n)+1 \})) \\
&\ge \beta [m(n)+1]\log(m(n)+1) \\
&\ge \beta [\alpha n/\log(n)]\log( \alpha n/\log(n) ) \\
&  = \alpha \beta n \frac{\log(n)-\log(\log(n))+\log(\alpha)}{\log(n)} \\
&\ge \frac{1}{2} \alpha \beta n
\end{align*}

Consequently, $\limsup \log(\pi({\cA_n}^c)) /n \le -\frac{1}{2}
\alpha\beta <0$.

Finally, recall that Lemma~\ref{l:pack}
shows that for any $\epsilon>0$, $\pi(K_\epsilon)>0$. Apply
Theorem~\ref{t:main} to complete the proof.

\newcommand{\dtv}{d_{TV}}

The statement about $|| \cdot ||_p$ follows from the equivalence of
the $\cL^p$ and Hellinger metrics for bounded densities. This fact and
other useful inequalities about Hellinger distance $d$ (aka Jeffrey's
distance) are reviewed in \cite[section 5.8 and excercise
  5.7]{devroye:lugosi} which states:
\begin{align}
d(f,g)^2 \le \norm{f-g}_1 \le 2 d(f,g)
\end{align}
Finally, the equivalence of the $\cL^1$-norm and the $\cL^p$ norm for
$1 \le p < \infty$ and for all densities uniformly bounded by a
constant, is well known.
\end{proof}
\comment{
 [P.S. Help! I don't know what
  measure structure to assign to $\pi$. It seems natural to use the
  $\sigma$-algebra generated by the rectangles $C$ that I explain just
  before equation 5.1, but the BSW theorem uses the Hellinger metric Borel
  field; are these the same??? Do I have to prove it? Is it even
  obvious that the $K_\epsilon$ sets are measurable under the
  Hellinger sigma field?]
}
\begin{remark}
Applying Lemma~\ref{l:poisstail} shows that the preceding theorem
applies to any hierarchy prior $\kappa$ whose tail behaves like that
of a $\Poisson(\lambda)$, for any $\lambda>0$. The theorem does not
apply to the case in which $\kappa$ is $\Geometric$.
\end{remark}
\begin{remark}
The restrictions on the tail of the prior occur because condition 3 in
Theorem~\ref{t:main} requires that the ``sieve'' sets $\cA_n$ do
not grow in ``size'' too quickly as $n$ grows. The choice
$\cA=\wcF_{m(n)}$ made in the proof of this theorem is (essentially)
the fastest rate of growth that this situation permits (in the absence
of better bracketing estimates). Accordingly, condition 2 of
Theorem~\ref{t:main} requires that the tail of the prior drops off
somewhat faster than a $\Geometric$ distribution.
\end{remark}
\begin{remark}
For further discussion of how to interpret this result, please see the
discussion in \autoref{c:discuss}.
\end{remark}


%% file: discuss.tex
\chapter{Discussion of Consistency Results}\label{s:priordisc}\label{c:discuss}

It is challenging to suggest what the practical consequences (if any)
of Theorem~\autoref{t:mine} are. It proves that under the modeling
set up with \iid\ observations that has been considered throughout
this thesis, that the posterior of the random split prior $\pi$ is a
consistent estimate of any measurable true regression function {\em if}
the tail of the prior decays at least as fast as $\exp(-\beta j
\log(j))$ for some $\beta>0$. Perhaps it would be wise not to over interpret
the result. After all, it only supplies a sufficient condition for consistency
and does not either establish that $\Geometric$ priors
lead to inconsistent estimators or that $\Poisson$-like priors are a
good (i.e. practical) idea. Additionally, it attempts to prove
consistency in a fairly strong sense: that the posterior mass
concentrates on Hellinger neighborhoods of the truth. Weaker
consistency results, say that the posterior mean be $\cL^2$
consistent, might go through under milder assumptions.

\section{Consideration of the Diaconis and Freedman Results}
Judging from the results of Diaconis and Freedman for their prior, as
discussed in \autoref{s:DFprior}, it might, in fact, be the case that
the posterior is consistent for any hierarchy prior $\kappa$ (with
full support), so long as $f_0$ is not the constant function
$\frac{1}{2}$. That is, it might be that the only situation in which
$\pi$ is inconsistent (even using a ``poor'' choice of $\kappa$) is
when there is no {\em real} pattern at all in the data because every
coin flip was fair. Note that if this were our only concern, namely
that the estimates might be inconsistent under some specific finite
collection of possible scenarios, this could be easily addressed
(albeit in a decidedly non-Bayesian way), by choosing a prior that
puts point mass on the troublesome cases. Consistency for these
exceptions would then be guaranteed by a suitable application of
Doob's result~\cite{doob}. Interestingly, it could destroy consistency
for other cases. An example of such a mixture is given by Diaconis and
Freedman~\cite{df:1984:discrete}. Of course, a true subjective
Bayesian would never change his or her prior in this way. Rather,
unless these cases actually are a subjective impossibility, a purist
would merely see this as an explication of why their {\em true} prior
(that they are perhaps still in the process of articulating) differs
from the former one.

\section{An Experiment to Check a Worrisome Case}
\begin{figure}
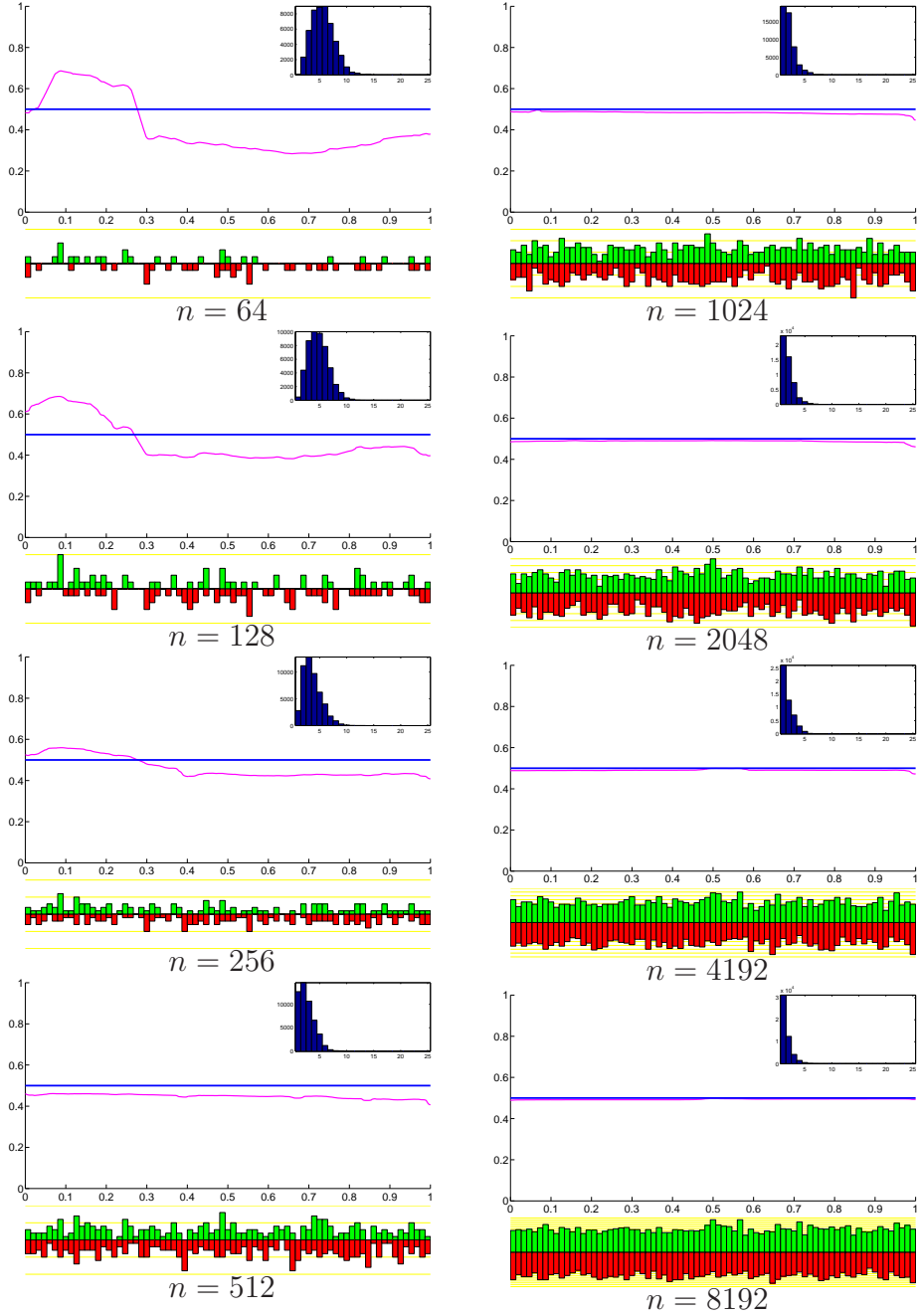

\centering
\begin{minipage}[t]{5.5in}
\begin{minipage}[t]{2.5in}
\centering
\minifigx{unifpost_playpoissnull_64}{-90}{$n=64$}\\
\minifigx{unifpost_playpoissnull_128}{-90}{$n=128$}\\
\minifigx{unifpost_playpoissnull_256}{-90}{$n=256$}\\
\minifigx{unifpost_playpoissnull_512}{-90}{$n=512$}
\end{minipage}
\begin{minipage}[t]{2.5in}
\centering
\minifigx{unifpost_playpoissnull_1024}{-90}{$n=1024$}\\
\minifigx{unifpost_playpoissnull_2048}{-90}{$n=2048$}\\
\minifigx{unifpost_playpoissnull_4096}{-90}{$n=4192$}\\
\minifigx{unifpost_playpoissnull_8192}{-90}{$n=8192$}
\end{minipage}
\end{minipage}
\caption[Null Case with $\Poisson(5)$ Prior]
{Null Case with $\Poisson(5)$ Prior}
\label{f:playpoissnull}
\end{figure}

\begin{figure}
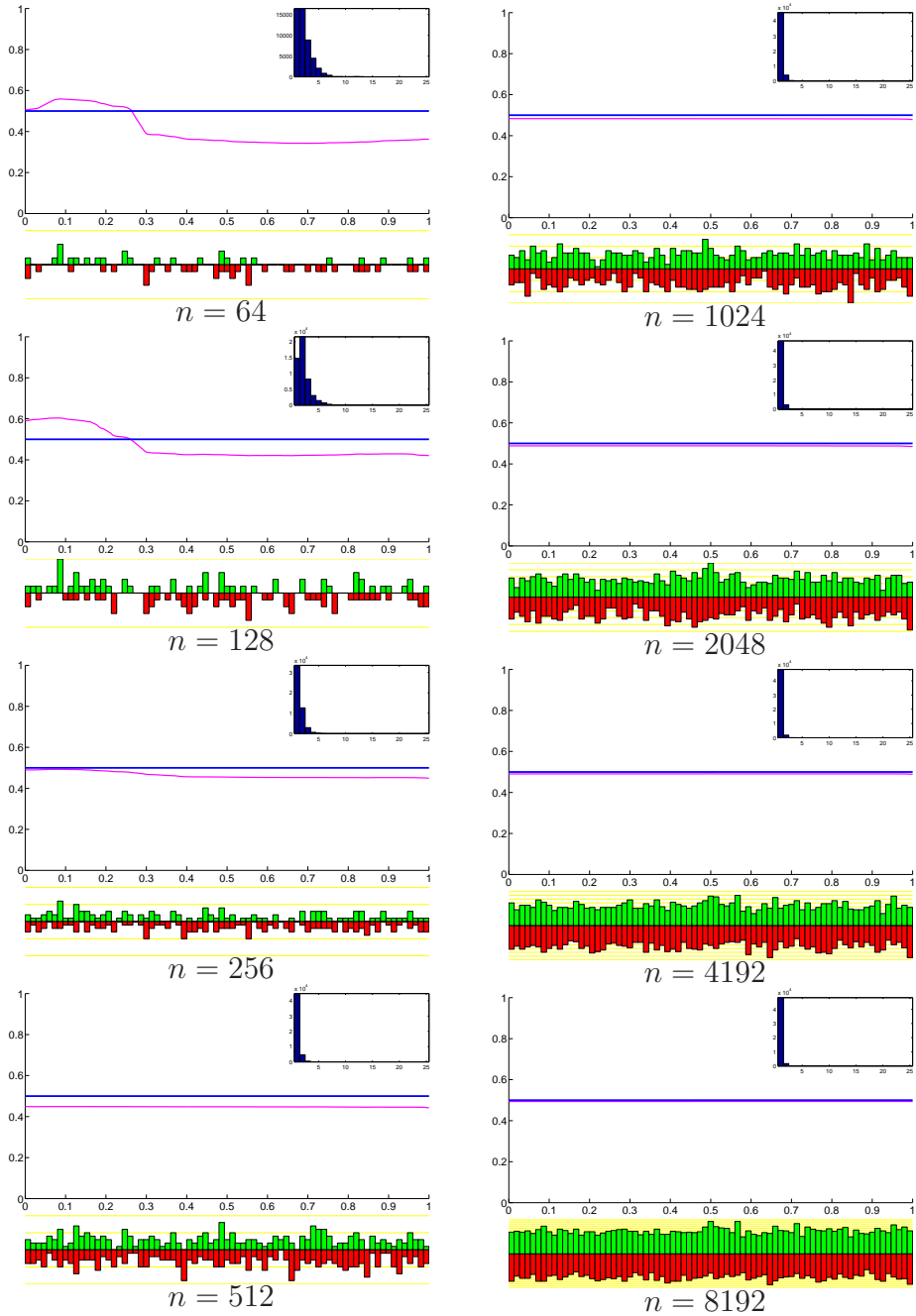

\centering
\begin{minipage}[t]{5.5in}
\begin{minipage}[t]{2.5in}
\centering
\minifigx{unifpost_playgeomnull_64}{-90}{$n=64$}\\
\minifigx{unifpost_playgeomnull_128}{-90}{$n=128$}\\
\minifigx{unifpost_playgeomnull_256}{-90}{$n=256$}\\
\minifigx{unifpost_playgeomnull_512}{-90}{$n=512$}
\end{minipage}
\begin{minipage}[t]{2.5in}
\centering
\minifigx{unifpost_playgeomnull_1024}{-90}{$n=1024$}\\
\minifigx{unifpost_playgeomnull_2048}{-90}{$n=2048$}\\
\minifigx{unifpost_playgeomnull_4096}{-90}{$n=4192$}\\
\minifigx{unifpost_playgeomnull_8192}{-90}{$n=8192$}
\end{minipage}
\end{minipage}
\caption[Null Case with $\Geometric(\frac{1}{2})$ prior]
{Null Case with $\Geometric(\frac{1}{2})$ prior}
\label{f:playgeomnull}
\end{figure}

It seems reasonable to try and test for the possibility of
inconsistency when the true function $f_0 \equiv \half$ by running a simulation
experiment. To do this, generate an increasing sequence of data sets
over a range of sizes that are drawn from the $f_0 \equiv \half$
distribution. That is, since there is no true signal at all in the
data, it will be interesting to see how the posterior mean
responds. For a $\Poisson$ prior, the posterior mean should (at least
eventually) settle down, but will this hold for a $\Geometric$ prior?
Indeed, the results in figures~\ref{f:playpoissnull}, and
\ref{f:playgeomnull} indicate that both estimates correctly identify
the null case, whether using a $\Poisson(5)$ prior or a
$\Geometric(\half)$ one respectively. A new feature of these figures
is the small histogram included on each one. It is a histogram of the
posterior on the number of steps $K$ in the unknown function. The
$\Poisson(5)$ prior starts out believing that there will be a good
number of steps in the data. This is clearly reflected in the
histograms in \autoref{f:playpoissnull} for $n=64$ and $n=128$. The
posterior mean in these cases has more ``wiggle'' than for the
$\Geometric(\half)$ prior even though both estimates see the same
data. For $n=512$ both estimates find that the posterior mean is
roughly constant, but somewhat lower than $\half$. Apparently, this is
a real feature of the data and not an artifact of either
prior. Eventually, for large $n$ all of these minor considerations
wash out and clear preference for very flat models has
triumphed. Interestingly, from $n=512$ or so onwards, looking at the
histogram of $K$, the $\Geometric$ prior has become convinced that the
model has only one split. The $\Poisson$ prior is only beginning to
reach this level of certainty about the truth as $n$ reaches $8192$.

\section{Theory versus Practice}
Nevertheless, if the user insists upon using an estimation procedure that has
been proven to be consistent (say, by the preceding theorem), he or she still
have a great deal of freedom. Roughly speaking, they can use any
hierarchy prior they like on the first $100!!!$ natural numbers and
append to it an arbitrary  $\Poisson$ tail. It is hard to
imagine that there would be any practical difference between the estimators
resulting from this ``extended'' prior and those resulting from the
unmodified one, at least for realistic sample sizes. There certainly
would not be any difference in practice, because in practice we only
approximately compute the posterior anyway and no ordinary Markov chain Monte
Carlo would ever be run long enough to notice the change.

\section{Heuristics about Poisson and Geometric Priors}
Despite these concerns, I still (tentatively) advocate the use of a
$\Geometric$ prior because of the following heuristic arguments: (a) it
favors simple models (b) its tail does not drop off ``too quickly'' so
that it will hopefully not require enormous amounts of evidence for
the data to ``overwhelm'' the prior: by ruling out simple models in
favor of better fitting models (c) its tail does drop steadily;
hopefully this will protect us from ``over-fitting'' (d) if we consider
the mode of the posterior (on a log scale), the geometric prior
penalizes each additional split by a constant ($\log(\alpha)$) so that
for the mode to shift to a model with an additional split we require a
commensurate improvement in the log-likelihood of the data.

If using a $\Poisson(\lambda)$ prior, I would be concerned that I
might have specified a parameter $\lambda$ that was too small. If the
true regression function were unexpectedly complicated, it might take
a large amount of data to overwhelm the prior. Interestingly, if we
attempt to remedy this by taking an exponential mixture of
$\Poisson$'s, we get back a $\Geometric$ prior.

Another objection of mine is that (for modestly large $\lambda$) the
$\Poisson$ prior puts less mass on models with one region than on
models with two regions. This makes sense if I actually expect that
the regression function will be fairly complicated, but it violates my
(frequentist) training to consider a complex model before
``eliminating'' the simpler one.

As a compromise I would propose a prior on $K$ that was $\Geometric$ until a
certain point $k_0$ and then decays like a $\Poisson$. On the other
hand, if, in fact,
$\Geometric$ priors prove reliable or even conservative, it might make
sense to consider an even heavier tailed distribution. A modest
proposal is to take a uniform mixture of $\Geometric(p)$ priors for
$p$ ranging from $0$ to $1$. This results in a prior whose tail decays like
$1/k$, the mass at $K=k$ being $[k(k+1)]^{-1}$.

\section{Conclusions}
None of these arguments is conclusive. In the absence of sound
theoretical arguments it is perhaps best to rely on experimental
evidence. From \autoref{c:examples} there is a good bit of evidence
that $\Geometric$ priors perform well. What do $\Poisson$ priors do
when applied to these data sets? In
\autoref{f:unifpost_playpoisson_all} the posterior mean resulting from
a $\Poisson$ prior for three values of $\lambda$ is plotted. The
results are comparable to what happened as $\alpha$ was varied in
\autoref{f:unifpost_plottalk}. The prior with the shorter tail
flattens out some bumps that the prior with the longer tail leaves in.

\mediumfig{unifpost_playpoisson_all}{-90}{Poisson Posterior
  Example}{{\em The Result of Using a Poisson Prior:} Consider
  again the data from experimental run 1, but apply a
  $\Poisson(\lambda)$ prior with $\lambda=1,5,$ or $10$\\
{\centering
Key:
 \blue{True $f$},
 \dcyan{$\Poisson(1)$ prior},
 \black{$\Poisson(5)$ prior},
 \green{$\Poisson(10)$ prior}
}
}

Still, for some applications especially, the practical question
remains: how to choose $\alpha$ (respectively $\lambda$)? Experience
in the problem domain is the only method I can readily
propose. Alternatives, like using cross-validation or an empirical
Bayes approach remain attractive, but unproven.


%% file: voronoi.tex
\chapter{Extensions}
\label{c:voronoi}
There are numerous ways in which to extend $\pi$, but perhaps the most
pressing is to extend $\pi$ to multi-dimensional data sets $D_n =\{
(X_i, B_i) \}_{i=1}^n$ where the predictor $X_i$ is in $\mathbf{R}^d$.
One route to extend $\pi$ to higher dimensional problems is to observe
that basically, all we need to consider is a suitable way to partition
the space randomly. If an interesting way to choose a partition at
random is found, then describe a new prior by saying: draw a partition
and give each region an independent uniform success probability. One
natural way to randomly partition $\mathbf{R}^d$ is to suppose that a
certain number of generating points are drawn from a Poisson process
with constant rate function $\lambda$ and to associate each point with
its Voronoi (nearest neighbor) region. Alternatively, one could select
a subset of the observed $x$-values at random and use their locations
to determine a Voronoi partition. This alternative is, unfortunately,
not a purely Bayesian proposal since the partitioning depends on the
data set given. On the other hand, it only depends on the $x$-values
of the data set, so that it remains a Bayesian procedure with respect
to the response data (the $y$-values).

To be specific, the prior I consider (call it $\pi^*$) can be
described by the following. Let $x_1, \dots, x_n$ denote the $n$
observed values of
the covariates in $\cX$ and let $\rho$ be a metric on
$\cX$. Proper choice of $\rho$ is essential to good performance in
complex applications, but using Euclidean distance should suffice for simple
problems. To any subset $x_{i_1}, \dots, x_{i_k}$ of the full list
associate the Voronoi partition of $\cX$. That is, say that a point
$x$ is in the $x_{i_j}$ cell if $\rho(x,x_{i_j}) \leq
\rho(x,x_{i_j'})$ for $j' \in \{1, \dots, k\}$. For definiteness, in
the case of ties say that $x$ is in the cell occurring first in the
original ordering of the $x$'s. Consequently, every point $x \in \cX$
is in exactly one cell.  Put a prior on these partitions of the space
$\cX$ by putting a prior on the (finite) set of all possible
(nonempty) subsets of the list $x_1, \dots, x_n$. Say that the prior
probability of a subset only depends on the size of the subset and
that the probability of a given size $k$, is proportional to the
probability that a $\Geometric(1-\alpha)$ random variable takes the
value $k$. Finally, having chosen a subset at random, and consequently
having fixed a partition of $\cX$, assign to each element of the
partition a success probability $s_i$ drawn
uniformly at random from $[0,1]$. The generates a multi-dimensional
regression function $f:\cX \mapsto [0,1]$ at random.

\mediumfig{talk_drawdata}{-90}{A Two Dimensional Data Set and Target
  Function}{{\em A Two Dimensional Data Set and Target Function:} The gray-scale
  on this figure depicts a certain regression function $f_0^*$ on a
  rectangle. $250$ data points are drawn by flipping an
  $f_0^*(x)$-coin at the points indicated. \red{Heads} tend to result
  when $f_0^*$ is large (white). \blue{Tails} tend to result when
  $f_0^*$ is small (black)}

To explore these ideas I simulated a two dimensional data set with 250
data points, illustrated in \autoref{f:talk_drawdata}. The $x$-values
were chosen randomly (albeit not uniformly) from the illustrated
rectangle; the $y$-values were drawn as independent $\Bernoulli$
random variables whose success probability is indicated by the
gray-scale in the figure. Points in the whiter regions have a higher
chance of being an ``x,'' (i.e. $y=1$) while points in the darker
regions have a higher chance of being an ``o'' (i.e. $y=0$). Jointly,
the $x$ and $y$ data was actually generated in my simulation in the
reverse manner: a fair coin was flipped to determine if $y$ will be
$1$ or $0$ and then (conditionally) an $x$-value was chosen. Suppose
$y$ came up as $1$, then with probability $1/2$, $x$ will be drawn
uniformly from the square on the right; otherwise, with probability
$1/2$ it will be drawn from a bivariate normal distribution with standard deviation
$0.1$ that is centered on the left-hand square. If $y$ came up
a $0$, the situation for $x$ would be reversed. These two descriptions
are essentially equivalent and the goal for the posterior mean
estimate is always the same: to estimate the conditional probability
of $y$ to be $1$ given $x$; i.e. to estimate the gray-scale image.

Sampling from the posterior of this prior is (theoretically at least)
quite simple. As before in \autoref{c:implement}, the success
probabilities can be integrated out analytically so that the posterior
probability of a particular (nonempty) subset of size $k$ is proportional to:

\begin{align}
\alpha^{k-1}\prod_{j=1}^k \frac{N^1_j! ~ N^0_j!}{(N^1_j+N^0_j+1)!}
\end{align}

Where $N^1_j$ and $N^0_j$ denote the number of $y$-values equal to $1$
or $0$, respectively, on partition element $j$. Conditionally, the
posterior distribution of a certain success probability $\S_j$ is
$\Beta(N^1_j+1,N^0_j+1)$. Since the number of possible subsets is
finite and all of them (except the empty subset) have positive
posterior probability, a standard Metropolis-Hastings type MCMC allows
us to sample from the posterior (at least in theory).  In practice,
though, the rate of mixing matters; in an effort to improve this
I have conducted preliminary work that employs the simulated tempering technique
developed by Geyer and Thompson~\cite{geyer:thompson:1993}.

All that remains to be specified is a transitive
random-walk on (nonempty) subsets of a set of $n$ elements which has a
known stationary distribution. This is easily done. Identify the class
of all subsets of a set of size $n$ with the class of binary vector of
length $n$, with each coordinate indicating the presence or absence of
a given element. Exclude the $0$-vector from this set. Consider the
random walk that picks a number $J$ randomly from $1$ to $n$ and then
proposes flipping the $J$'th bit. The proposal is not allowed if it
would create the $0$-vector; hold in this case. This Markov chain is
easily seen to sample uniformly from the class of all non-empty
subsets and consequently it is easy to modify it with the
Metropolis-Hastings ratio in order to sample from the posterior.

\begin{figure}
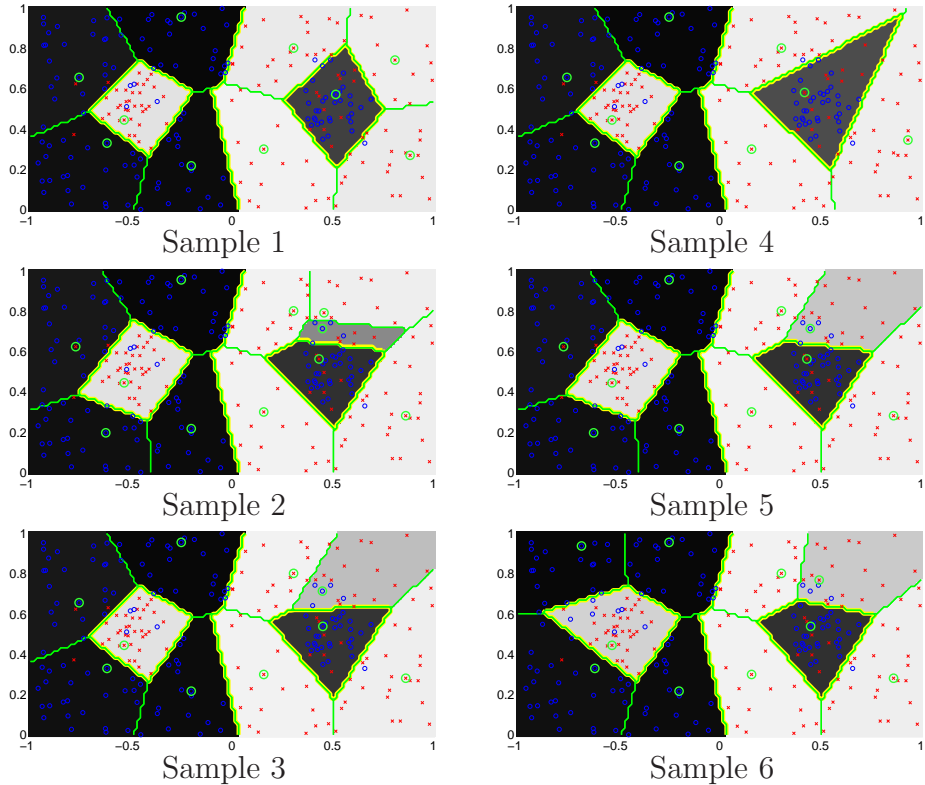

\centering
\begin{minipage}[t]{5.5in}
\begin{minipage}[t]{2.5in}
\centering
\minifigx{show_1}{-90}{Sample 1}\\
\minifigx{show_2}{-90}{Sample 2}\\
\minifigx{show_3}{-90}{Sample 3}
\end{minipage}
\begin{minipage}[t]{2.5in}
\centering
\minifigx{show_4}{-90}{Sample 4}\\
\minifigx{show_5}{-90}{Sample 5}\\
\minifigx{show_6}{-90}{Sample 6}
\end{minipage}
\end{minipage}
\caption[Modal Samples: Voronoi Posterior]{{\em Modal Samples: Voronoi
    Posterior} The pictured samples were the ones most frequently
    occurring in a sample drawn from the Voronoi posterior
    $\pi^*$. Each sample is equivalent to a certain subset of the
    original covariate list. The included points determine the
    partition and are drawn with a green circle. The gray-scale
    on a given partition element represents the posterior mean of the
    corresponding success probability parameter. Notice how
    parsimoniously the circled points determine regions that isolate
    out the two clumps of data.}
\label{f:modalsamples}
\end{figure}

\comment{
substantially and applies to any partition-based prior. To improve the
mixing rate of the MCMC, I employ a simulated tempering algorithm. An
example illustrating the posterior mean under $\pi''$ is given in
figure~\ref{f:temper2_100K}.
}

\mediumfig{temper2_100K}{-90}{A Bivariate Posterior Mean Example}{\emph{A Bivariate Example.} $250$ red and blue markers
  were put down randomly as shown. The posterior-mean estimate of $f_0$
  under $\pi^*$, i.e. the estimated conditional probability of red at
  each position, is shown in gray; notice how the modal samples from
  \autoref{f:modalsamples} are incorporated into the posterior.}

\mediumfig{vnw1_200K_250_05}{-90}{Weighted Voronoi Posterior}{{\em Weighted
  Voronoi Posterior:} Redefining the Voronoi cells to use random
  weights allows for elliptical arcs and lines to be used in the
  partition and eliminates some artifacts.}

I have also proposed extensions of this technique to put a prior on
smooth functions; under this proposal the prior concentrates on
functions which softly partition the space using weighted Voronoi
regions, but provide for smooth transitions between regions.

A simple extension of this technique that removes many of the
artifacts that are otherwise present because of the dependence on the
particular locations of the covariates can be made by using a randomly
weighted Voronoi partitions~\cite{voronoibook}. To make such a prior,
simply augment each cell with a random variable $W_j$ which is an
independent $\Gamma(\gamma, \gamma^{-1})$ \apriori. Then, redefine the
Voronoi partition so that cells with higher weight tend to be
bigger. Specifically, consider the point $x$ to be in cell $j$ rather
than cell $k$ whenever: $\rho(x,x_j)/w_j \leq
\rho(x,x_k)/w_k$. Intuitively, if an ordinary Voronoi partition can be
understood by supposing that ``crystals'' grow out radially from the
generating seeds until they hit other growing crystals, then this
weighted Voronoi partition allows for the crystals to grow at
different rates. It is also simple to allow the crystals to start
growing at different times (by simply adding a different random offset
to $\rho$ for each cell). The overall scale of the prior on $W$ is irrelevant
so this parameterization sets the mean to $1$. For the experiment
shown in \autoref{f:vnw1_200K_250_05} I used $\gamma=5$. As can be
seen, this simple modification allows the posterior to choose neat
balls and lines with which to isolate out the different contours of
the data. The posterior is computed using standard MCMC techniques.

Another route to extend $\pi$ is to utilize the observed connection
between $\pi$ and $\piDF$ (for details, please refer to
\autoref{s:DFprior}). To make that connection, I employed binary
random variables $\eta_i(u)$ that indicated if a test point $u$ was or
was not above a certain random threshold $V_i$. To generalize, then,
we can take $\eta_i(x)$ to indicate if $x$ is in a certain random
half-space $H_i$. This would be similar to a version of CART in which
we split first into two halves via $H_1$, and then split each half via
$H_2$, and so on. We could make something close to ordinary CART if we
utilized random coordinate aligned half-spaces and if we employed a
suitably ``regularized'' $\piDF$ which did not assign an independent
uniform to all possible $2^k$ binary $k$-tuples at level $k$. More
generally, one can invent other ways of ``regularizing'' $\piDF$ so that
numerous binary tuples are tied together.

\mediumfig{c_bag2d_load}{-90}{Bagged CART in 2d}{{\em Bagged CART in
    2d:} Running bagged CART on this two dimensional data yields a
    reasonably good estimate with an interesting horizontal and vertical blurring
    pattern}
Finally, for comparison, \autoref{f:c_bag2d_load} the result of
bagging ordinary two-dimensional CART on this same data is shown. It
has a clear advantage on the vertical split down the center that
happens to be coordinate aligned. It does a decent job of isolating
the two clumps of data. Interestingly, there is a clear horizontal and vertical
blurring pattern that arises from the use of partitions that only
partially isolated the clumped data: it is isolated in one coordinate,
but not in the other.


%% file: afterword.tex
\chapter{Afterword}\label{c:afterword}

As statisticians, we analyze data, formulate models, estimate
parameters, and use these models to form predictions about future
data. Broadly speaking, then, our business is inference. We go to a
lot of trouble to formulate good models and we have spent a great deal
of effort debating about the details of how to make inference within a
model -- e.g. Bayesian versus frequentist inference.

Generally, though, the way we select our model remains an art.
``Non-parametric'' methods are a step forward here, because, generally
speaking, they at least prescribe how to select a variety of
``smoothing-parameters'' which essentially determine which model among
some class of models we actually apply. The main topic of this thesis
is of this sort: A Bayesian approach to the question of how to
estimate the number and location of change-points or, more generally, how
to choose a partition of the data into approximately exchangeable subsets.

In this afterword, I would like to step back from the details of this
subject and address the more general problem of how we
\emph{formulate} a model. Sometimes, we prefer to shunt responsibility
for this and appeal to the scientist for help; but, fundamentally, all
inference comes back to data eventually -- how else did the scientists
discover their model? Considering, then, that the formulation of a
model is essentially ``an art,'' and that we know that the results of
our analyses are not absolute, but \emph{relative} to the modeling
choices we have made, how can we be so bold as to expect that reality
will conform to our model-specific ``confidence intervals?'' This is
not to say that statistical practice does not work, or that the formal
properties of our analyses under our stated assumptions are invalid,
but simply to remind the reader that the thread that connects the
prescribed model with reality has no \emph{formal} basis.

In practice, a good rule of thumb is to consider several models and,
consciously or not, select one which is simple and which we expect will
fit the data at least nearly as well as other more complicated models
that we might prefer not to have to consider. In practice, we look at
the data \emph{ourselves} before selecting a model and build in any
particular types of regularity that we happen to \emph{notice} that it
possesses into our model -- at least if we think it will affect our
conclusions.

Bayesian analysts do not escape these problems. They may subjectively
allow for a mixture of several different models and for a range of
``hyper-parameters'' but this only ameliorates the core problem because
no-one ever accounts for all the possible regularities that might be found.

To further make my point, consider the following, admittedly fanciful,
thought experiment. An alien race comes to earth and challenges
mankind to an intelligence test. We are given a binary time-series
to analyze one-hundred bits at a time. It begins innocuously enough:
\begin{verbatim}
     01000001000000100000110111101110000000100101000011
     10101010000000001110010010000000001011001101101010
\end{verbatim}

What statistical models shall we consider? We could try $\iid$
Bernoulli, or perhaps a hidden Markov model. If ambitious, we might
let the data choose the order of our model. Upon doing so, we find
that the null model fits best, despite the superficial appearance of
runs, and we model the data as random coins with success probability
$0.35$. The aliens ask us to give a confidence interval on the
fraction of $1$'s in the next 100 bits and to estimate the probability
that the 50'th bit of this new data will be 1. Most statisticians fall
back on the classics and use something approximately like $0.35 \pm
0.1$ for both answers. Others who used a richer model suggest wider
intervals and the ordinary debate ensues.

So far so good, right? Or should we worry that our glance at the data
and default choice of classical model may miss structure that we
failed to notice? Naah! The aliens couldn't be that tricky.... Then
the next 100 bits arrive.

\begin{verbatim}
     00000000000000000000000000000000000000000000000000
     00000000000000000000000000000000000000000000000000
\end{verbatim}

A dramatic failure, but no problem, the advocates of the HMM model
were ready for this sort of thing. We have, they argue, simply
encountered a hidden state which always produces 0. This explains the
data well enough, but some remain skeptical. They propose that the
character of the data may change with every new segment of 100 bits so
that the effective size of our data set is only 2. The next 100 bits arrive.

\begin{verbatim}
     00011100011100011100011100011100011100011100011100
     01110001110001110001110001110001110001110001110001
\end{verbatim}

And again we are surprised. Some HMM advocates insist we just need to
add a number of new hidden states. Others extend the HMM model to
favor this sort of cycle-like behavior. The next 100 bits arrive.

\begin{verbatim}
     00100100001111110110101010001000100001011010001100
     00100011010011000100110001100110001010001011100000
\end{verbatim}

And most are satisfied that the data has returned to $\iid$, but with
the world attention that this situation has generated someone notices
that these are, in fact, the first 100 bits in the expansion for the
fractional part of $\pi$. (C.f. \texttt{http://www.algonet.se/\breakpt{\textasciitilde}eliasb/pi/binpi.html}) 
Oh dear -- we certainly hadn't planned on this -- but we ignore this
regularity in the data at our peril. It couldn't happen by chance,
could it?

As the test continues, we continue to be surprised by the patterns we
are sent. By now, we have learned a lot: we know less than we think
about the future. Every 100 bits we have been presented with a pattern
that we hadn't expected.  Finally, we are given sequence after
sequence that we can't explain.  Eventually the aliens conclude that,
however feebly, we are, at least, a modestly intelligent form of life;
and, taking pity on us, they decide to reveal the patterns we hadn't
discovered. The last sequences, for example, were actually
Shakespeare, encoded by simple alien cipher that no human had ever
considered. Furthermore, the first sequence wasn't actually
``random:'' to generate it, all we had to do was start \texttt{matlab}
or an equivalent (alien) computational program and type:
\begin{verbatim}
                x=rand(1,100)<=0.397;
                sprintf('%c',x+'0')
\end{verbatim}

Perhaps you object to my example. You prefer regression to time-series
analysis and are content to consider data for which no-one would object
to the model that the responses are independent given the predictors.
Perhaps the problem of extrapolating a non-stationary time-series
seems far too lofty to you. But you haven't escaped it by wishing it
away; in fact, the time-series problem can be embedded in the
regression problem. We need only suppose that the covariates are
tested one by one in some fixed designed fashion so that we see the
responses sequentially. If we know that the regression function ``smoothly''
depends on the covariates, present methods can be expected to work;
but, if the dependence is sufficiently complicated, each new data point
tells us something entirely new, just like in my fictitious time-series.

Even if the data is generated $\iid$ -- the regression case considered
in most of my theoretical work -- there is plenty of room for
improvement.  In this situation, we are, indeed, much better off -- we
can make rigorous probability statements about the quality of our
predictions \emph{on average}. Even so, as more and more data come in,
and the general shape of the regression function becomes more tightly
resolved, the ``knowledge we gain'' itself comes to us in a time-series
fashion. Because of this, we cannot easily make predictions about the
regression function at some \emph{fixed} point.

For example, suppose the unit interval were divided into subintervals
of size $\frac{1}{2}$, $\frac{1}{4}$, $\frac{1}{8}$, etc... and
suppose that the regression function takes a different value on each
piece. In this way, we would quickly have enough data to estimate
large-scale features, like the value of the regression function on the
larger pieces, but if we are asked to make a prediction on one of the
very small pieces, what are we to do? Perhaps, if we were smart
enough, we wouldn't blithely approximate the regression function as
``smooth,'' we would look for patterns in the regression values on the
large pieces to help extrapolate to the smaller pieces.  If, in fact,
the regression values were seen to alternate between high and low, we
would be silly to treat the function as ``smooth.'' Instead, let us
hope that were are lucky enough that it is quite ``regular.''

What, then are we to do? Consider again the alien's test, and the
complex sequence of modeling decisions that we needed to make along
the way. How can we summarize the thought process that we went through
as surprising data continued to come in?  What possible prior on models
could our analysis (even approximately) conform to?  Our only
recourse, it seems, is to formalize the idea of regularity and make
explicit the manner in which we choose a model to accord with the
regularities apparent in the data. 

A reasonable defining property of regularity is that a distribution is
\emph{regular} if it can be (approximately) reproduced within a
certain budget of time by applying some modestly short computer
program to the data that we have previously seen and a ``random
sequence.'' Roughly speaking, then, we can put a prior on models
and/or regularities by putting a suitable prior on computer programs.
Alternatively, we can select among the computer programs in a manner
conforming more closely to the ``method of maximum likelihood.''
Formal versions of these ideas have been
proposed~\cite{donoho_KC},\cite{solomonoff},\cite{chaitin},\cite{kolmogorov_KC2},\cite{kolmogorov_KC1},\cite{cover},\cite{schmidhuber}
though much work is needed to formulate methods that are ready for actual
use. Still, it seems to this author that the further development of
some version of these ideas is an essential, natural, and unavoidable
step in the progression of statistical thinking.
